%
%
\font \aufont=  cmr10 at 14pt 
\font\titfont=  cmr10 at 18pt 
\font\footfont=cmr10 at 8  pt
   

%
%

\def\fmtversion{2.0}
\catcode`\@=11
\ifx\amstexloaded@\relax\catcode`\@=\active
 \endinput\else\let\amstexloaded@\relax\fi
\def\W@{\immediate\write\sixt@@n}
\def\CR@{\W@{}\W@{AmS-TeX - Version \fmtversion}\W@{}
\W@{COPYRIGHT 1985, 1990 - AMERICAN MATHEMATICAL SOCIETY}
\W@{Use of this macro package is not restricted provided}
\W@{each use is acknowledged upon publication.}\W@{}}
\CR@
\everyjob{\CR@}
\toksdef\toks@@=2
\long\def\rightappend@#1\to#2{\toks@{\\{#1}}\toks@@
 =\expandafter{#2}\xdef#2{\the\toks@@\the\toks@}\toks@{}\toks@@{}}
\def\alloclist@{}
\newif\ifalloc@
\def\showallocations{{\def\\{\immediate\write\m@ne}\alloclist@}\alloc@true}
\def\alloc@#1#2#3#4#5{\global\advance\count1#1by\@ne
 \ch@ck#1#4#2\allocationnumber=\count1#1
 \global#3#5=\allocationnumber
 \edef\next@{\string#5=\string#2\the\allocationnumber}%
 \expandafter\rightappend@\next@\to\alloclist@}
\newcount\count@@
\newcount\count@@@
\def\FN@{\futurelet\next}
\def\DN@{\def\next@}
\def\DNii@{\def\nextii@}
\def\RIfM@{\relax\ifmmode}
\def\RIfMIfI@{\relax\ifmmode\ifinner}
\def\setboxz@h{\setbox\z@\hbox}
\def\wdz@{\wd\z@}
\def\boxz@{\box\z@}
\def\setbox@ne{\setbox\@ne}
\def\wd@ne{\wd\@ne}
\def\iterate{\body\expandafter\iterate\else\fi}
\newlinechar=`\^^J
\def\err@#1{\errmessage{AmS-TeX error: #1}}
\newhelp\defaulthelp@{Sorry, I already gave what help I could...^^J
Maybe you should try asking a human?^^J
An error might have occurred before I noticed any problems.^^J
``If all else fails, read the instructions.''}
\def\Err@#1{\errhelp\defaulthelp@\errmessage{AmS-TeX error: #1}}
\def\eat@#1{}
\def\in@#1#2{\def\in@@##1#1##2##3\in@@{\ifx\in@##2\in@false\else\in@true\fi}%
 \in@@#2#1\in@\in@@}
\newif\ifin@
\def\space@.{\futurelet\space@\relax}
\space@. %
\newhelp\athelp@
{Only certain combinations beginning with @ make sense to me.^^J
Perhaps you wanted \string\@\space for a printed @?^^J
I've ignored the character or group after @.}
\def\futureletnextat@{\futurelet\next\at@}
{\catcode`\@=\active
\lccode`\Z=`\@ \lccode`\I=`\I \lowercase{%
\gdef@{\csname futureletnextatZ\endcsname}
\expandafter\gdef\csname atZ\endcsname                                      
 {\ifcat\noexpand\next a\def\next{\csname atZZ\endcsname}\else
 \ifcat\noexpand\next0\def\next{\csname atZZ\endcsname}\else
 \ifcat\noexpand\next\relax\def\next{\csname atZZZ\endcsname}\else
 \def\next{\csname atZZZZ\endcsname}\fi\fi\fi\next}
\expandafter\gdef\csname atZZ\endcsname#1{\expandafter                      
 \ifx\csname #1Zat\endcsname\relax\def\next
 {\errhelp\expandafter=\csname athelpZ\endcsname
 \csname errZ\endcsname{Invalid use of \string@}}\else
 \def\next{\csname #1Zat\endcsname}\fi\next}
\expandafter\gdef\csname atZZZ\endcsname#1{\expandafter                     
 \ifx\csname \string#1ZZat\endcsname\relax\def\next
 {\errhelp\expandafter=\csname athelpZ\endcsname
 \csname errZ\endcsname{Invalid use of \string@}}\else
 \def\next{\csname \string#1ZZat\endcsname}\fi\next}
\expandafter\gdef\csname atZZZZ\endcsname#1{\errhelp                        
 \expandafter=\csname athelpZ\endcsname
 \csname errZ\endcsname{Invalid use of \string@}}}}                         
\def\atdef@#1{\expandafter\def\csname #1@at\endcsname}
\def\atdef@@#1{\expandafter\def\csname \string#1@@at\endcsname}
\newhelp\defahelp@{If you typed \string\define\space cs instead of
\string\define\string\cs\space^^J
I've substituted an inaccessible control sequence so that your^^J
definition will be completed without mixing me up too badly.^^J
If you typed \string\define{\string\cs} the inaccessible control sequence^^J
was defined to be \string\cs, and the rest of your^^J
definition appears as input.}
\newhelp\defbhelp@{I've ignored your definition, because it might^^J
conflict with other uses that are important to me.}
\def\define{\FN@\define@}
\def\define@{\ifcat\noexpand\next\relax
 \expandafter\define@@\else\errhelp\defahelp@                               
 \err@{\string\define\space must be followed by a control
 sequence}\expandafter\def\expandafter\nextii@\fi}                          
\def\undefined@@@@@@@@@@{}
\def\preloaded@@@@@@@@@@{}
\def\next@@@@@@@@@@{}
\def\define@@#1{\ifx#1\relax\errhelp\defbhelp@                              
 \err@{\string#1\space is already defined}\DN@{\DNii@}\else
 \expandafter\ifx\csname\expandafter\eat@\string                            
 #1@@@@@@@@@@\endcsname\undefined@@@@@@@@@@\errhelp\defbhelp@
 \err@{\string#1\space can't be defined}\DN@{\DNii@}\else
 \expandafter\ifx\csname\expandafter\eat@\string#1\endcsname\relax          
 \global\let#1\undefined\DN@{\def#1}\else\errhelp\defbhelp@
 \err@{\string#1\space is already defined}\DN@{\DNii@}\fi
 \fi\fi\next@}
\let\redefine\def
\def\predefine#1#2{\let#1#2}
\def\undefine#1{\let#1\undefined}
\newdimen\captionwidth@
\captionwidth@\hsize
\advance\captionwidth@-1.5in
\def\pagewidth#1{\hsize#1\relax
 \captionwidth@\hsize\advance\captionwidth@-1.5in}
\def\pageheight#1{\vsize#1\relax}
\def\hcorrection#1{\advance\hoffset#1\relax}
\def\vcorrection#1{\advance\voffset#1\relax}

\let\graveaccent\`
\let\acuteaccent\'
\let\tildeaccent\~
\let\hataccent\^
\let\underscore\_
\let\B\=
\let\D\.
\let\ic@\/
\def\/{\unskip\ic@}
\def\textfonti{\the\textfont\@ne}
\def\t#1#2{{\edef\next@{\the\font}\textfonti\accent"7F \next@#1#2}}
\def~{\unskip\nobreak\ \ignorespaces}
\def\.{.\spacefactor\@m}
\atdef@;{\leavevmode\null;}
\atdef@:{\leavevmode\null:}
\atdef@?{\leavevmode\null?}
\def\@{\char64 }
\def\relaxnext@{\let\next\relax}
\atdef@-{\relaxnext@\leavevmode
 \DN@{\ifx\next-\DN@-{\FN@\nextii@}\else
  \DN@{\leavevmode\hbox{-}}\fi\next@}%
 \DNii@{\ifx\next-\DN@-{\leavevmode\hbox{---}}\else
  \DN@{\leavevmode\hbox{--}}\fi\next@}%
 \FN@\next@}
\def\srdr@{\kern.16667em}
\def\drsr@{\kern.02778em}
\def\sldl@{\kern.02778em}
\def\dlsl@{\kern.16667em}
\atdef@"{\unskip\relaxnext@
 \DN@{\ifx\next\space@\DN@. {\FN@\nextii@}\else
  \DN@.{\FN@\nextii@}\fi\next@.}%
 \DNii@{\ifx\next`\DN@`{\FN@\nextiii@}\else
  \ifx\next\lq\DN@\lq{\FN@\nextiii@}\else
  \DN@####1{\FN@\nextiv@}\fi\fi\next@}%
 \def\nextiii@{\ifx\next`\DN@`{\sldl@``}\else\ifx\next\lq
  \DN@\lq{\sldl@``}\else\DN@{\dlsl@`}\fi\fi\next@}%
 \def\nextiv@{\ifx\next'\DN@'{\srdr@''}\else
  \ifx\next\rq\DN@\rq{\srdr@''}\else\DN@{\drsr@'}\fi\fi\next@}%
 \FN@\next@}

\def\textfontii{\the\textfont\tw@}
\def\lbrace@{\delimiter"4266308 }
\def\rbrace@{\delimiter"5267309 }
\def\{{\RIfM@\lbrace@\else{\textfontii f}\spacefactor\@m\fi}
\def\}{\RIfM@\rbrace@\else
 \let\@sf\empty\ifhmode\edef\@sf{\spacefactor\the\spacefactor}\fi
 {\textfontii g}\@sf\relax\fi}
\let\lbrace\{
\let\rbrace\}
\def\AmSTeX{{\textfontii A}\kern-.1667em\lower.5ex\hbox
 {\textfontii M}\kern-.125em{\textfontii S}-\TeX}
\def\vmodeerr@#1{\Err@{\string#1\space not allowed between paragraphs}}
\def\mathmodeerr@#1{\Err@{\string#1\space not allowed in math mode}}
\def\linebreak{\RIfM@\mathmodeerr@\linebreak\else
 \ifhmode\unskip\unkern\break\else\vmodeerr@\linebreak\fi\fi}

\newskip\saveskip@
\def\allowlinebreak{\RIfM@\mathmodeerr@\allowlinebreak\else
 \ifhmode\saveskip@\lastskip\unskip
 \allowbreak\ifdim\saveskip@>\z@\hskip\saveskip@\fi
 \else\vmodeerr@\allowlinebreak\fi\fi}
\def\nolinebreak{\RIfM@\mathmodeerr@\nolinebreak\else
 \ifhmode\saveskip@\lastskip\unskip
 \nobreak\ifdim\saveskip@>\z@\hskip\saveskip@\fi
 \else\vmodeerr@\nolinebreak\fi\fi}
\def\newline{\relaxnext@
 \DN@{\RIfM@\expandafter\mathmodeerr@\expandafter\newline\else
  \ifhmode\ifx\next\par\else
  \expandafter\unskip\expandafter\null\expandafter\hfill\expandafter\break\fi
  \else
  \expandafter\vmodeerr@\expandafter\newline\fi\fi}%
 \FN@\next@}
\def\dmatherr@#1{\Err@{\string#1\space not allowed in display math mode}}
\def\nondmatherr@#1{\Err@{\string#1\space not allowed in non-display math
 mode}}
\def\onlydmatherr@#1{\Err@{\string#1\space allowed only in display math mode}}
\def\nonmatherr@#1{\Err@{\string#1\space allowed only in math mode}}
\def\mathbreak{\RIfMIfI@\break\else
 \dmatherr@\mathbreak\fi\else\nonmatherr@\mathbreak\fi}
\def\nomathbreak{\RIfMIfI@\nobreak\else
 \dmatherr@\nomathbreak\fi\else\nonmatherr@\nomathbreak\fi}
\def\allowmathbreak{\RIfMIfI@\allowbreak\else
 \dmatherr@\allowmathbreak\fi\else\nonmatherr@\allowmathbreak\fi}
\def\pagebreak{\RIfM@
 \ifinner\nondmatherr@\pagebreak\else\postdisplaypenalty-\@M\fi
 \else\ifvmode\removelastskip\break\else\vadjust{\break}\fi\fi}
\def\nopagebreak{\RIfM@
 \ifinner\nondmatherr@\nopagebreak\else\postdisplaypenalty\@M\fi
 \else\ifvmode\nobreak\else\vadjust{\nobreak}\fi\fi}
\def\nonvmodeerr@#1{\Err@{\string#1\space not allowed within a paragraph
 or in math}}
\def\vnonvmode@#1#2{\relaxnext@\DNii@{\ifx\next\par\DN@{#1}\else
 \DN@{#2}\fi\next@}%
 \ifvmode\DN@{#1}\else
 \DN@{\FN@\nextii@}\fi\next@}
\def\newpage{\vnonvmode@{\vfill\break}{\nonvmodeerr@\newpage}}
\def\smallpagebreak{\vnonvmode@\smallbreak{\nonvmodeerr@\smallpagebreak}}
\def\medpagebreak{\vnonvmode@\medbreak{\nonvmodeerr@\medpagebreak}}
\def\bigpagebreak{\vnonvmode@\bigbreak{\nonvmodeerr@\bigpagebreak}}
\def\NoBlackBoxes{\global\overfullrule\z@}
\def\BlackBoxes{\global\overfullrule5\p@}
\def\Invalid@#1{\def#1{\Err@{\Invalid@@\string#1}}}
\def\Invalid@@{Invalid use of }
\Invalid@\caption
\Invalid@\captionwidth
\newdimen\smallcaptionwidth@
\def\topspace{\mid@false\ins@}
\def\midspace{\mid@true\ins@}
\newif\ifmid@
\def\captionfont@{}
\def\ins@#1{\relaxnext@\allowbreak
 \smallcaptionwidth@\captionwidth@\gdef\thespace@{#1}%
 \DN@{\ifx\next\space@\DN@. {\FN@\nextii@}\else
  \DN@.{\FN@\nextii@}\fi\next@.}%
 \DNii@{\ifx\next\caption\DN@\caption{\FN@\nextiii@}%
  \else\let\next@\nextiv@\fi\next@}%
 \def\nextiv@{\vnonvmode@
  {\ifmid@\expandafter\midinsert\else\expandafter\topinsert\fi
   \vbox to\thespace@{}\endinsert}
  {\ifmid@\nonvmodeerr@\midspace\else\nonvmodeerr@\topspace\fi}}%
 \def\nextiii@{\ifx\next\captionwidth\expandafter\nextv@
  \else\expandafter\nextvi@\fi}%
 \def\nextv@\captionwidth##1##2{\smallcaptionwidth@##1\relax\nextvi@{##2}}%
 \def\nextvi@##1{\def\thecaption@{\captionfont@##1}%
  \DN@{\ifx\next\space@\DN@. {\FN@\nextvii@}\else
   \DN@.{\FN@\nextvii@}\fi\next@.}%
  \FN@\next@}%
 \def\nextvii@{\vnonvmode@
  {\ifmid@\expandafter\midinsert\else
  \expandafter\topinsert\fi\vbox to\thespace@{}\nobreak\smallskip
  \setboxz@h{\noindent\ignorespaces\thecaption@\unskip}%
  \ifdim\wdz@>\smallcaptionwidth@\centerline{\vbox{\hsize\smallcaptionwidth@
   \noindent\ignorespaces\thecaption@\unskip}}%
  \else\centerline{\boxz@}\fi\endinsert}
  {\ifmid@\nonvmodeerr@\midspace
  \else\nonvmodeerr@\topspace\fi}}%
 \FN@\next@}
\def\newcodes@{\catcode`\\=12 \catcode`\{=12 \catcode`\}=12 \catcode`\#=12
 \catcode`\%=12\relax}
\def\oldcodes@{\catcode`\\=0 \catcode`\{=1 \catcode`\}=2 \catcode`\#=6
 \catcode`\%=14\relax}
\def\comment{\newcodes@\endlinechar=10 \comment@}
{\lccode`\0=`\\
\lowercase{\gdef\comment@#1^^J{\comment@@#10endcomment\comment@@@}%
\gdef\comment@@#10endcomment{\FN@\comment@@@}%
\gdef\comment@@@#1\comment@@@{\ifx\next\comment@@@\let\next\comment@
 \else\def\next{\oldcodes@\endlinechar=`\^^M\relax}%
 \fi\next}}}
\def\pr@m@s{\ifx'\next\DN@##1{\prim@s}\else\let\next@\egroup\fi\next@}
\def\prime{{\null\prime@\null}}
\mathchardef\prime@="0230
\let\dsize\displaystyle

\let\ssize\scriptstyle

\def\,{\RIfM@\mskip\thinmuskip\relax\else\kern.16667em\fi}
\def\!{\RIfM@\mskip-\thinmuskip\relax\else\kern-.16667em\fi}
\let\thinspace\,
\let\negthinspace\!
\def\medspace{\RIfM@\mskip\medmuskip\relax\else\kern.222222em\fi}
\def\negmedspace{\RIfM@\mskip-\medmuskip\relax\else\kern-.222222em\fi}
\def\thickspace{\RIfM@\mskip\thickmuskip\relax\else\kern.27777em\fi}
\let\;\thickspace
\def\negthickspace{\RIfM@\mskip-\thickmuskip\relax\else
 \kern-.27777em\fi}
\atdef@,{\RIfM@\mskip.1\thinmuskip\else\leavevmode\null,\fi}
\atdef@!{\RIfM@\mskip-.1\thinmuskip\else\leavevmode\null!\fi}
\atdef@.{\RIfM@&&\else\leavevmode.\spacefactor3000 \fi}
\def\and{\DOTSB\;\mathchar"3026 \;}

\def\frac#1#2{{#1\over#2}}

\newdimen\ex@
\ex@.2326ex
\Invalid@\thickness
\def\thickfrac{\relaxnext@
 \DN@{\ifx\next\thickness\let\next@\nextii@\else
 \DN@{\nextii@\thickness1}\fi\next@}%
 \DNii@\thickness##1##2##3{{##2\above##1\ex@##3}}%
 \FN@\next@}

\def\thickfracwithdelims#1#2{\relaxnext@\def\ldelim@{#1}\def\rdelim@{#2}%
 \DN@{\ifx\next\thickness\let\next@\nextii@\else
 \DN@{\nextii@\thickness1}\fi\next@}%
 \DNii@\thickness##1##2##3{{##2\abovewithdelims
 \ldelim@\rdelim@##1\ex@##3}}%
 \FN@\next@}

\def\:{\nobreak\hskip.1111em\mathpunct{}\nonscript\mkern-\thinmuskip{:}\hskip
 .3333emplus.0555em\relax}
\def\snug{\unskip\kern-\mathsurround}
\def\topsmash{\top@true\bot@false\smash@}
\def\botsmash{\top@false\bot@true\smash@}
\newif\iftop@
\newif\ifbot@
\def\smash{\top@true\bot@true\smash@}
\def\smash@{\RIfM@\expandafter\mathpalette\expandafter\mathsm@sh\else
 \expandafter\makesm@sh\fi}
\def\finsm@sh{\iftop@\ht\z@\z@\fi\ifbot@\dp\z@\z@\fi\leavevmode\boxz@}
\def\LimitsOnSums{\global\let\slimits@\displaylimits}
\def\NoLimitsOnSums{\global\let\slimits@\nolimits}
\LimitsOnSums
\mathchardef\coprod@="1360       \def\coprod{\DOTSB\coprod@\slimits@}
\mathchardef\bigvee@="1357       \def\bigvee{\DOTSB\bigvee@\slimits@}
\mathchardef\bigwedge@="1356     \def\bigwedge{\DOTSB\bigwedge@\slimits@}
\mathchardef\biguplus@="1355     \def\biguplus{\DOTSB\biguplus@\slimits@}
\mathchardef\bigcap@="1354       \def\bigcap{\DOTSB\bigcap@\slimits@}
\mathchardef\bigcup@="1353       \def\bigcup{\DOTSB\bigcup@\slimits@}
\mathchardef\prod@="1351         \def\prod{\DOTSB\prod@\slimits@}
\mathchardef\sum@="1350          \def\sum{\DOTSB\sum@\slimits@}
\mathchardef\bigotimes@="134E    \def\bigotimes{\DOTSB\bigotimes@\slimits@}
\mathchardef\bigoplus@="134C     \def\bigoplus{\DOTSB\bigoplus@\slimits@}
\mathchardef\bigodot@="134A      \def\bigodot{\DOTSB\bigodot@\slimits@}
\mathchardef\bigsqcup@="1346     \def\bigsqcup{\DOTSB\bigsqcup@\slimits@}
\def\LimitsOnInts{\global\let\ilimits@\displaylimits}
\def\NoLimitsOnInts{\global\let\ilimits@\nolimits}
\NoLimitsOnInts
\def\int{\DOTSI\intop\ilimits@}
\def\oint{\DOTSI\ointop\ilimits@}
\def\intic@{\mathchoice{\hskip.5em}{\hskip.4em}{\hskip.4em}{\hskip.4em}}
\def\negintic@{\mathchoice
 {\hskip-.5em}{\hskip-.4em}{\hskip-.4em}{\hskip-.4em}}
\def\intkern@{\mathchoice{\!\!\!}{\!\!}{\!\!}{\!\!}}
\def\intdots@{\mathchoice{\plaincdots@}
 {{\cdotp}\mkern1.5mu{\cdotp}\mkern1.5mu{\cdotp}}
 {{\cdotp}\mkern1mu{\cdotp}\mkern1mu{\cdotp}}
 {{\cdotp}\mkern1mu{\cdotp}\mkern1mu{\cdotp}}}
\newcount\intno@
\def\iint{\DOTSI\intno@\tw@\FN@\ints@}
\def\iiint{\DOTSI\intno@\thr@@\FN@\ints@}
\def\iiiint{\DOTSI\intno@4 \FN@\ints@}
\def\idotsint{\DOTSI\intno@\z@\FN@\ints@}
\def\ints@{\findlimits@\ints@@}
\newif\iflimtoken@
\newif\iflimits@
\def\findlimits@{\limtoken@true\ifx\next\limits\limits@true
 \else\ifx\next\nolimits\limits@false\else
 \limtoken@false\ifx\ilimits@\nolimits\limits@false\else
 \ifinner\limits@false\else\limits@true\fi\fi\fi\fi}
\def\multint@{\int\ifnum\intno@=\z@\intdots@                                
 \else\intkern@\fi                                                          
 \ifnum\intno@>\tw@\int\intkern@\fi                                         
 \ifnum\intno@>\thr@@\int\intkern@\fi                                       
 \int}                                                                      
\def\multintlimits@{\intop\ifnum\intno@=\z@\intdots@\else\intkern@\fi
 \ifnum\intno@>\tw@\intop\intkern@\fi
 \ifnum\intno@>\thr@@\intop\intkern@\fi\intop}
\def\ints@@{\iflimtoken@                                                    
 \def\ints@@@{\iflimits@\negintic@\mathop{\intic@\multintlimits@}\limits    
  \else\multint@\nolimits\fi                                                
  \eat@}                                                                    
 \else                                                                      
 \def\ints@@@{\iflimits@\negintic@
  \mathop{\intic@\multintlimits@}\limits\else
  \multint@\nolimits\fi}\fi\ints@@@}
\def\LimitsOnNames{\global\let\nlimits@\displaylimits}
\def\NoLimitsOnNames{\global\let\nlimits@\nolimits@}
\LimitsOnNames
\def\nolimits@{\relaxnext@
 \DN@{\ifx\next\limits\DN@\limits{\nolimits}\else
  \let\next@\nolimits\fi\next@}%
 \FN@\next@}
\def\newmcodes@{\mathcode`\'="0027 \mathcode`\*="002A \mathcode`\.="613A
 \mathcode`\-="002D \mathcode`\/="002F \mathcode`\:="603A }
\def\operatorname#1{\mathop{\newmcodes@\kern\z@\fam\z@#1}\nolimits@}
\def\operatornamewithlimits#1{\mathop{\newmcodes@\kern\z@\fam\z@#1}\nlimits@}
\def\qopname@#1{\mathop{\fam\z@#1}\nolimits@}
\def\qopnamewl@#1{\mathop{\fam\z@#1}\nlimits@}
\def\arccos{\qopname@{arccos}}
\def\arcsin{\qopname@{arcsin}}
\def\arctan{\qopname@{arctan}}
\def\arg{\qopname@{arg}}
\def\cos{\qopname@{cos}}
\def\cosh{\qopname@{cosh}}
\def\cot{\qopname@{cot}}
\def\coth{\qopname@{coth}}
\def\csc{\qopname@{csc}}
\def\deg{\qopname@{deg}}
\def\det{\qopnamewl@{det}}
\def\dim{\qopname@{dim}}
\def\exp{\qopname@{exp}}
\def\gcd{\qopnamewl@{gcd}}
\def\hom{\qopname@{hom}}
\def\inf{\qopnamewl@{inf}}
\def\injlim{\qopnamewl@{inj\,lim}}
\def\ker{\qopname@{ker}}
\def\lg{\qopname@{lg}}
\def\lim{\qopnamewl@{lim}}
\def\liminf{\qopnamewl@{lim\,inf}}
\def\limsup{\qopnamewl@{lim\,sup}}
\def\ln{\qopname@{ln}}
\def\log{\qopname@{log}}
\def\max{\qopnamewl@{max}}
\def\min{\qopnamewl@{min}}
\def\Pr{\qopnamewl@{Pr}}
\def\projlim{\qopnamewl@{proj\,lim}}
\def\sec{\qopname@{sec}}
\def\sin{\qopname@{sin}}
\def\sinh{\qopname@{sinh}}
\def\sup{\qopnamewl@{sup}}
\def\tan{\qopname@{tan}}
\def\tanh{\qopname@{tanh}}
\def\varinjlim{\mathop{\vtop{\ialign{##\crcr
 \hfil\rm lim\hfil\crcr\noalign{\nointerlineskip}\rightarrowfill\crcr
 \noalign{\nointerlineskip\kern-\ex@}\crcr}}}}
\def\varprojlim{\mathop{\vtop{\ialign{##\crcr
 \hfil\rm lim\hfil\crcr\noalign{\nointerlineskip}\leftarrowfill\crcr
 \noalign{\nointerlineskip\kern-\ex@}\crcr}}}}
\def\varliminf{\mathop{\underline{\vrule height\z@ depth.2exwidth\z@
 \hbox{\rm lim}}}}

\newdimen\buffer@
\buffer@\fontdimen13 \tenex
\newdimen\buffer
\buffer\buffer@

\def\ResetBuffer{\fontdimen13 \tenex\buffer@\global\buffer\buffer@}
\def\shave#1{\mathop{\hbox{$\m@th\fontdimen13 \tenex\z@                     
 \displaystyle{#1}$}}\fontdimen13 \tenex\buffer}

\Invalid@\\
\def\Let@{\relax\iffalse{\fi\let\\=\cr\iffalse}\fi}
\Invalid@\vspace
\def\vspace@{\def\vspace##1{\crcr\noalign{\vskip##1\relax}}}
\def\multilimits@{\bgroup\vspace@\Let@
 \baselineskip\fontdimen10 \scriptfont\tw@
 \advance\baselineskip\fontdimen12 \scriptfont\tw@
 \lineskip\thr@@\fontdimen8 \scriptfont\thr@@
 \lineskiplimit\lineskip
 \vbox\bgroup\ialign\bgroup\hfil$\m@th\scriptstyle{##}$\hfil\crcr}
\def\Sb{_\multilimits@}
\def\endSb{\crcr\egroup\egroup\egroup}
\def\Sp{^\multilimits@}

\def\spreadlines#1{\RIfMIfI@\onlydmatherr@\spreadlines\else
 \openup#1\relax\fi\else\onlydmatherr@\spreadlines\fi}
\def\Mathstrut@{\copy\Mathstrutbox@}
\newbox\Mathstrutbox@
\setbox\Mathstrutbox@\null
\setboxz@h{$\m@th($}
\ht\Mathstrutbox@\ht\z@
\dp\Mathstrutbox@\dp\z@
\newdimen\spreadmlines@
\def\spreadmatrixlines#1{\RIfMIfI@
 \onlydmatherr@\spreadmatrixlines\else
 \spreadmlines@#1\relax\fi\else\onlydmatherr@\spreadmatrixlines\fi}
\def\matrix{\null\,\vcenter\bgroup\Let@\vspace@
 \normalbaselines\openup\spreadmlines@\ialign
 \bgroup\hfil$\m@th##$\hfil&&\quad\hfil$\m@th##$\hfil\crcr
 \Mathstrut@\crcr\noalign{\kern-\baselineskip}}
\def\endmatrix{\crcr\Mathstrut@\crcr\noalign{\kern-\baselineskip}\egroup
 \egroup\,}
\def\format{\crcr\egroup\iffalse{\fi\ifnum`}=0 \fi\format@}
\newtoks\hashtoks@
\hashtoks@{#}
\def\format@#1\\{\def\preamble@{#1}%
 \def\l{$\m@th\the\hashtoks@$\hfil}%
 \def\c{\hfil$\m@th\the\hashtoks@$\hfil}%
 \def\r{\hfil$\m@th\the\hashtoks@$}%
 \edef\Preamble@{\preamble@}\ifnum`{=0 \fi\iffalse}\fi
 \ialign\bgroup\span\Preamble@\crcr}
\def\smallmatrix{\null\,\vcenter\bgroup\vspace@\Let@
 \baselineskip9\ex@\lineskip\ex@
 \ialign\bgroup\hfil$\m@th\scriptstyle{##}$\hfil&&\thickspace\hfil
 $\m@th\scriptstyle{##}$\hfil\crcr}
\def\endsmallmatrix{\crcr\egroup\egroup\,}

\newmuskip\dotsspace@
\dotsspace@1.5mu
\def\strip@#1 {#1}
\def\spacehdots#1\for#2{\multispan{#2}\xleaders
 \hbox{$\m@th\mkern\strip@#1 \dotsspace@.\mkern\strip@#1 \dotsspace@$}\hfill}
\def\hdotsfor#1{\spacehdots\@ne\for{#1}}
\def\multispan@#1{\omit\mscount#1\unskip\loop\ifnum\mscount>\@ne\sp@n\repeat}
\def\spaceinnerhdots#1\for#2\after#3{\multispan@{\strip@#2 }#3\xleaders
 \hbox{$\m@th\mkern\strip@#1 \dotsspace@.\mkern\strip@#1 \dotsspace@$}\hfill}
\def\innerhdotsfor#1\after#2{\spaceinnerhdots\@ne\for#1\after{#2}}
\def\cases{\bgroup\spreadmlines@\jot\left\{\,\matrix\format\l&\quad\l\\}
\def\endcases{\endmatrix\right.\egroup}
\newif\ifinany@
\newif\ifinalign@
\newif\ifingather@
\def\strut@{\copy\strutbox@}
\newbox\strutbox@
\setbox\strutbox@\hbox{\vrule height8\p@ depth3\p@ width\z@}
\def\topaligned{\null\,\vtop\aligned@}
\def\botaligned{\null\,\vbox\aligned@}
\def\aligned{\null\,\vcenter\aligned@}
\def\aligned@{\bgroup\vspace@\Let@
 \ifinany@\else\openup\jot\fi\ialign
 \bgroup\hfil\strut@$\m@th\displaystyle{##}$&
 $\m@th\displaystyle{{}##}$\hfil\crcr}
\def\endaligned{\crcr\egroup\egroup}

\def\alignedat#1{\null\,\vcenter\bgroup\doat@{#1}\vspace@\Let@
 \ifinany@\else\openup\jot\fi\ialign\bgroup\span\preamble@@\crcr}
\newcount\atcount@
\def\doat@#1{\toks@{\hfil\strut@$\m@th
 \displaystyle{\the\hashtoks@}$&$\m@th\displaystyle
 {{}\the\hashtoks@}$\hfil}
 \atcount@#1\relax\advance\atcount@\m@ne                                    
 \loop\ifnum\atcount@>\z@\toks@=\expandafter{\the\toks@&\hfil$\m@th
 \displaystyle{\the\hashtoks@}$&$\m@th
 \displaystyle{{}\the\hashtoks@}$\hfil}\advance
  \atcount@\m@ne\repeat                                                     
 \xdef\preamble@{\the\toks@}\xdef\preamble@@{\preamble@}}

\def\gathered{\null\,\vcenter\bgroup\vspace@\Let@
 \ifinany@\else\openup\jot\fi\ialign
 \bgroup\hfil\strut@$\m@th\displaystyle{##}$\hfil\crcr}
\def\endgathered{\crcr\egroup\egroup}
\newif\iftagsleft@
\def\TagsOnLeft{\global\tagsleft@true}
\def\TagsOnRight{\global\tagsleft@false}
\TagsOnLeft
\newif\ifmathtags@
\def\TagsAsMath{\global\mathtags@true}
\def\TagsAsText{\global\mathtags@false}
\TagsAsText
\def\tagform@#1{\hbox{\rm(\ignorespaces#1\unskip)}}
\def\thetag{\leavevmode\tagform@}
\def\tag#1$${\iftagsleft@\leqno\else\eqno\fi                                
 \maketag@#1\maketag@                                                       
 $$}                                                                        
\def\maketag@{\FN@\maketag@@}
\def\maketag@@{\ifx\next"\expandafter\maketag@@@\else\expandafter\maketag@@@@
 \fi}
\def\maketag@@@"#1"#2\maketag@{\hbox{\rm#1}}                                
\def\maketag@@@@#1\maketag@{\ifmathtags@\tagform@{$\m@th#1$}\else
 \tagform@{#1}\fi}
\interdisplaylinepenalty\@M
\def\allowdisplaybreaks{\RIfMIfI@
 \onlydmatherr@\allowdisplaybreaks\else
 \interdisplaylinepenalty\z@\fi\else\onlydmatherr@\allowdisplaybreaks\fi}
\Invalid@\allowdisplaybreak
\Invalid@\displaybreak
\Invalid@\intertext
\def\allowdisplaybreak@{\def\allowdisplaybreak{\crcr\noalign{\allowbreak}}}
\def\displaybreak@{\def\displaybreak{\crcr\noalign{\break}}}
\def\intertext@{\def\intertext##1{\crcr\noalign{\vskip\belowdisplayskip
 \vbox{\normalbaselines\noindent##1}\vskip\abovedisplayskip}}}
\newskip\centering@
\centering@\z@ plus\@m\p@
\def\align{\relax\ifingather@\DN@{\csname align (in
  \string\gather)\endcsname}\else
 \ifmmode\ifinner\DN@{\onlydmatherr@\align}\else
  \let\next@\align@\fi
 \else\DN@{\onlydmatherr@\align}\fi\fi\next@}
\newhelp\andhelp@
{An extra & here is so disastrous that you should probably exit^^J
and fix things up.}
\newif\iftag@
\newcount\and@
\def\align@{\inalign@true\inany@true
 \vspace@\allowdisplaybreak@\displaybreak@\intertext@
 \def\tag{\global\tag@true\ifnum\and@=\z@\DN@{&&}\else
          \DN@{&}\fi\next@}%
 \iftagsleft@\DN@{\csname align \endcsname}\else
  \DN@{\csname align \space\endcsname}\fi\next@}
\def\Tag@{\iftag@\else\errhelp\andhelp@\err@{Extra & on this line}\fi}
\newdimen\lwidth@
\newdimen\rwidth@
\newdimen\maxlwidth@
\newdimen\maxrwidth@
\newdimen\totwidth@
\def\measure@#1\endalign{\lwidth@\z@\rwidth@\z@\maxlwidth@\z@\maxrwidth@\z@
 \global\and@\z@                                                            
 \setbox@ne\vbox                                                            
  {\everycr{\noalign{\global\tag@false\global\and@\z@}}\Let@                
  \halign{\setboxz@h{$\m@th\displaystyle{\@lign##}$}
   \global\lwidth@\wdz@                                                     
   \ifdim\lwidth@>\maxlwidth@\global\maxlwidth@\lwidth@\fi                  
   \global\advance\and@\@ne                                                 
   &\setboxz@h{$\m@th\displaystyle{{}\@lign##}$}\global\rwidth@\wdz@        
   \ifdim\rwidth@>\maxrwidth@\global\maxrwidth@\rwidth@\fi                  
   \global\advance\and@\@ne                                                
   &\Tag@
   \eat@{##}\crcr#1\crcr}}
 \totwidth@\maxlwidth@\advance\totwidth@\maxrwidth@}                       
\def\displ@y@{\global\dt@ptrue\openup\jot
 \everycr{\noalign{\global\tag@false\global\and@\z@\ifdt@p\global\dt@pfalse
 \vskip-\lineskiplimit\vskip\normallineskiplimit\else
 \penalty\interdisplaylinepenalty\fi}}}
\def\black@#1{\noalign{\ifdim#1>\displaywidth
 \dimen@\prevdepth\nointerlineskip                                          
 \vskip-\ht\strutbox@\vskip-\dp\strutbox@                                   
 \vbox{\noindent\hbox to#1{\strut@\hfill}}
 \prevdepth\dimen@                                                          
 \fi}}
\expandafter\def\csname align \space\endcsname#1\endalign
 {\measure@#1\endalign\global\and@\z@                                       
 \ifingather@\everycr{\noalign{\global\and@\z@}}\else\displ@y@\fi           
 \Let@\tabskip\centering@                                                   
 \halign to\displaywidth
  {\hfil\strut@\setboxz@h{$\m@th\displaystyle{\@lign##}$}
  \global\lwidth@\wdz@\boxz@\global\advance\and@\@ne                        
  \tabskip\z@skip                                                           
  &\setboxz@h{$\m@th\displaystyle{{}\@lign##}$}
  \global\rwidth@\wdz@\boxz@\hfill\global\advance\and@\@ne                  
  \tabskip\centering@                                                       
  &\setboxz@h{\@lign\strut@\maketag@##\maketag@}
  \dimen@\displaywidth\advance\dimen@-\totwidth@
  \divide\dimen@\tw@\advance\dimen@\maxrwidth@\advance\dimen@-\rwidth@     
  \ifdim\dimen@<\tw@\wdz@\llap{\vtop{\normalbaselines\null\boxz@}}
  \else\llap{\boxz@}\fi                                                    
  \tabskip\z@skip                                                          
  \crcr#1\crcr                                                             
  \black@\totwidth@}}                                                      
\newdimen\lineht@
\expandafter\def\csname align \endcsname#1\endalign{\measure@#1\endalign
 \global\and@\z@
 \ifdim\totwidth@>\displaywidth\let\displaywidth@\totwidth@\else
  \let\displaywidth@\displaywidth\fi                                        
 \ifingather@\everycr{\noalign{\global\and@\z@}}\else\displ@y@\fi
 \Let@\tabskip\centering@\halign to\displaywidth
  {\hfil\strut@\setboxz@h{$\m@th\displaystyle{\@lign##}$}%
  \global\lwidth@\wdz@\global\lineht@\ht\z@                                 
  \boxz@\global\advance\and@\@ne
  \tabskip\z@skip&\setboxz@h{$\m@th\displaystyle{{}\@lign##}$}%
  \global\rwidth@\wdz@\ifdim\ht\z@>\lineht@\global\lineht@\ht\z@\fi         
  \boxz@\hfil\global\advance\and@\@ne
  \tabskip\centering@&\kern-\displaywidth@                                  
  \setboxz@h{\@lign\strut@\maketag@##\maketag@}%
  \dimen@\displaywidth\advance\dimen@-\totwidth@
  \divide\dimen@\tw@\advance\dimen@\maxlwidth@\advance\dimen@-\lwidth@
  \ifdim\dimen@<\tw@\wdz@
   \rlap{\vbox{\normalbaselines\boxz@\vbox to\lineht@{}}}\else
   \rlap{\boxz@}\fi
  \tabskip\displaywidth@\crcr#1\crcr\black@\totwidth@}}
\expandafter\def\csname align (in \string\gather)\endcsname
  #1\endalign{\vcenter{\align@#1\endalign}}
\Invalid@\endalign
\newif\ifxat@
\def\alignat{\RIfMIfI@\DN@{\onlydmatherr@\alignat}\else
 \DN@{\csname alignat \endcsname}\fi\else
 \DN@{\onlydmatherr@\alignat}\fi\next@}
\newif\ifmeasuring@
\newbox\savealignat@
\expandafter\def\csname alignat \endcsname#1#2\endalignat                   
 {\inany@true\xat@false
 \def\tag{\global\tag@true\count@#1\relax\multiply\count@\tw@
  \xdef\tag@{}\loop\ifnum\count@>\and@\xdef\tag@{&\tag@}\advance\count@\m@ne
  \repeat\tag@}%
 \vspace@\allowdisplaybreak@\displaybreak@\intertext@
 \displ@y@\measuring@true                                                   
 \setbox\savealignat@\hbox{$\m@th\displaystyle\Let@
  \attag@{#1}
  \vbox{\halign{\span\preamble@@\crcr#2\crcr}}$}%
 \measuring@false                                                           
 \Let@\attag@{#1}
 \tabskip\centering@\halign to\displaywidth
  {\span\preamble@@\crcr#2\crcr                                             
  \black@{\wd\savealignat@}}}                                               
\Invalid@\endalignat
\def\xalignat{\RIfMIfI@
 \DN@{\onlydmatherr@\xalignat}\else
 \DN@{\csname xalignat \endcsname}\fi\else
 \DN@{\onlydmatherr@\xalignat}\fi\next@}
\expandafter\def\csname xalignat \endcsname#1#2\endxalignat
 {\inany@true\xat@true
 \def\tag{\global\tag@true\def\tag@{}\count@#1\relax\multiply\count@\tw@
  \loop\ifnum\count@>\and@\xdef\tag@{&\tag@}\advance\count@\m@ne\repeat\tag@}%
 \vspace@\allowdisplaybreak@\displaybreak@\intertext@
 \displ@y@\measuring@true\setbox\savealignat@\hbox{$\m@th\displaystyle\Let@
 \attag@{#1}\vbox{\halign{\span\preamble@@\crcr#2\crcr}}$}%
 \measuring@false\Let@
 \attag@{#1}\tabskip\centering@\halign to\displaywidth
 {\span\preamble@@\crcr#2\crcr\black@{\wd\savealignat@}}}
\def\attag@#1{\let\Maketag@\maketag@\let\TAG@\Tag@                          
 \let\Tag@=0\let\maketag@=0
 \ifmeasuring@\def\llap@##1{\setboxz@h{##1}\hbox to\tw@\wdz@{}}%
  \def\rlap@##1{\setboxz@h{##1}\hbox to\tw@\wdz@{}}\else
  \let\llap@\llap\let\rlap@\rlap\fi                                         
 \toks@{\hfil\strut@$\m@th\displaystyle{\@lign\the\hashtoks@}$\tabskip\z@skip
  \global\advance\and@\@ne&$\m@th\displaystyle{{}\@lign\the\hashtoks@}$\hfil
  \ifxat@\tabskip\centering@\fi\global\advance\and@\@ne}
 \iftagsleft@
  \toks@@{\tabskip\centering@&\Tag@\kern-\displaywidth
   \rlap@{\@lign\maketag@\the\hashtoks@\maketag@}%
   \global\advance\and@\@ne\tabskip\displaywidth}\else
  \toks@@{\tabskip\centering@&\Tag@\llap@{\@lign\maketag@
   \the\hashtoks@\maketag@}\global\advance\and@\@ne\tabskip\z@skip}\fi      
 \atcount@#1\relax\advance\atcount@\m@ne
 \loop\ifnum\atcount@>\z@
 \toks@=\expandafter{\the\toks@&\hfil$\m@th\displaystyle{\@lign
  \the\hashtoks@}$\global\advance\and@\@ne
  \tabskip\z@skip&$\m@th\displaystyle{{}\@lign\the\hashtoks@}$\hfil\ifxat@
  \tabskip\centering@\fi\global\advance\and@\@ne}\advance\atcount@\m@ne
 \repeat                                                                    
 \xdef\preamble@{\the\toks@\the\toks@@}
 \xdef\preamble@@{\preamble@}
 \let\maketag@\Maketag@\let\Tag@\TAG@}                                      
\Invalid@\endxalignat
\def\xxalignat{\RIfMIfI@
 \DN@{\onlydmatherr@\xxalignat}\else\DN@{\csname xxalignat
  \endcsname}\fi\else
 \DN@{\onlydmatherr@\xxalignat}\fi\next@}
\expandafter\def\csname xxalignat \endcsname#1#2\endxxalignat{\inany@true
 \vspace@\allowdisplaybreak@\displaybreak@\intertext@
 \displ@y\setbox\savealignat@\hbox{$\m@th\displaystyle\Let@
 \xxattag@{#1}\vbox{\halign{\span\preamble@@\crcr#2\crcr}}$}%
 \Let@\xxattag@{#1}\tabskip\z@skip\halign to\displaywidth
 {\span\preamble@@\crcr#2\crcr\black@{\wd\savealignat@}}}
\def\xxattag@#1{\toks@{\tabskip\z@skip\hfil\strut@
 $\m@th\displaystyle{\the\hashtoks@}$&%
 $\m@th\displaystyle{{}\the\hashtoks@}$\hfil\tabskip\centering@&}%
 \atcount@#1\relax\advance\atcount@\m@ne\loop\ifnum\atcount@>\z@
 \toks@=\expandafter{\the\toks@&\hfil$\m@th\displaystyle{\the\hashtoks@}$%
  \tabskip\z@skip&$\m@th\displaystyle{{}\the\hashtoks@}$\hfil
  \tabskip\centering@}\advance\atcount@\m@ne\repeat
 \xdef\preamble@{\the\toks@\tabskip\z@skip}\xdef\preamble@@{\preamble@}}
\Invalid@\endxxalignat
\newdimen\gwidth@
\newdimen\gmaxwidth@
\def\gmeasure@#1\endgather{\gwidth@\z@\gmaxwidth@\z@\setbox@ne\vbox{\Let@
 \halign{\setboxz@h{$\m@th\displaystyle{##}$}\global\gwidth@\wdz@
 \ifdim\gwidth@>\gmaxwidth@\global\gmaxwidth@\gwidth@\fi
 &\eat@{##}\crcr#1\crcr}}}
\def\gather{\RIfMIfI@\DN@{\onlydmatherr@\gather}\else
 \ingather@true\inany@true\def\tag{&}%
 \vspace@\allowdisplaybreak@\displaybreak@\intertext@
 \displ@y\Let@
 \iftagsleft@\DN@{\csname gather \endcsname}\else
  \DN@{\csname gather \space\endcsname}\fi\fi
 \else\DN@{\onlydmatherr@\gather}\fi\next@}
\expandafter\def\csname gather \space\endcsname#1\endgather
 {\gmeasure@#1\endgather\tabskip\centering@
 \halign to\displaywidth{\hfil\strut@\setboxz@h{$\m@th\displaystyle{##}$}%
 \global\gwidth@\wdz@\boxz@\hfil&
 \setboxz@h{\strut@{\maketag@##\maketag@}}%
 \dimen@\displaywidth\advance\dimen@-\gwidth@
 \ifdim\dimen@>\tw@\wdz@\llap{\boxz@}\else
 \llap{\vtop{\normalbaselines\null\boxz@}}\fi
 \tabskip\z@skip\crcr#1\crcr\black@\gmaxwidth@}}
\newdimen\glineht@
\expandafter\def\csname gather \endcsname#1\endgather{\gmeasure@#1\endgather
 \ifdim\gmaxwidth@>\displaywidth\let\gdisplaywidth@\gmaxwidth@\else
 \let\gdisplaywidth@\displaywidth\fi\tabskip\centering@\halign to\displaywidth
 {\hfil\strut@\setboxz@h{$\m@th\displaystyle{##}$}%
 \global\gwidth@\wdz@\global\glineht@\ht\z@\boxz@\hfil&\kern-\gdisplaywidth@
 \setboxz@h{\strut@{\maketag@##\maketag@}}%
 \dimen@\displaywidth\advance\dimen@-\gwidth@
 \ifdim\dimen@>\tw@\wdz@\rlap{\boxz@}\else
 \rlap{\vbox{\normalbaselines\boxz@\vbox to\glineht@{}}}\fi
 \tabskip\gdisplaywidth@\crcr#1\crcr\black@\gmaxwidth@}}
\newif\ifctagsplit@
\def\CenteredTagsOnSplits{\global\ctagsplit@true}
\def\TopOrBottomTagsOnSplits{\global\ctagsplit@false}
\TopOrBottomTagsOnSplits
\def\split{\relax\ifinany@\let\next@\insplit@\else
 \ifmmode\ifinner\def\next@{\onlydmatherr@\split}\else
 \let\next@\outsplit@\fi\else
 \def\next@{\onlydmatherr@\split}\fi\fi\next@}
\def\insplit@{\global\setbox\z@\vbox\bgroup\vspace@\Let@\ialign\bgroup
 \hfil\strut@$\m@th\displaystyle{##}$&$\m@th\displaystyle{{}##}$\hfill\crcr}
\def\endsplit{\crcr\egroup\egroup\iftagsleft@\expandafter\lendsplit@\else
 \expandafter\rendsplit@\fi}
\def\rendsplit@{\global\setbox9 \vbox
 {\unvcopy\z@\global\setbox8 \lastbox\unskip}
 \setbox@ne\hbox{\unhcopy8 \unskip\global\setbox\tw@\lastbox
 \unskip\global\setbox\thr@@\lastbox}
 \global\setbox7 \hbox{\unhbox\tw@\unskip}
 \ifinalign@\ifctagsplit@                                                   
  \gdef\split@{\hbox to\wd\thr@@{}&
   \vcenter{\vbox{\moveleft\wd\thr@@\boxz@}}}
 \else\gdef\split@{&\vbox{\moveleft\wd\thr@@\box9}\crcr
  \box\thr@@&\box7}\fi                                                      
 \else                                                                      
  \ifctagsplit@\gdef\split@{\vcenter{\boxz@}}\else
  \gdef\split@{\box9\crcr\hbox{\box\thr@@\box7}}\fi
 \fi
 \split@}                                                                   
\def\lendsplit@{\global\setbox9\vtop{\unvcopy\z@}
 \setbox@ne\vbox{\unvcopy\z@\global\setbox8\lastbox}
 \setbox@ne\hbox{\unhcopy8\unskip\setbox\tw@\lastbox
  \unskip\global\setbox\thr@@\lastbox}
 \ifinalign@\ifctagsplit@                                                   
  \gdef\split@{\hbox to\wd\thr@@{}&
  \vcenter{\vbox{\moveleft\wd\thr@@\box9}}}
  \else                                                                     
  \gdef\split@{\hbox to\wd\thr@@{}&\vbox{\moveleft\wd\thr@@\box9}}\fi
 \else
  \ifctagsplit@\gdef\split@{\vcenter{\box9}}\else
  \gdef\split@{\box9}\fi
 \fi\split@}
\def\outsplit@#1$${\align\insplit@#1\endalign$$}
\newdimen\multlinegap@
\multlinegap@1em
\newdimen\multlinetaggap@
\multlinetaggap@1em
\def\MultlineGap#1{\global\multlinegap@#1\relax}
\def\multlinegap#1{\RIfMIfI@\onlydmatherr@\multlinegap\else
 \multlinegap@#1\relax\fi\else\onlydmatherr@\multlinegap\fi}
\def\nomultlinegap{\multlinegap{\z@}}
\def\multline{\RIfMIfI@
 \DN@{\onlydmatherr@\multline}\else
 \DN@{\multline@}\fi\else
 \DN@{\onlydmatherr@\multline}\fi\next@}
\newif\iftagin@
\def\tagin@#1{\tagin@false\in@\tag{#1}\ifin@\tagin@true\fi}
\def\multline@#1$${\inany@true\vspace@\allowdisplaybreak@\displaybreak@
 \tagin@{#1}\iftagsleft@\DN@{\multline@l#1$$}\else
 \DN@{\multline@r#1$$}\fi\next@}
\newdimen\mwidth@
\def\rmmeasure@#1\endmultline{%
 \def\shoveleft##1{##1}\def\shoveright##1{##1}
 \setbox@ne\vbox{\Let@\halign{\setboxz@h
  {$\m@th\@lign\displaystyle{}##$}\global\mwidth@\wdz@
  \crcr#1\crcr}}}
\newdimen\mlineht@
\newif\ifzerocr@
\newif\ifonecr@
\def\lmmeasure@#1\endmultline{\global\zerocr@true\global\onecr@false
 \everycr{\noalign{\ifonecr@\global\onecr@false\fi
  \ifzerocr@\global\zerocr@false\global\onecr@true\fi}}
  \def\shoveleft##1{##1}\def\shoveright##1{##1}%
 \setbox@ne\vbox{\Let@\halign{\setboxz@h
  {$\m@th\@lign\displaystyle{}##$}\ifonecr@\global\mwidth@\wdz@
  \global\mlineht@\ht\z@\fi\crcr#1\crcr}}}
\newbox\mtagbox@
\newdimen\ltwidth@
\newdimen\rtwidth@
\def\multline@l#1$${\iftagin@\DN@{\lmultline@@#1$$}\else
 \DN@{\setbox\mtagbox@\null\ltwidth@\z@\rtwidth@\z@
  \lmultline@@@#1$$}\fi\next@}
\def\lmultline@@#1\endmultline\tag#2$${%
 \setbox\mtagbox@\hbox{\maketag@#2\maketag@}
 \lmmeasure@#1\endmultline\dimen@\mwidth@\advance\dimen@\wd\mtagbox@
 \advance\dimen@\multlinetaggap@                                            
 \ifdim\dimen@>\displaywidth\ltwidth@\z@\else\ltwidth@\wd\mtagbox@\fi       
 \lmultline@@@#1\endmultline$$}
\def\lmultline@@@{\displ@y
 \def\shoveright##1{##1\hfilneg\hskip\multlinegap@}%
 \def\shoveleft##1{\setboxz@h{$\m@th\displaystyle{}##1$}%
  \setbox@ne\hbox{$\m@th\displaystyle##1$}%
  \hfilneg
  \iftagin@
   \ifdim\ltwidth@>\z@\hskip\ltwidth@\hskip\multlinetaggap@\fi
  \else\hskip\multlinegap@\fi\hskip.5\wd@ne\hskip-.5\wdz@##1}
  \halign\bgroup\Let@\hbox to\displaywidth
   {\strut@$\m@th\displaystyle\hfil{}##\hfil$}\crcr
   \hfilneg                                                                 
   \iftagin@                                                                
    \ifdim\ltwidth@>\z@                                                     
     \box\mtagbox@\hskip\multlinetaggap@                                    
    \else
     \rlap{\vbox{\normalbaselines\hbox{\strut@\box\mtagbox@}%
     \vbox to\mlineht@{}}}\fi                                               
   \else\hskip\multlinegap@\fi}                                             
\def\multline@r#1$${\iftagin@\DN@{\rmultline@@#1$$}\else
 \DN@{\setbox\mtagbox@\null\ltwidth@\z@\rtwidth@\z@
  \rmultline@@@#1$$}\fi\next@}
\def\rmultline@@#1\endmultline\tag#2$${\ltwidth@\z@
 \setbox\mtagbox@\hbox{\maketag@#2\maketag@}%
 \rmmeasure@#1\endmultline\dimen@\mwidth@\advance\dimen@\wd\mtagbox@
 \advance\dimen@\multlinetaggap@
 \ifdim\dimen@>\displaywidth\rtwidth@\z@\else\rtwidth@\wd\mtagbox@\fi
 \rmultline@@@#1\endmultline$$}
\def\rmultline@@@{\displ@y
 \def\shoveright##1{##1\hfilneg\iftagin@\ifdim\rtwidth@>\z@
  \hskip\rtwidth@\hskip\multlinetaggap@\fi\else\hskip\multlinegap@\fi}%
 \def\shoveleft##1{\setboxz@h{$\m@th\displaystyle{}##1$}%
  \setbox@ne\hbox{$\m@th\displaystyle##1$}%
  \hfilneg\hskip\multlinegap@\hskip.5\wd@ne\hskip-.5\wdz@##1}%
 \halign\bgroup\Let@\hbox to\displaywidth
  {\strut@$\m@th\displaystyle\hfil{}##\hfil$}\crcr
 \hfilneg\hskip\multlinegap@}
\def\endmultline{\iftagsleft@\expandafter\lendmultline@\else
 \expandafter\rendmultline@\fi}
\def\lendmultline@{\hfilneg\hskip\multlinegap@\crcr\egroup}
\def\rendmultline@{\iftagin@                                                
 \ifdim\rtwidth@>\z@                                                        
  \hskip\multlinetaggap@\box\mtagbox@                                       
 \else\llap{\vtop{\normalbaselines\null\hbox{\strut@\box\mtagbox@}}}\fi     
 \else\hskip\multlinegap@\fi                                                
 \hfilneg\crcr\egroup}
\def\bmod{\mskip-\medmuskip\mkern5mu\mathbin{\fam\z@ mod}\penalty900
 \mkern5mu\mskip-\medmuskip}
\def\pmod#1{\allowbreak\ifinner\mkern8mu\else\mkern18mu\fi
 ({\fam\z@ mod}\,\,#1)}
\def\pod#1{\allowbreak\ifinner\mkern8mu\else\mkern18mu\fi(#1)}
\def\mod#1{\allowbreak\ifinner\mkern12mu\else\mkern18mu\fi{\fam\z@ mod}\,\,#1}
\newcount\cfraccount@
\def\cfrac{\bgroup\bgroup\advance\cfraccount@\@ne\strut
 \iffalse{\fi\def\\{\over\displaystyle}\iffalse}\fi}
\def\lcfrac{\bgroup\bgroup\advance\cfraccount@\@ne\strut
 \iffalse{\fi\def\\{\hfill\over\displaystyle}\iffalse}\fi}
\def\rcfrac{\bgroup\bgroup\advance\cfraccount@\@ne\strut\hfill
 \iffalse{\fi\def\\{\over\displaystyle}\iffalse}\fi}
\def\gloop@#1\repeat{\gdef\body{#1}\iterate}
\def\endcfrac{\gloop@\ifnum\cfraccount@>\z@\global\advance\cfraccount@\m@ne
 \egroup\hskip-\nulldelimiterspace\egroup\repeat}
\def\binrel@#1{\setboxz@h{\thinmuskip0mu
  \medmuskip\m@ne mu\thickmuskip\@ne mu$#1\m@th$}%
 \setbox@ne\hbox{\thinmuskip0mu\medmuskip\m@ne mu\thickmuskip
  \@ne mu${}#1{}\m@th$}%
 \setbox\tw@\hbox{\hskip\wd@ne\hskip-\wdz@}}
\def\overset#1\to#2{\binrel@{#2}\ifdim\wd\tw@<\z@
 \mathbin{\mathop{\kern\z@#2}\limits^{#1}}\else\ifdim\wd\tw@>\z@
 \mathrel{\mathop{\kern\z@#2}\limits^{#1}}\else
 {\mathop{\kern\z@#2}\limits^{#1}}{}\fi\fi}
\def\underset#1\to#2{\binrel@{#2}\ifdim\wd\tw@<\z@
 \mathbin{\mathop{\kern\z@#2}\limits_{#1}}\else\ifdim\wd\tw@>\z@
 \mathrel{\mathop{\kern\z@#2}\limits_{#1}}\else
 {\mathop{\kern\z@#2}\limits_{#1}}{}\fi\fi}
\def\oversetbrace#1\to#2{\overbrace{#2}^{#1}}
\def\undersetbrace#1\to#2{\underbrace{#2}_{#1}}
\def\sideset#1\and#2\to#3{%
 \setbox@ne\hbox{$\dsize{\vphantom{#3}}#1{#3}\m@th$}%
 \setbox\tw@\hbox{$\dsize{#3}#2\m@th$}%
 \hskip\wd@ne\hskip-\wd\tw@\mathop{\hskip\wd\tw@\hskip-\wd@ne
  {\vphantom{#3}}#1{#3}#2}}
\def\rightarrowfill@#1{$#1\m@th\mathord-\mkern-6mu\cleaders
 \hbox{$#1\mkern-2mu\mathord-\mkern-2mu$}\hfill
 \mkern-6mu\mathord\rightarrow$}
\def\leftarrowfill@#1{$#1\m@th\mathord\leftarrow\mkern-6mu\cleaders
 \hbox{$#1\mkern-2mu\mathord-\mkern-2mu$}\hfill\mkern-6mu\mathord-$}
\def\leftrightarrowfill@#1{$#1\m@th\mathord\leftarrow\mkern-6mu\cleaders
 \hbox{$#1\mkern-2mu\mathord-\mkern-2mu$}\hfill
 \mkern-6mu\mathord\rightarrow$}
\def\overrightarrow{\mathpalette\overrightarrow@}
\def\overrightarrow@#1#2{\vbox{\ialign{##\crcr\rightarrowfill@#1\crcr
 \noalign{\kern-\ex@\nointerlineskip}$\m@th\hfil#1#2\hfil$\crcr}}}

\def\overleftarrow{\mathpalette\overleftarrow@}
\def\overleftarrow@#1#2{\vbox{\ialign{##\crcr\leftarrowfill@#1\crcr
 \noalign{\kern-\ex@\nointerlineskip}$\m@th\hfil#1#2\hfil$\crcr}}}
\def\overleftrightarrow{\mathpalette\overleftrightarrow@}
\def\overleftrightarrow@#1#2{\vbox{\ialign{##\crcr\leftrightarrowfill@#1\crcr
 \noalign{\kern-\ex@\nointerlineskip}$\m@th\hfil#1#2\hfil$\crcr}}}
\def\underrightarrow{\mathpalette\underrightarrow@}
\def\underrightarrow@#1#2{\vtop{\ialign{##\crcr$\m@th\hfil#1#2\hfil$\crcr
 \noalign{\nointerlineskip}\rightarrowfill@#1\crcr}}}

\def\underleftarrow{\mathpalette\underleftarrow@}
\def\underleftarrow@#1#2{\vtop{\ialign{##\crcr$\m@th\hfil#1#2\hfil$\crcr
 \noalign{\nointerlineskip}\leftarrowfill@#1\crcr}}}
\def\underleftrightarrow{\mathpalette\underleftrightarrow@}
\def\underleftrightarrow@#1#2{\vtop{\ialign{##\crcr$\m@th\hfil#1#2\hfil$\crcr
 \noalign{\nointerlineskip}\leftrightarrowfill@#1\crcr}}}
\let\DOTSI\relax
\let\DOTSB\relax

\newif\ifmath@
{\uccode`7=`\\ \uccode`8=`m \uccode`9=`a \uccode`0=`t \uccode`!=`h
 \uppercase{\gdef\math@#1#2#3#4#5#6\math@{\global\math@false\ifx 7#1\ifx 8#2%
 \ifx 9#3\ifx 0#4\ifx !#5\xdef\meaning@{#6}\global\math@true\fi\fi\fi\fi\fi}}}
\newif\ifmathch@
{\uccode`7=`c \uccode`8=`h \uccode`9=`\"
 \uppercase{\gdef\mathch@#1#2#3#4#5#6\mathch@{\global\mathch@false
  \ifx 7#1\ifx 8#2\ifx 9#5\global\mathch@true\xdef\meaning@{9#6}\fi\fi\fi}}}
\newcount\classnum@
\def\getmathch@#1.#2\getmathch@{\classnum@#1 \divide\classnum@4096
 \ifcase\number\classnum@\or\or\gdef\thedots@{\dotsb@}\or
 \gdef\thedots@{\dotsb@}\fi}
\newif\ifmathbin@
{\uccode`4=`b \uccode`5=`i \uccode`6=`n
 \uppercase{\gdef\mathbin@#1#2#3{\relaxnext@
  \DNii@##1\mathbin@{\ifx\space@\next\global\mathbin@true\fi}%
 \global\mathbin@false\DN@##1\mathbin@{}%
 \ifx 4#1\ifx 5#2\ifx 6#3\DN@{\FN@\nextii@}\fi\fi\fi\next@}}}
\newif\ifmathrel@
{\uccode`4=`r \uccode`5=`e \uccode`6=`l
 \uppercase{\gdef\mathrel@#1#2#3{\relaxnext@
  \DNii@##1\mathrel@{\ifx\space@\next\global\mathrel@true\fi}%
 \global\mathrel@false\DN@##1\mathrel@{}%
 \ifx 4#1\ifx 5#2\ifx 6#3\DN@{\FN@\nextii@}\fi\fi\fi\next@}}}
\newif\ifmacro@
{\uccode`5=`m \uccode`6=`a \uccode`7=`c
 \uppercase{\gdef\macro@#1#2#3#4\macro@{\global\macro@false
  \ifx 5#1\ifx 6#2\ifx 7#3\global\macro@true
  \xdef\meaning@{\macro@@#4\macro@@}\fi\fi\fi}}}
\def\macro@@#1->#2\macro@@{#2}
\newif\ifDOTS@
\newcount\DOTSCASE@
{\uccode`6=`\\ \uccode`7=`D \uccode`8=`O \uccode`9=`T \uccode`0=`S
 \uppercase{\gdef\DOTS@#1#2#3#4#5{\global\DOTS@false\DN@##1\DOTS@{}%
  \ifx 6#1\ifx 7#2\ifx 8#3\ifx 9#4\ifx 0#5\let\next@\DOTS@@\fi\fi\fi\fi\fi
  \next@}}}
{\uccode`3=`B \uccode`4=`I \uccode`5=`X
 \uppercase{\gdef\DOTS@@#1{\relaxnext@
  \DNii@##1\DOTS@{\ifx\space@\next\global\DOTS@true\fi}%
  \DN@{\FN@\nextii@}%
  \ifx 3#1\global\DOTSCASE@\z@\else
  \ifx 4#1\global\DOTSCASE@\@ne\else
  \ifx 5#1\global\DOTSCASE@\tw@\else\DN@##1\DOTS@{}%
  \fi\fi\fi\next@}}}
\newif\ifnot@
{\uccode`5=`\\ \uccode`6=`n \uccode`7=`o \uccode`8=`t
 \uppercase{\gdef\not@#1#2#3#4{\relaxnext@
  \DNii@##1\not@{\ifx\space@\next\global\not@true\fi}%
 \global\not@false\DN@##1\not@{}%
 \ifx 5#1\ifx 6#2\ifx 7#3\ifx 8#4\DN@{\FN@\nextii@}\fi\fi\fi
 \fi\next@}}}
\newif\ifkeybin@
\def\keybin@{\keybin@true
 \ifx\next+\else\ifx\next=\else\ifx\next<\else\ifx\next>\else\ifx\next-\else
 \ifx\next*\else\ifx\next:\else\keybin@false\fi\fi\fi\fi\fi\fi\fi}
\def\dots{\RIfM@\expandafter\mdots@\else\expandafter\tdots@\fi}
\def\tdots@{\unskip\relaxnext@
 \DN@{$\m@th\mathinner{\ldotp\ldotp\ldotp}\,
   \ifx\next,\,$\else\ifx\next.\,$\else\ifx\next;\,$\else\ifx\next:\,$\else
   \ifx\next?\,$\else\ifx\next!\,$\else$ \fi\fi\fi\fi\fi\fi}%
 \ \FN@\next@}
\def\mdots@{\FN@\mdots@@}
\def\mdots@@{\gdef\thedots@{\dotso@}
 \ifx\next\boldkey\gdef\thedots@\boldkey{\boldkeydots@}\else                
 \ifx\next\boldsymbol\gdef\thedots@\boldsymbol{\boldsymboldots@}\else       
 \ifx,\next\gdef\thedots@{\dotsc}
 \else\ifx\not\next\gdef\thedots@{\dotsb@}
 \else\keybin@
 \ifkeybin@\gdef\thedots@{\dotsb@}
 \else\xdef\meaning@{\meaning\next..........}\xdef\meaning@@{\meaning@}
  \expandafter\math@\meaning@\math@
  \ifmath@
   \expandafter\mathch@\meaning@\mathch@
   \ifmathch@\expandafter\getmathch@\meaning@\getmathch@\fi                 
  \else\expandafter\macro@\meaning@@\macro@                                 
  \ifmacro@                                                                
   \expandafter\not@\meaning@\not@\ifnot@\gdef\thedots@{\dotsb@}
  \else\expandafter\DOTS@\meaning@\DOTS@
  \ifDOTS@
   \ifcase\number\DOTSCASE@\gdef\thedots@{\dotsb@}%
    \or\gdef\thedots@{\dotsi}\else\fi                                      
  \else\expandafter\math@\meaning@\math@                                   
  \ifmath@\expandafter\mathbin@\meaning@\mathbin@
  \ifmathbin@\gdef\thedots@{\dotsb@}
  \else\expandafter\mathrel@\meaning@\mathrel@
  \ifmathrel@\gdef\thedots@{\dotsb@}
  \fi\fi\fi\fi\fi\fi\fi\fi\fi\fi\fi\fi
 \thedots@}
\def\plainldots@{\mathinner{\ldotp\ldotp\ldotp}}
\def\plaincdots@{\mathinner{\cdotp\cdotp\cdotp}}
\def\dotsi{\!\plaincdots@}
\let\dotsb@\plaincdots@
\newif\ifextra@
\newif\ifrightdelim@
\def\rightdelim@{\global\rightdelim@true                                    
 \ifx\next)\else                                                            
 \ifx\next]\else
 \ifx\next\rbrack\else
 \ifx\next\}\else
 \ifx\next\rbrace\else
 \ifx\next\rangle\else
 \ifx\next\rceil\else
 \ifx\next\rfloor\else
 \ifx\next\rgroup\else
 \ifx\next\rmoustache\else
 \ifx\next\right\else
 \ifx\next\bigr\else
 \ifx\next\biggr\else
 \ifx\next\Bigr\else                                                        
 \ifx\next\Biggr\else\global\rightdelim@false
 \fi\fi\fi\fi\fi\fi\fi\fi\fi\fi\fi\fi\fi\fi\fi}
\def\extra@{%
 \global\extra@false\rightdelim@\ifrightdelim@\global\extra@true            
 \else\ifx\next$\global\extra@true                                          
 \else\xdef\meaning@{\meaning\next..........}
 \expandafter\macro@\meaning@\macro@\ifmacro@                               
 \expandafter\DOTS@\meaning@\DOTS@
 \ifDOTS@
 \ifnum\DOTSCASE@=\tw@\global\extra@true                                    
 \fi\fi\fi\fi\fi}
\newif\ifbold@
\def\dotso@{\relaxnext@
 \ifbold@
  \let\next\delayed@
  \DNii@{\extra@\plainldots@\ifextra@\,\fi}%
 \else
  \DNii@{\DN@{\extra@\plainldots@\ifextra@\,\fi}\FN@\next@}%
 \fi
 \nextii@}
\def\extrap@#1{%
 \ifx\next,\DN@{#1\,}\else
 \ifx\next;\DN@{#1\,}\else
 \ifx\next.\DN@{#1\,}\else\extra@
 \ifextra@\DN@{#1\,}\else
 \let\next@#1\fi\fi\fi\fi\next@}
\def\ldots{\DN@{\extrap@\plainldots@}%
 \FN@\next@}
\def\cdots{\DN@{\extrap@\plaincdots@}%
 \FN@\next@}

\def\dotsc{\relaxnext@
 \DN@{\ifx\next;\plainldots@\,\else
  \ifx\next.\plainldots@\,\else\extra@\plainldots@
  \ifextra@\,\fi\fi\fi}%
 \FN@\next@}
\def\cdot{\mathchar"2201 }
\def\longrightarrow{\DOTSB\relbar\joinrel\rightarrow}

\def\mapsto{\DOTSB\mapstochar\rightarrow}

\def\dddot#1{{\mathop{#1}\limits^{\vbox to-1.4\ex@{\kern-\tw@\ex@
 \hbox{\rm...}\vss}}}}
\def\ddddot#1{{\mathop{#1}\limits^{\vbox to-1.4\ex@{\kern-\tw@\ex@
 \hbox{\rm....}\vss}}}}
\def\sphat{^{\mathchoice{}{}%
 {\,\,\botsmash{\hbox{\lower4\ex@\hbox{$\m@th\widehat{\null}$}}}}%
 {\,\botsmash{\hbox{\lower3\ex@\hbox{$\m@th\hat{\null}$}}}}}}

\def\spacute{^{\!\botsmash{\hbox{\lower\@ne ex\hbox{\'{}}}}}}
\def\spgrave{^{\mathchoice{}{}{}{\!}%
 \botsmash{\hbox{\lower\@ne ex\hbox{\`{}}}}}}
\def\spdot{^{\hbox{\raise\ex@\hbox{\rm.}}}}
\def\spddot{^{\hbox{\raise\ex@\hbox{\rm..}}}}
\def\spdddot{^{\hbox{\raise\ex@\hbox{\rm...}}}}
\def\spddddot{^{\hbox{\raise\ex@\hbox{\rm....}}}}
\def\spbreve{^{\!\botsmash{\hbox{\lower4\ex@\hbox{\u{}}}}}}

\def\textonlyfont@#1#2{\def#1{\RIfM@
 \Err@{Use \string#1\space only in text}\else#2\fi}}
\textonlyfont@\rm\tenrm
\textonlyfont@\it\tenit
\textonlyfont@\sl\tensl
\textonlyfont@\bf\tenbf
\def\oldnos#1{\RIfM@{\mathcode`\,="013B \fam\@ne#1}\else
 \leavevmode\hbox{$\m@th\mathcode`\,="013B \fam\@ne#1$}\fi}
\def\text{\RIfM@\expandafter\text@\else\expandafter\text@@\fi}
\def\text@@#1{\leavevmode\hbox{#1}}
\def\mathhexbox@#1#2#3{\text{$\m@th\mathchar"#1#2#3$}}
\def\dag{{\mathhexbox@279}}
\def\ddag{{\mathhexbox@27A}}
\def\S{{\mathhexbox@278}}
\def\P{{\mathhexbox@27B}}
\newif\iffirstchoice@
\firstchoice@true
\def\text@#1{\mathchoice
 {\hbox{\everymath{\displaystyle}\def\textfonti{\the\textfont\@ne}%
  \def\textfontii{\the\textfont\tw@}\textdef@@ T#1}}
 {\hbox{\firstchoice@false
  \everymath{\textstyle}\def\textfonti{\the\textfont\@ne}%
  \def\textfontii{\the\textfont\tw@}\textdef@@ T#1}}
 {\hbox{\firstchoice@false
  \everymath{\scriptstyle}\def\textfonti{\the\scriptfont\@ne}%
  \def\textfontii{\the\scriptfont\tw@}\textdef@@ S\rm#1}}
 {\hbox{\firstchoice@false
  \everymath{\scriptscriptstyle}\def\textfonti
  {\the\scriptscriptfont\@ne}%
  \def\textfontii{\the\scriptscriptfont\tw@}\textdef@@ s\rm#1}}}
\def\textdef@@#1{\textdef@#1\rm\textdef@#1\bf\textdef@#1\sl\textdef@#1\it}
\def\rmfam{0}
\def\textdef@#1#2{%
 \DN@{\csname\expandafter\eat@\string#2fam\endcsname}%
 \if S#1\edef#2{\the\scriptfont\next@\relax}%
 \else\if s#1\edef#2{\the\scriptscriptfont\next@\relax}%
 \else\edef#2{\the\textfont\next@\relax}\fi\fi}
\scriptfont\itfam\tenit \scriptscriptfont\itfam\tenit
\scriptfont\slfam\tensl \scriptscriptfont\slfam\tensl
\newif\iftopfolded@
\newif\ifbotfolded@
\def\topfoldedtext{\topfolded@true\botfolded@false\foldedtext@}
\def\botfoldedtext{\botfolded@true\topfolded@false\foldedtext@}
\def\foldedtext{\topfolded@false\botfolded@false\foldedtext@}
\Invalid@\foldedwidth
\def\foldedtext@{\relaxnext@
 \DN@{\ifx\next\foldedwidth\let\next@\nextii@\else
  \DN@{\nextii@\foldedwidth{.3\hsize}}\fi\next@}%
 \DNii@\foldedwidth##1##2{\setbox\z@\vbox
  {\normalbaselines\hsize##1\relax
  \tolerance1600 \noindent\ignorespaces##2}\ifbotfolded@\boxz@\else
  \iftopfolded@\vtop{\unvbox\z@}\else\vcenter{\boxz@}\fi\fi}%
 \FN@\next@}
\def\bold{\RIfM@\expandafter\bold@\else
 \expandafter\nonmatherr@\expandafter\bold\fi}
\def\bold@#1{{\bold@@{#1}}}
\def\bold@@#1{\fam\bffam\relax#1}
\def\slanted{\RIfM@\expandafter\slanted@\else
 \expandafter\nonmatherr@\expandafter\slanted\fi}
\def\slanted@#1{{\slanted@@{#1}}}
\def\slanted@@#1{\fam\slfam\relax#1}
\def\roman{\RIfM@\expandafter\roman@\else
 \expandafter\nonmatherr@\expandafter\roman\fi}
\def\roman@#1{{\roman@@{#1}}}
\def\roman@@#1{\fam\rmfam\relax#1}
\def\italic{\RIfM@\expandafter\italic@\else
 \expandafter\nonmatherr@\expandafter\italic\fi}
\def\italic@#1{{\italic@@{#1}}}
\def\italic@@#1{\fam\itfam\relax#1}
\def\Cal{\RIfM@\expandafter\Cal@\else
 \expandafter\nonmatherr@\expandafter\Cal\fi}
\def\Cal@#1{{\Cal@@{#1}}}
\def\Cal@@#1{\noaccents@\fam\tw@#1}
\mathchardef\Gamma="0000
\mathchardef\Delta="0001
\mathchardef\Theta="0002
\mathchardef\Lambda="0003
\mathchardef\Xi="0004
\mathchardef\Pi="0005
\mathchardef\Sigma="0006
\mathchardef\Upsilon="0007
\mathchardef\Phi="0008
\mathchardef\Psi="0009
\mathchardef\Omega="000A
\mathchardef\varGamma="0100
\mathchardef\varDelta="0101
\mathchardef\varTheta="0102
\mathchardef\varLambda="0103
\mathchardef\varXi="0104
\mathchardef\varPi="0105
\mathchardef\varSigma="0106
\mathchardef\varUpsilon="0107
\mathchardef\varPhi="0108
\mathchardef\varPsi="0109
\mathchardef\varOmega="010A
\newif\ifmsamloaded@
\newif\ifmsbmloaded@
\newif\ifeufmloaded@
\let\alloc@@\alloc@
\def\hexnumber@#1{\ifcase#1 0\or 1\or 2\or 3\or 4\or 5\or 6\or 7\or 8\or
 9\or A\or B\or C\or D\or E\or F\fi}
\edef\bffam@{\hexnumber@\bffam}
\def\loadmsam{\msamloaded@true
 \font@\tenmsa=msam10
 \font@\sevenmsa=msam7
 \font@\fivemsa=msam5
 \alloc@@8\fam\chardef\sixt@@n\msafam
 \textfont\msafam=\tenmsa
 \scriptfont\msafam=\sevenmsa
 \scriptscriptfont\msafam=\fivemsa
 \edef\msafam@{\hexnumber@\msafam}%
 \mathchardef\dabar@"0\msafam@39
 \def\dashrightarrow{\mathrel{\dabar@\dabar@\mathchar"0\msafam@4B}}%
 \def\dashleftarrow{\mathrel{\mathchar"0\msafam@4C\dabar@\dabar@}}%
 \let\dasharrow\dashrightarrow
 \def\ulcorner{\delimiter"4\msafam@70\msafam@70 }
 \def\urcorner{\delimiter"5\msafam@71\msafam@71 }
 \def\llcorner{\delimiter"4\msafam@78\msafam@78 }
 \def\lrcorner{\delimiter"5\msafam@79\msafam@79 }
 \def\yen{{\mathhexbox@\msafam@55 }}
 \def\checkmark{{\mathhexbox@\msafam@58 }}
 \def\circledR{{\mathhexbox@\msafam@72 }}
 \def\maltese{{\mathhexbox@\msafam@7A }}}
\def\loadmsbm{\msbmloaded@true
 \font@\tenmsb=msbm10
 \font@\sevenmsb=msbm7
 \font@\fivemsb=msbm5
 \alloc@@8\fam\chardef\sixt@@n\msbfam
 \textfont\msbfam=\tenmsb
 \scriptfont\msbfam=\sevenmsb
 \scriptscriptfont\msbfam=\fivemsb
 \edef\msbfam@{\hexnumber@\msbfam}%
 }
\def\widehat#1{\ifmsbmloaded@
  \setboxz@h{$\m@th#1$}\ifdim\wdz@>\tw@ em\mathaccent"0\msbfam@5B{#1}\else
  \mathaccent"0362{#1}\fi
 \else\mathaccent"0362{#1}\fi}
\def\widetilde#1{\ifmsbmloaded@
  \setboxz@h{$\m@th#1$}\ifdim\wdz@>\tw@ em\mathaccent"0\msbfam@5D{#1}\else
  \mathaccent"0365{#1}\fi
 \else\mathaccent"0365{#1}\fi}
\def\newsymbol#1#2#3#4#5{\define#1{}\let\next@\relax
 \ifnum#2=\@ne\ifmsamloaded@\let\next@\msafam@\fi\else
 \ifnum#2=\tw@\ifmsbmloaded@\let\next@\msbfam@\fi\fi\fi
 \ifx\next@\relax
  \ifnum#2>\tw@\Err@{\Invalid@@\string\newsymbol}\else
  \ifnum#2=\@ne\Err@{You must first \string\loadmsam}\else
   \Err@{You must first \string\loadmsbm}\fi\fi
 \else
  \mathchardef#1="#3\next@#4#5
 \fi}

\def\Bbb{\RIfM@\expandafter\Bbb@\else
 \expandafter\nonmatherr@\expandafter\Bbb\fi}
\def\Bbb@#1{{\Bbb@@{#1}}}
\def\Bbb@@#1{\noaccents@\fam\msbfam\relax#1}
\def\loadeufm{\eufmloaded@true
 \font@\teneufm=eufm10
 \font@\seveneufm=eufm7
 \font@\fiveeufm=eufm5
 \alloc@@8\fam\chardef\sixt@@n\eufmfam
 \textfont\eufmfam=\teneufm
 \scriptfont\eufmfam=\seveneufm
 \scriptscriptfont\eufmfam=\fiveeufm}
\def\frak{\RIfM@\expandafter\frak@\else
 \expandafter\nonmatherr@\expandafter\frak\fi}
\def\frak@#1{{\frak@@{#1}}}
\def\frak@@#1{\fam\eufmfam\relax#1}

\newif\ifcmmibloaded@
\newif\ifcmbsyloaded@
\def\loadbold{\cmmibloaded@true\cmbsyloaded@true
 \font@\tencmmib=cmmib10 \font@\sevencmmib=cmmib7 \font@\fivecmmib=cmmib5
 \skewchar\tencmmib='177 \skewchar\sevencmmib='177 \skewchar\fivecmmib='177
 \alloc@@8\fam\chardef\sixt@@n\cmmibfam
 \textfont\cmmibfam=\tencmmib
 \scriptfont\cmmibfam=\sevencmmib
 \scriptscriptfont\cmmibfam=\fivecmmib
 \edef\cmmibfam@{\hexnumber@\cmmibfam}%
 \font@\tencmbsy=cmbsy10 \font@\sevencmbsy=cmbsy7 \font@\fivecmbsy=cmbsy5
 \skewchar\tencmbsy='60 \skewchar\sevencmbsy='60 \skewchar\fivecmbsy='60
 \alloc@@8\fam\chardef\sixt@@n\cmbsyfam
 \textfont\cmbsyfam=\tencmbsy
 \scriptfont\cmbsyfam=\sevencmbsy
 \scriptscriptfont\cmbsyfam=\fivecmbsy
 \edef\cmbsyfam@{\hexnumber@\cmbsyfam}}
\def\mathchari@#1#2#3{\ifcmmibloaded@\mathchar"#1\cmmibfam@#2#3 \else
 \Err@{First bold symbol font not loaded}\fi}
\def\mathcharii@#1#2#3{\ifcmbsyloaded@\mathchar"#1\cmbsyfam@#2#3 \else
 \Err@{Second bold symbol font not loaded}\fi}
\def\boldkey#1{\ifcat\noexpand#1A%
  \ifcmmibloaded@{\fam\cmmibfam#1}\else
   \Err@{First bold symbol font not loaded}\fi
 \else
 \ifx#1!\mathchar"5\bffam@21 \else
 \ifx#1(\mathchar"4\bffam@28 \else\ifx#1)\mathchar"5\bffam@29 \else
 \ifx#1+\mathchar"2\bffam@2B \else\ifx#1:\mathchar"3\bffam@3A \else
 \ifx#1;\mathchar"6\bffam@3B \else\ifx#1=\mathchar"3\bffam@3D \else
 \ifx#1?\mathchar"5\bffam@3F \else\ifx#1[\mathchar"4\bffam@5B \else
 \ifx#1]\mathchar"5\bffam@5D \else
 \ifx#1,\mathchari@63B \else
 \ifx#1-\mathcharii@200 \else
 \ifx#1.\mathchari@03A \else
 \ifx#1/\mathchari@03D \else
 \ifx#1<\mathchari@33C \else
 \ifx#1>\mathchari@33E \else
 \ifx#1*\mathcharii@203 \else
 \ifx#1|\mathcharii@06A \else
 \ifx#10\bold0\else\ifx#11\bold1\else\ifx#12\bold2\else\ifx#13\bold3\else
 \ifx#14\bold4\else\ifx#15\bold5\else\ifx#16\bold6\else\ifx#17\bold7\else
 \ifx#18\bold8\else\ifx#19\bold9\else
  \Err@{\string\boldkey\space can't be used with #1}%
 \fi\fi\fi\fi\fi\fi\fi\fi\fi\fi\fi\fi\fi\fi\fi
 \fi\fi\fi\fi\fi\fi\fi\fi\fi\fi\fi\fi\fi\fi}
\def\boldsymbol#1{%
 \DN@{\Err@{You can't use \string\boldsymbol\space with \string#1}#1}%
 \ifcat\noexpand#1A%
   \let\next@\relax
  \ifcmmibloaded@{\fam\cmmibfam#1}\else\Err@{First bold symbol
   font not loaded}\fi
 \else
  \xdef\meaning@{\meaning#1.........}%
  \expandafter\math@\meaning@\math@
  \ifmath@
   \expandafter\mathch@\meaning@\mathch@
   \ifmathch@
    \expandafter\boldsymbol@@\meaning@\boldsymbol@@
   \fi
  \else
   \expandafter\macro@\meaning@\macro@
   \expandafter\delim@\meaning@\delim@
   \ifdelim@
    \expandafter\delim@@\meaning@\delim@@
   \else
    \boldsymbol@{#1}%
   \fi
  \fi
 \fi
 \next@}
\def\mathhexboxii@#1#2{\ifcmbsyloaded@\mathhexbox@{\cmbsyfam@}{#1}{#2}\else
  \Err@{Second bold symbol font not loaded}\fi}
\def\boldsymbol@#1{\let\next@\relax\let\next#1%
 \ifx\next\cdot\mathcharii@201 \else
 \ifx\next\prime{{\null\mathcharii@030 \null}}\else
 \ifx\next\lbrack\mathchar"4\bffam@5B \else
 \ifx\next\rbrack\mathchar"5\bffam@5D \else
 \ifx\next\{\mathcharii@466 \else
 \ifx\next\lbrace\mathcharii@466 \else
 \ifx\next\}\mathcharii@567 \else
 \ifx\next\rbrace\mathcharii@567 \else
 \ifx\next\surd{{\mathcharii@170}}\else
 \ifx\next\S{{\mathhexboxii@78}}\else
 \ifx\next\P{{\mathhexboxii@7B}}\else
 \ifx\next\dag{{\mathhexboxii@79}}\else
 \ifx\next\ddag{{\mathhexboxii@7A}}\else
 \DN@{\Err@{You can't use \string\boldsymbol\space with \string#1}#1}%
 \fi\fi\fi\fi\fi\fi\fi\fi\fi\fi\fi\fi\fi}
\def\boldsymbol@@#1.#2\boldsymbol@@{\classnum@#1 \count@@@\classnum@        
 \divide\classnum@4096 \count@\classnum@                                    
 \multiply\count@4096 \advance\count@@@-\count@ \count@@\count@@@           
 \divide\count@@@\@cclvi \count@\count@@                                    
 \multiply\count@@@\@cclvi \advance\count@@-\count@@@                       
 \divide\count@@@\@cclvi                                                    
 \multiply\classnum@4096 \advance\classnum@\count@@                         
 \ifnum\count@@@=\z@                                                        
  \count@"\bffam@ \multiply\count@\@cclvi
  \advance\classnum@\count@
  \DN@{\mathchar\number\classnum@}%
 \else
  \ifnum\count@@@=\@ne                                                      
   \ifcmmibloaded@
   \count@"\cmmibfam@ \multiply\count@\@cclvi
   \advance\classnum@\count@
   \DN@{\mathchar\number\classnum@}%
   \else\DN@{\Err@{First bold symbol font not loaded}}\fi
  \else
   \ifnum\count@@@=\tw@                                                    
  \ifcmbsyloaded@
    \count@"\cmbsyfam@ \multiply\count@\@cclvi
    \advance\classnum@\count@
    \DN@{\mathchar\number\classnum@}%
  \else\DN@{\Err@{Second bold symbol font not loaded}}\fi
  \fi
 \fi
\fi}
\newif\ifdelim@
\newcount\delimcount@
{\uccode`6=`\\ \uccode`7=`d \uccode`8=`e \uccode`9=`l
 \uppercase{\gdef\delim@#1#2#3#4#5\delim@
  {\delim@false\ifx 6#1\ifx 7#2\ifx 8#3\ifx 9#4\delim@true
   \xdef\meaning@{#5}\fi\fi\fi\fi}}}
\def\delim@@#1"#2#3#4#5#6\delim@@{\if#32%
\let\next@\relax
 \ifcmbsyloaded@
 \mathcharii@#2#4#5 \else\Err@{Second bold family not loaded}\fi\fi}
\def\vert{\delimiter"026A30C }
\def\Vert{\delimiter"026B30D }
\let\|\Vert

\def\boldkeydots@#1{\bold@true\let\next=#1\let\delayed@=#1\mdots@@
 \boldkey#1\bold@false}  
\def\boldsymboldots@#1{\bold@true\let\next#1\let\delayed@#1\mdots@@
 \boldsymbol#1\bold@false}
\newif\ifeufbloaded@
\def\loadeufb{\eufbloaded@true
 \font@\teneufb=eufb10
 \font@\seveneufb=eufb7
 \font@\fiveeufb=eufb5
 \alloc@@8\fam\chardef\sixt@@n\eufbfam
 \textfont\eufbfam=\teneufb
 \scriptfont\eufbfam=\seveneufb
 \scriptscriptfont\eufbfam=\fiveeufb
 \edef\eufbfam@{\hexnumber@\eufbfam}}
\newif\ifeusmloaded@
\def\loadeusm{\eusmloaded@true
 \font@\teneusm=eusm10
 \font@\seveneusm=eusm7
 \font@\fiveeusm=eusm5
 \alloc@@8\fam\chardef\sixt@@n\eusmfam
 \textfont\eusmfam=\teneusm
 \scriptfont\eusmfam=\seveneusm
 \scriptscriptfont\eusmfam=\fiveeusm
 \edef\eusmfam@{\hexnumber@\eusmfam}}
\newif\ifeusbloaded@
\def\loadeusb{\eusbloaded@true
 \font@\teneusb=eusb10
 \font@\seveneusb=eusb7
 \font@\fiveeusb=eusb5
 \alloc@@8\fam\chardef\sixt@@n\eusbfam
 \textfont\eusbfam=\teneusb
 \scriptfont\eusbfam=\seveneusb
 \scriptscriptfont\eusbfam=\fiveeusb
 \edef\eusbfam@{\hexnumber@\eusbfam}}
\newif\ifeurmloaded@
\def\loadeurm{\eurmloaded@true
 \font@\teneurm=eurm10
 \font@\seveneurm=eurm7
 \font@\fiveeurm=eurm5
 \alloc@@8\fam\chardef\sixt@@n\eurmfam
 \textfont\eurmfam=\teneurm
 \scriptfont\eurmfam=\seveneurm
 \scriptscriptfont\eurmfam=\fiveeurm
 \edef\eurmfam@{\hexnumber@\eurmfam}}
\newif\ifeurbloaded@
\def\loadeurb{\eurbloaded@true
 \font@\teneurb=eurb10
 \font@\seveneurb=eurb7
 \font@\fiveeurb=eurb5
 \alloc@@8\fam\chardef\sixt@@n\eurbfam
 \textfont\eurbfam=\teneurb
 \scriptfont\eurbfam=\seveneurb
 \scriptscriptfont\eurbfam=\fiveeurb
 \edef\eurbfam@{\hexnumber@\eurbfam}}
\def\accentclass@{7}
\def\noaccents@{\def\accentclass@{0}}
\def\makeacc@#1#2{\def#1{\mathaccent"\accentclass@#2 }}
\makeacc@\hat{05E}
\makeacc@\check{014}
\makeacc@\tilde{07E}
\makeacc@\acute{013}
\makeacc@\grave{012}
\makeacc@\dot{05F}
\makeacc@\ddot{07F}
\makeacc@\breve{015}
\makeacc@\bar{016}

\newcount\skewcharcount@
\newcount\familycount@
\def\theskewchar@{\familycount@\@ne
 \global\skewcharcount@\the\skewchar\textfont\@ne                           
 \ifnum\fam>\m@ne\ifnum\fam<16
  \global\familycount@\the\fam\relax
  \global\skewcharcount@\the\skewchar\textfont\the\fam\relax\fi\fi          
 \ifnum\skewcharcount@>\m@ne
  \ifnum\skewcharcount@<128
  \multiply\familycount@256
  \global\advance\skewcharcount@\familycount@
  \global\advance\skewcharcount@28672
  \mathchar\skewcharcount@\else
  \global\skewcharcount@\m@ne\fi\else
 \global\skewcharcount@\m@ne\fi}                                            
\newcount\pointcount@
\def\getpoints@#1.#2\getpoints@{\pointcount@#1 }
\newdimen\accentdimen@
\newcount\accentmu@
\def\dimentomu@{\multiply\accentdimen@ 100
 \expandafter\getpoints@\the\accentdimen@\getpoints@
 \multiply\pointcount@18
 \divide\pointcount@\@m
 \global\accentmu@\pointcount@}
\def\Makeacc@#1#2{\def#1{\RIfM@\DN@{\mathaccent@
 {"\accentclass@#2 }}\else\DN@{\nonmatherr@{#1}}\fi\next@}}
\def\unbracefonts@{\let\Cal@\Cal@@\let\roman@\roman@@\let\bold@\bold@@
 \let\slanted@\slanted@@}
\def\mathaccent@#1#2{\ifnum\fam=\m@ne\xdef\thefam@{1}\else
 \xdef\thefam@{\the\fam}\fi                                                 
 \accentdimen@\z@                                                           
 \setboxz@h{\unbracefonts@$\m@th\fam\thefam@\relax#2$}
 \ifdim\accentdimen@=\z@\DN@{\mathaccent#1{#2}}
  \setbox@ne\hbox{\unbracefonts@$\m@th\fam\thefam@\relax#2\theskewchar@$}
  \setbox\tw@\hbox{$\m@th\ifnum\skewcharcount@=\m@ne\else
   \mathchar\skewcharcount@\fi$}
  \global\accentdimen@\wd@ne\global\advance\accentdimen@-\wdz@
  \global\advance\accentdimen@-\wd\tw@                                     
  \global\multiply\accentdimen@\tw@
  \dimentomu@\global\advance\accentmu@\@ne                                 
 \else\DN@{{\mathaccent#1{#2\mkern\accentmu@ mu}%
    \mkern-\accentmu@ mu}{}}\fi                                             
 \next@}\Makeacc@\Hat{05E}
\Makeacc@\Check{014}
\Makeacc@\Tilde{07E}
\Makeacc@\Acute{013}
\Makeacc@\Grave{012}
\Makeacc@\Dot{05F}
\Makeacc@\Ddot{07F}
\Makeacc@\Breve{015}
\Makeacc@\Bar{016}
\def\Vec{\RIfM@\DN@{\mathaccent@{"017E }}\else
 \DN@{\nonmatherr@\Vec}\fi\next@}
\def\newbox@{\alloc@4\box\chardef\insc@unt}
\def\accentedsymbol#1#2{\expandafter\newbox@\csname\expandafter
  \eat@\string#1@box\endcsname
 \expandafter\setbox\csname\expandafter\eat@
  \string#1@box\endcsname\hbox{$\m@th#2$}\define
  #1{\expandafter\copy\csname\expandafter\eat@\string#1@box\endcsname{}}}
\def\sqrt#1{\radical"270370 {#1}}
\let\underline@\underline
\let\overline@\overline
\def\underline#1{\underline@{#1}}
\def\overline#1{\overline@{#1}}
\Invalid@\leftroot
\Invalid@\uproot
\newcount\uproot@
\newcount\leftroot@
\def\root{\relaxnext@
  \DN@{\ifx\next\uproot\let\next@\nextii@\else
   \ifx\next\leftroot\let\next@\nextiii@\else
   \let\next@\plainroot@\fi\fi\next@}%
  \DNii@\uproot##1{\uproot@##1\relax\FN@\nextiv@}%
  \def\nextiv@{\ifx\next\space@\DN@. {\FN@\nextv@}\else
   \DN@.{\FN@\nextv@}\fi\next@.}%
  \def\nextv@{\ifx\next\leftroot\let\next@\nextvi@\else
   \let\next@\plainroot@\fi\next@}%
  \def\nextvi@\leftroot##1{\leftroot@##1\relax\plainroot@}%
   \def\nextiii@\leftroot##1{\leftroot@##1\relax\FN@\nextvii@}%
  \def\nextvii@{\ifx\next\space@
   \DN@. {\FN@\nextviii@}\else
   \DN@.{\FN@\nextviii@}\fi\next@.}%
  \def\nextviii@{\ifx\next\uproot\let\next@\nextix@\else
   \let\next@\plainroot@\fi\next@}%
  \def\nextix@\uproot##1{\uproot@##1\relax\plainroot@}%
  \bgroup\uproot@\z@\leftroot@\z@\FN@\next@}
\def\plainroot@#1\of#2{\setbox\rootbox\hbox{$\m@th\scriptscriptstyle{#1}$}%
 \mathchoice{\r@@t\displaystyle{#2}}{\r@@t\textstyle{#2}}
 {\r@@t\scriptstyle{#2}}{\r@@t\scriptscriptstyle{#2}}\egroup}
\def\r@@t#1#2{\setboxz@h{$\m@th#1\sqrt{#2}$}%
 \dimen@\ht\z@\advance\dimen@-\dp\z@
 \setbox@ne\hbox{$\m@th#1\mskip\uproot@ mu$}\advance\dimen@ by1.667\wd@ne
 \mkern-\leftroot@ mu\mkern5mu\raise.6\dimen@\copy\rootbox
 \mkern-10mu\mkern\leftroot@ mu\boxz@}
\def\boxed#1{\setboxz@h{$\m@th\displaystyle{#1}$}\dimen@.4\ex@
 \advance\dimen@3\ex@\advance\dimen@\dp\z@
 \hbox{\lower\dimen@\hbox{%
 \vbox{\hrule height.4\ex@
 \hbox{\vrule width.4\ex@\hskip3\ex@\vbox{\vskip3\ex@\boxz@\vskip3\ex@}%
 \hskip3\ex@\vrule width.4\ex@}\hrule height.4\ex@}%
 }}}
\let\ampersand@\relax
\newdimen\minaw@
\minaw@11.11128\ex@
\newdimen\minCDaw@
\minCDaw@2.5pc
\def\minCDarrowwidth#1{\RIfMIfI@\onlydmatherr@\minCDarrowwidth
 \else\minCDaw@#1\relax\fi\else\onlydmatherr@\minCDarrowwidth\fi}
\newif\ifCD@
\def\CD{\bgroup\vspace@\relax\let\ampersand@&\iffalse}\fi
 \CD@true\vcenter\bgroup\Let@\tabskip\z@skip\baselineskip20\ex@
 \lineskip3\ex@\lineskiplimit3\ex@\halign\bgroup
 &\hfill$\m@th##$\hfill\crcr}
\def\endCD{\crcr\egroup\egroup\egroup}
\newdimen\bigaw@
\atdef@>#1>#2>{\ampersand@                                                  
 \setboxz@h{$\m@th\ssize\;{#1}\;\;$}
 \setbox@ne\hbox{$\m@th\ssize\;{#2}\;\;$}
 \setbox\tw@\hbox{$\m@th#2$}
 \ifCD@\global\bigaw@\minCDaw@\else\global\bigaw@\minaw@\fi                 
 \ifdim\wdz@>\bigaw@\global\bigaw@\wdz@\fi
 \ifdim\wd@ne>\bigaw@\global\bigaw@\wd@ne\fi                                
 \ifCD@\hskip.5em\fi                                                        
 \ifdim\wd\tw@>\z@
  \mathrel{\mathop{\hbox to\bigaw@{\rightarrowfill}}\limits^{#1}_{#2}}
 \else\mathrel{\mathop{\hbox to\bigaw@{\rightarrowfill}}\limits^{#1}}\fi    
 \ifCD@\hskip.5em\fi                                                       
 \ampersand@}                                                              
\atdef@<#1<#2<{\ampersand@\setboxz@h{$\m@th\ssize\;\;{#1}\;$}%
 \setbox@ne\hbox{$\m@th\ssize\;\;{#2}\;$}\setbox\tw@\hbox{$\m@th#2$}%
 \ifCD@\global\bigaw@\minCDaw@\else\global\bigaw@\minaw@\fi
 \ifdim\wdz@>\bigaw@\global\bigaw@\wdz@\fi
 \ifdim\wd@ne>\bigaw@\global\bigaw@\wd@ne\fi
 \ifCD@\hskip.5em\fi
 \ifdim\wd\tw@>\z@
  \mathrel{\mathop{\hbox to\bigaw@{\leftarrowfill}}\limits^{#1}_{#2}}\else
  \mathrel{\mathop{\hbox to\bigaw@{\leftarrowfill}}\limits^{#1}}\fi
 \ifCD@\hskip.5em\fi\ampersand@}
\atdef@)#1)#2){\ampersand@
 \setboxz@h{$\m@th\ssize\;{#1}\;\;$}%
 \setbox@ne\hbox{$\m@th\ssize\;{#2}\;\;$}%
 \setbox\tw@\hbox{$\m@th#2$}%
 \ifCD@
 \global\bigaw@\minCDaw@\else\global\bigaw@\minaw@\fi
 \ifdim\wdz@>\bigaw@\global\bigaw@\wdz@\fi
 \ifdim\wd@ne>\bigaw@\global\bigaw@\wd@ne\fi
 \ifCD@\hskip.5em\fi
 \ifdim\wd\tw@>\z@
  \mathrel{\mathop{\hbox to\bigaw@{\rightarrowfill}}\limits^{#1}_{#2}}%
 \else\mathrel{\mathop{\hbox to\bigaw@{\rightarrowfill}}\limits^{#1}}\fi
 \ifCD@\hskip.5em\fi
 \ampersand@}
\atdef@(#1(#2({\ampersand@\setboxz@h{$\m@th\ssize\;\;{#1}\;$}%
 \setbox@ne\hbox{$\m@th\ssize\;\;{#2}\;$}\setbox\tw@\hbox{$\m@th#2$}%
 \ifCD@\global\bigaw@\minCDaw@\else\global\bigaw@\minaw@\fi
 \ifdim\wdz@>\bigaw@\global\bigaw@\wdz@\fi
 \ifdim\wd@ne>\bigaw@\global\bigaw@\wd@ne\fi
 \ifCD@\hskip.5em\fi
 \ifdim\wd\tw@>\z@
  \mathrel{\mathop{\hbox to\bigaw@{\leftarrowfill}}\limits^{#1}_{#2}}\else
  \mathrel{\mathop{\hbox to\bigaw@{\leftarrowfill}}\limits^{#1}}\fi
 \ifCD@\hskip.5em\fi\ampersand@}
\atdef@ A#1A#2A{\llap{$\m@th\vcenter{\hbox
 {$\ssize#1$}}$}\Big\uparrow\rlap{$\m@th\vcenter{\hbox{$\ssize#2$}}$}&&}
\atdef@ V#1V#2V{\llap{$\m@th\vcenter{\hbox
 {$\ssize#1$}}$}\Big\downarrow\rlap{$\m@th\vcenter{\hbox{$\ssize#2$}}$}&&}
\atdef@={&\hskip.5em\mathrel
 {\vbox{\hrule width\minCDaw@\vskip3\ex@\hrule width
 \minCDaw@}}\hskip.5em&}
\atdef@|{\Big\Vert&&}
\atdef@@\vert{\Big\Vert&&}
\def\pretend#1\haswidth#2{\setboxz@h{$\m@th\scriptstyle{#2}$}\hbox
 to\wdz@{\hfill$\m@th\scriptstyle{#1}$\hfill}}
\def\pmb{\RIfM@\expandafter\mathpalette\expandafter\pmb@\else
 \expandafter\pmb@@\fi}
\def\pmb@@#1{\leavevmode\setboxz@h{#1}\kern-.025em\copy\z@\kern-\wdz@
 \kern-.05em\copy\z@\kern-\wdz@\kern-.025em\raise.0433em\boxz@}
\def\binrel@@#1{\ifdim\wd2<\z@\mathbin{#1}\else\ifdim\wd\tw@>\z@
 \mathrel{#1}\else{#1}\fi\fi}
\newdimen\pmbraise@
\def\pmb@#1#2{\setbox\thr@@\hbox{$\m@th#1{#2}$}%
 \setbox4 \hbox{$\m@th#1\mkern.7794mu$}\pmbraise@\wd4
 \binrel@{#2}\binrel@@{\mkern-.45mu\copy\thr@@\kern-\wd\thr@@
 \mkern-.9mu\copy\thr@@\kern-\wd\thr@@\mkern-.45mu\raise\pmbraise@\box\thr@@}}
\def\documentstyle#1{\input #1.sty\relax}
\font\dummyft@=dummy
\fontdimen1 \dummyft@=\z@
\fontdimen2 \dummyft@=\z@
\fontdimen3 \dummyft@=\z@
\fontdimen4 \dummyft@=\z@
\fontdimen5 \dummyft@=\z@
\fontdimen6 \dummyft@=\z@
\fontdimen7 \dummyft@=\z@
\fontdimen8 \dummyft@=\z@
\fontdimen9 \dummyft@=\z@
\fontdimen10 \dummyft@=\z@
\fontdimen11 \dummyft@=\z@
\fontdimen12 \dummyft@=\z@
\fontdimen13 \dummyft@=\z@
\fontdimen14 \dummyft@=\z@
\fontdimen15 \dummyft@=\z@
\fontdimen16 \dummyft@=\z@
\fontdimen17 \dummyft@=\z@
\fontdimen18 \dummyft@=\z@
\fontdimen19 \dummyft@=\z@
\fontdimen20 \dummyft@=\z@
\fontdimen21 \dummyft@=\z@
\fontdimen22 \dummyft@=\z@
\def\fontlist@{\\{\tenrm}\\{\sevenrm}\\{\fiverm}\\{\teni}\\{\seveni}%
 \\{\fivei}\\{\tensy}\\{\sevensy}\\{\fivesy}\\{\tenex}\\{\tenbf}\\{\sevenbf}%
 \\{\fivebf}\\{\tensl}\\{\tenit}}
\def\font@#1=#2 {\rightappend@#1\to\fontlist@\font#1=#2 }
\def\dodummy@{{\def\\##1{\global\let##1\dummyft@}\fontlist@}}
\def\nopages@{\output={\setbox\z@\box255 \deadcycles\z@}%
 \alloc@5\toks\toksdef\@cclvi\output}
\let\galleys\nopages@
\newif\ifsyntax@
\newcount\countxviii@
\def\syntax{\syntax@true\dodummy@\countxviii@\count18
 \loop\ifnum\countxviii@>\m@ne\textfont\countxviii@=\dummyft@
 \scriptfont\countxviii@=\dummyft@\scriptscriptfont\countxviii@=\dummyft@
 \advance\countxviii@\m@ne\repeat                                           
 \dummyft@\tracinglostchars\z@\nopages@\frenchspacing\hbadness\@M}
\def\S@{S } \def\G@{G } \def\P@{P }
\newif\ifbadans@
\def\printoptions{\W@{Do you want S(yntax check),
  G(alleys) or P(ages)?^^JType S, G or P, follow by <return>: }\loop
 \read\m@ne to\ans@
 \xdef\next@{\def\noexpand\Ans@{\ans@}}\uppercase\expandafter{\next@}
 \ifx\Ans@\S@\badans@false\syntax\else
 \ifx\Ans@\G@\badans@false\galleys\else
 \ifx\Ans@\P@\badans@false\else
 \badans@true\fi\fi\fi
 \ifbadans@\W@{Type S, G or P, follow by <return>: }%
 \repeat}
\def\alloc@#1#2#3#4#5{\global\advance\count1#1by\@ne
 \ch@ck#1#4#2\allocationnumber=\count1#1
 \global#3#5=\allocationnumber
 \ifalloc@\wlog{\string#5=\string#2\the\allocationnumber}\fi}
\def\document{\def\alloclist@{}\def\fontlist@{}}
\let\enddocument\bye

\let\proclaim\undefined
\let\footnote\undefined
\let\=\undefined
\let\>\undefined

\catcode`\@=\active

%
%
\def\next{AMSPPT}\ifx\styname\next \endinput\fi
\catcode`\@=11
\def\styname{AMSPPT}
\def\styversion{2.0}
{\W@{\styname.STY - Version \styversion}\W@{}}
\hyphenation{acad-e-my acad-e-mies af-ter-thought anom-aly anom-alies
an-ti-deriv-a-tive an-tin-o-my an-tin-o-mies apoth-e-o-ses apoth-e-o-sis
ap-pen-dix ar-che-typ-al as-sign-a-ble as-sist-ant-ship as-ymp-tot-ic
asyn-chro-nous at-trib-uted at-trib-ut-able bank-rupt bank-rupt-cy
bi-dif-fer-en-tial blue-print busier busiest cat-a-stroph-ic
cat-a-stroph-i-cally con-gress cross-hatched data-base de-fin-i-tive
de-riv-a-tive dis-trib-ute dri-ver dri-vers eco-nom-ics econ-o-mist
elit-ist equi-vari-ant ex-quis-ite ex-tra-or-di-nary flow-chart
for-mi-da-ble forth-right friv-o-lous ge-o-des-ic ge-o-det-ic geo-met-ric
griev-ance griev-ous griev-ous-ly hexa-dec-i-mal ho-lo-no-my ho-mo-thetic
ideals idio-syn-crasy in-fin-ite-ly in-fin-i-tes-i-mal ir-rev-o-ca-ble
key-stroke lam-en-ta-ble light-weight mal-a-prop-ism man-u-script
mar-gin-al meta-bol-ic me-tab-o-lism meta-lan-guage me-trop-o-lis
met-ro-pol-i-tan mi-nut-est mol-e-cule mono-chrome mono-pole mo-nop-oly
mono-spline mo-not-o-nous mul-ti-fac-eted mul-ti-plic-able non-euclid-ean
non-iso-mor-phic non-smooth par-a-digm par-a-bol-ic pa-rab-o-loid
pa-ram-e-trize para-mount pen-ta-gon phe-nom-e-non post-script pre-am-ble
pro-ce-dur-al pro-hib-i-tive pro-hib-i-tive-ly pseu-do-dif-fer-en-tial
pseu-do-fi-nite pseu-do-nym qua-drat-ics quad-ra-ture qua-si-smooth
qua-si-sta-tion-ary qua-si-tri-an-gu-lar quin-tes-sence quin-tes-sen-tial
re-arrange-ment rec-tan-gle ret-ri-bu-tion retro-fit retro-fit-ted
right-eous right-eous-ness ro-bot ro-bot-ics sched-ul-ing se-mes-ter
semi-def-i-nite semi-ho-mo-thet-ic set-up se-vere-ly side-step sov-er-eign
spe-cious spher-oid spher-oid-al star-tling star-tling-ly
sta-tis-tics sto-chas-tic straight-est strange-ness strat-a-gem strong-hold
sum-ma-ble symp-to-matic syn-chro-nous topo-graph-i-cal tra-vers-a-ble
tra-ver-sal tra-ver-sals treach-ery turn-around un-at-tached un-err-ing-ly
white-space wide-spread wing-spread wretch-ed wretch-ed-ly Brown-ian
Eng-lish Euler-ian Feb-ru-ary Gauss-ian Grothen-dieck Hamil-ton-ian
Her-mit-ian Jan-u-ary Japan-ese Kor-te-weg Le-gendre Lip-schitz
Lip-schitz-ian Mar-kov-ian Noe-ther-ian No-vem-ber Rie-mann-ian
Schwarz-schild Sep-tem-ber}
\Invalid@\nofrills
\Invalid@\usualspace
\newif\ifnofrills@
\def\nofrills@#1#2{\relaxnext@
  \DN@{\ifx\next\nofrills
    \nofrills@true\let#2\relax\DN@\nofrills{\nextii@}%
  \else
    \nofrills@false\def#2{#1}\let\next@\nextii@\fi
\next@}}
\def\usualspace@#1{\ifnofrills@\def\usualspace{#1}\fi}
\def\addto#1#2{\csname \expandafter\eat@\string#1@\endcsname
  \expandafter{\the\csname \expandafter\eat@\string#1@\endcsname#2}}
\newdimen\bigsize@
\def\big@#1#2{{\hbox{$\left#2\vcenter to#1\bigsize@{}%
  \right.\nulldelimiterspace\z@\m@th$}}}
\def\big{\big@\@ne}
\def\Big{\big@{1.5}}
\def\bigg{\big@\tw@}
\def\Bigg{\big@{2.5}}
\def\raggedcenter@{\leftskip\z@ plus.4\hsize \rightskip\leftskip
 \parfillskip\z@ \parindent\z@ \spaceskip.3333em \xspaceskip.5em
 \pretolerance9999\tolerance9999 \exhyphenpenalty\@M
 \hyphenpenalty\@M \let\\\linebreak}
\def\upperspecialchars{\def\ss{SS}\let\i=I\let\j=J\let\ae\AE\let\oe\OE
  \let\o\O\let\aa\AA\let\l\L}
\def\uppercasetext@#1{%
  {\spaceskip1.2\fontdimen2\the\font plus1.2\fontdimen3\the\font
   \upperspecialchars\uctext@#1$\m@th\aftergroup\eat@$}}
\def\uctext@#1$#2${\endash@#1-\endash@$#2$\uctext@}
\def\endash@#1-#2\endash@{\uppercase{#1}\if\notempty{#2}--\endash@#2\endash@\fi}
\def\runaway@#1{\DN@{#1}\ifx\envir@\next@
  \Err@{You seem to have a missing or misspelled \string\end#1 ...}%
  \let\envir@\empty\fi}
\newif\iftemp@
\def\notempty#1{TT\fi\def\test@{#1}\ifx\test@\empty\temp@false
  \else\temp@true\fi \iftemp@}
\font@\tensmc=cmcsc10
\font@\sevenex=cmex7
\font@\sevenit=cmti7
\font@\eightrm=cmr8 
\font@\sixrm=cmr6 
\font@\eighti=cmmi8     \skewchar\eighti='177 
\font@\sixi=cmmi6       \skewchar\sixi='177   
\font@\eightsy=cmsy8    \skewchar\eightsy='60 
\font@\sixsy=cmsy6      \skewchar\sixsy='60   
\font@\eightex=cmex8
\font@\eightbf=cmbx8 
\font@\sixbf=cmbx6   
\font@\eightit=cmti8 
\font@\eightsl=cmsl8 
\font@\eightsmc=cmcsc8
\font@\eighttt=cmtt8 
\loadmsam
\loadmsbm
\loadeufm
\input amssym.tex\relax
\newtoks\tenpoint@
\def\tenpoint{\normalbaselineskip12\p@
 \abovedisplayskip12\p@ plus3\p@ minus9\p@
 \belowdisplayskip\abovedisplayskip
 \abovedisplayshortskip\z@ plus3\p@
 \belowdisplayshortskip7\p@ plus3\p@ minus4\p@
 \textonlyfont@\rm\tenrm \textonlyfont@\it\tenit
 \textonlyfont@\sl\tensl \textonlyfont@\bf\tenbf
 \textonlyfont@\smc\tensmc \textonlyfont@\tt\tentt
 \ifsyntax@ \def\big##1{{\hbox{$\left##1\right.$}}}%
  \let\Big\big \let\bigg\big \let\Bigg\big
 \else
  \textfont\z@=\tenrm  \scriptfont\z@=\sevenrm  \scriptscriptfont\z@=\fiverm
  \textfont\@ne=\teni  \scriptfont\@ne=\seveni  \scriptscriptfont\@ne=\fivei
  \textfont\tw@=\tensy \scriptfont\tw@=\sevensy \scriptscriptfont\tw@=\fivesy
  \textfont\thr@@=\tenex \scriptfont\thr@@=\sevenex
        \scriptscriptfont\thr@@=\sevenex
  \textfont\itfam=\tenit \scriptfont\itfam=\sevenit
        \scriptscriptfont\itfam=\sevenit
  \textfont\bffam=\tenbf \scriptfont\bffam=\sevenbf
        \scriptscriptfont\bffam=\fivebf
  \setbox\strutbox\hbox{\vrule height8.5\p@ depth3.5\p@ width\z@}%
  \setbox\strutbox@\hbox{\lower.5\normallineskiplimit\vbox{%
        \kern-\normallineskiplimit\copy\strutbox}}%
 \setbox\z@\vbox{\hbox{$($}\kern\z@}\bigsize@=1.2\ht\z@
 \fi
 \normalbaselines\rm\ex@.2326ex\jot3\ex@\the\tenpoint@}
\newtoks\eightpoint@
\def\eightpoint{\normalbaselineskip10\p@
 \abovedisplayskip10\p@ plus2.4\p@ minus7.2\p@
 \belowdisplayskip\abovedisplayskip
 \abovedisplayshortskip\z@ plus2.4\p@
 \belowdisplayshortskip5.6\p@ plus2.4\p@ minus3.2\p@
 \textonlyfont@\rm\eightrm \textonlyfont@\it\eightit
 \textonlyfont@\sl\eightsl \textonlyfont@\bf\eightbf
 \textonlyfont@\smc\eightsmc \textonlyfont@\tt\eighttt
 \ifsyntax@\def\big##1{{\hbox{$\left##1\right.$}}}%
  \let\Big\big \let\bigg\big \let\Bigg\big
 \else
  \textfont\z@=\eightrm \scriptfont\z@=\sixrm \scriptscriptfont\z@=\fiverm
  \textfont\@ne=\eighti \scriptfont\@ne=\sixi \scriptscriptfont\@ne=\fivei
  \textfont\tw@=\eightsy \scriptfont\tw@=\sixsy \scriptscriptfont\tw@=\fivesy
  \textfont\thr@@=\eightex \scriptfont\thr@@=\sevenex
   \scriptscriptfont\thr@@=\sevenex
  \textfont\itfam=\eightit \scriptfont\itfam=\sevenit
   \scriptscriptfont\itfam=\sevenit
  \textfont\bffam=\eightbf \scriptfont\bffam=\sixbf
   \scriptscriptfont\bffam=\fivebf
 \setbox\strutbox\hbox{\vrule height7\p@ depth3\p@ width\z@}%
 \setbox\strutbox@\hbox{\raise.5\normallineskiplimit\vbox{%
   \kern-\normallineskiplimit\copy\strutbox}}%
 \setbox\z@\vbox{\hbox{$($}\kern\z@}\bigsize@=1.2\ht\z@
 \fi
 \normalbaselines\eightrm\ex@.2326ex\jot3\ex@\the\eightpoint@}
\parindent1pc
\normallineskiplimit\p@
\newdimen\indenti \indenti=2pc
\def\pageheight#1{\vsize#1}
\def\pagewidth#1{\hsize#1%
   \captionwidth@\hsize \advance\captionwidth@-2\indenti}
\pagewidth{30pc} \pageheight{47pc}
\def\topmatter{%
 \ifx\undefined\msafam
 \else\font@\eightmsa=msam8 \font@\sixmsa=msam6
   \ifsyntax@\else \addto\tenpoint{\textfont\msafam=\tenmsa
              \scriptfont\msafam=\sevenmsa \scriptscriptfont\msafam=\fivemsa}%
     \addto\eightpoint{\textfont\msafam=\eightmsa \scriptfont\msafam=\sixmsa
              \scriptscriptfont\msafam=\fivemsa}%
   \fi
 \fi
 \ifx\undefined\msbfam
 \else\font@\eightmsb=msbm8 \font@\sixmsb=msbm6
   \ifsyntax@\else \addto\tenpoint{\textfont\msbfam=\tenmsb
         \scriptfont\msbfam=\sevenmsb \scriptscriptfont\msbfam=\fivemsb}%
     \addto\eightpoint{\textfont\msbfam=\eightmsb \scriptfont\msbfam=\sixmsb
         \scriptscriptfont\msbfam=\fivemsb}%
   \fi
 \fi
 \ifx\undefined\eufmfam
 \else \font@\eighteufm=eufm8 \font@\sixeufm=eufm6
   \ifsyntax@\else \addto\tenpoint{\textfont\eufmfam=\teneufm
       \scriptfont\eufmfam=\seveneufm \scriptscriptfont\eufmfam=\fiveeufm}%
     \addto\eightpoint{\textfont\eufmfam=\eighteufm
       \scriptfont\eufmfam=\sixeufm \scriptscriptfont\eufmfam=\fiveeufm}%
   \fi
 \fi
 \ifx\undefined\eufbfam
 \else \font@\eighteufb=eufb8 \font@\sixeufb=eufb6
   \ifsyntax@\else \addto\tenpoint{\textfont\eufbfam=\teneufb
      \scriptfont\eufbfam=\seveneufb \scriptscriptfont\eufbfam=\fiveeufb}%
    \addto\eightpoint{\textfont\eufbfam=\eighteufb
      \scriptfont\eufbfam=\sixeufb \scriptscriptfont\eufbfam=\fiveeufb}%
   \fi
 \fi
 \ifx\undefined\eusmfam
 \else \font@\eighteusm=eusm8 \font@\sixeusm=eusm6
   \ifsyntax@\else \addto\tenpoint{\textfont\eusmfam=\teneusm
       \scriptfont\eusmfam=\seveneusm \scriptscriptfont\eusmfam=\fiveeusm}%
     \addto\eightpoint{\textfont\eusmfam=\eighteusm
       \scriptfont\eusmfam=\sixeusm \scriptscriptfont\eusmfam=\fiveeusm}%
   \fi
 \fi
 \ifx\undefined\eusbfam
 \else \font@\eighteusb=eusb8 \font@\sixeusb=eusb6
   \ifsyntax@\else \addto\tenpoint{\textfont\eusbfam=\teneusb
       \scriptfont\eusbfam=\seveneusb \scriptscriptfont\eusbfam=\fiveeusb}%
     \addto\eightpoint{\textfont\eusbfam=\eighteusb
       \scriptfont\eusbfam=\sixeusb \scriptscriptfont\eusbfam=\fiveeusb}%
   \fi
 \fi
 \ifx\undefined\eurmfam
 \else \font@\eighteurm=eurm8 \font@\sixeurm=eurm6
   \ifsyntax@\else \addto\tenpoint{\textfont\eurmfam=\teneurm
       \scriptfont\eurmfam=\seveneurm \scriptscriptfont\eurmfam=\fiveeurm}%
     \addto\eightpoint{\textfont\eurmfam=\eighteurm
       \scriptfont\eurmfam=\sixeurm \scriptscriptfont\eurmfam=\fiveeurm}%
   \fi
 \fi
 \ifx\undefined\eurbfam
 \else \font@\eighteurb=eurb8 \font@\sixeurb=eurb6
   \ifsyntax@\else \addto\tenpoint{\textfont\eurbfam=\teneurb
       \scriptfont\eurbfam=\seveneurb \scriptscriptfont\eurbfam=\fiveeurb}%
    \addto\eightpoint{\textfont\eurbfam=\eighteurb
       \scriptfont\eurbfam=\sixeurb \scriptscriptfont\eurbfam=\fiveeurb}%
   \fi
 \fi
 \ifx\undefined\cmmibfam
 \else \font@\eightcmmib=cmmib8 \font@\sixcmmib=cmmib6
   \ifsyntax@\else \addto\tenpoint{\textfont\cmmibfam=\tencmmib
       \scriptfont\cmmibfam=\sevencmmib \scriptscriptfont\cmmibfam=\fivecmmib}%
    \addto\eightpoint{\textfont\cmmibfam=\eightcmmib
       \scriptfont\cmmibfam=\sixcmmib \scriptscriptfont\cmmibfam=\fivecmmib}%
   \fi
 \fi
 \ifx\undefined\cmbsyfam
 \else \font@\eightcmbsy=cmbsy8 \font@\sixcmbsy=cmbsy6
   \ifsyntax@\else \addto\tenpoint{\textfont\cmbsyfam=\tencmbsy
      \scriptfont\cmbsyfam=\sevencmbsy \scriptscriptfont\cmbsyfam=\fivecmbsy}%
    \addto\eightpoint{\textfont\cmbsyfam=\eightcmbsy
      \scriptfont\cmbsyfam=\sixcmbsy \scriptscriptfont\cmbsyfam=\fivecmbsy}%
   \fi
 \fi
 \let\topmatter\relax}
\def\chapterno@{\uppercase\expandafter{\romannumeral\chaptercount@}}
\newcount\chaptercount@
\def\chapter{\nofrills@{\afterassignment\chapterno@
                        CHAPTER \global\chaptercount@=}\chapter@
 \DNii@##1{\leavevmode\hskip-\leftskip
   \rlap{\vbox to\z@{\vss\centerline{\eightpoint
   \chapter@##1\unskip}\baselineskip2pc\null}}\hskip\leftskip
   \nofrills@false}%
 \FN@\next@}
\newbox\titlebox@
\def\title{\nofrills@{\uppercasetext@}\title@%
 \DNii@##1\endtitle{\global\setbox\titlebox@\vtop{\tenpoint\bf
 \raggedcenter@\ignorespaces
 \baselineskip1.3\baselineskip\title@{##1}\endgraf}%
 \ifmonograph@ \edef\next{\the\leftheadtoks}\ifx\next\empty
    \leftheadtext{##1}\fi
 \fi
 \edef\next{\the\rightheadtoks}\ifx\next\empty \rightheadtext{##1}\fi
 }\FN@\next@}
\newbox\authorbox@
\def\author#1\endauthor{\global\setbox\authorbox@
 \vbox{\tenpoint\smc\raggedcenter@\ignorespaces
 #1\endgraf}\relaxnext@ \edef\next{\the\leftheadtoks}%
 \ifx\next\empty\leftheadtext{#1}\fi}
\newbox\affilbox@
\def\affil#1\endaffil{\global\setbox\affilbox@
 \vbox{\tenpoint\raggedcenter@\ignorespaces#1\endgraf}}
\newcount\addresscount@
\addresscount@\z@
\def\address#1\endaddress{\global\advance\addresscount@\@ne
  \expandafter\gdef\csname address\number\addresscount@\endcsname
  {\vskip12\p@ minus6\p@\noindent\eightpoint\smc\ignorespaces#1\par}}
\def\email{\nofrills@{\eightpoint{\it E-mail\/}:\enspace}\email@
  \DNii@##1\endemail{%
  \expandafter\gdef\csname email\number\addresscount@\endcsname
  {\def\usualspace{{\it\enspace}}\smallskip\noindent\eightpoint\email@
  \ignorespaces##1\par}}%
 \FN@\next@}
\def\thedate@{}
\def\date#1\enddate{\gdef\thedate@{\tenpoint\ignorespaces#1\unskip}}
\def\thethanks@{}
\def\thanks#1\endthanks{\gdef\thethanks@{\eightpoint\ignorespaces#1.\unskip}}
\def\thekeywords@{}
\def\keywords{\nofrills@{{\it Key words and phrases.\enspace}}\keywords@
 \DNii@##1\endkeywords{\def\thekeywords@{\def\usualspace{{\it\enspace}}%
 \eightpoint\keywords@\ignorespaces##1\unskip.}}%
 \FN@\next@}
\def\thesubjclass@{}
\def\subjclass{\nofrills@{{\rm1980 {\it Mathematics Subject
   Classification\/} (1985 {\it Revision\/}).\enspace}}\subjclass@
 \DNii@##1\endsubjclass{\def\thesubjclass@{\def\usualspace
  {{\rm\enspace}}\eightpoint\subjclass@\ignorespaces##1\unskip.}}%
 \FN@\next@}
\newbox\abstractbox@
\def\abstract{\nofrills@{{\smc Abstract.\enspace}}\abstract@
 \DNii@{\setbox\abstractbox@\vbox\bgroup\noindent$$\vbox\bgroup
  \def\envir@{abstract}\advance\hsize-2\indenti
  \usualspace@{{\enspace}}\eightpoint \noindent\abstract@\ignorespaces}%
 \FN@\next@}
\def\endabstract{\par\unskip\egroup$$\egroup}
\def\widestnumber#1#2{\begingroup\let\head\null\let\subhead\empty
   \let\subsubhead\subhead
   \ifx#1\head\global\setbox\tocheadbox@\hbox{#2.\enspace}%
   \else\ifx#1\subhead\global\setbox\tocsubheadbox@\hbox{#2.\enspace}%
   \else\ifx#1\key\bgroup\let\endrefitem@\egroup
        \key#2\endrefitem@\global\refindentwd\wd\keybox@
   \else\ifx#1\no\bgroup\let\endrefitem@\egroup
        \no#2\endrefitem@\global\refindentwd\wd\nobox@
   \else\ifx#1\page\global\setbox\pagesbox@\hbox{\quad\bf#2}%
   \else\ifx#1\item\setboxz@h{#2}\global\rosteritemwd\wdz@
        \global\advance\rosteritemwd by.5\parindent
   \else\message{\string\widestnumber is not defined for this option
   (\string#1)}%
\fi\fi\fi\fi\fi\fi\endgroup}
\newif\ifmonograph@
\def\Monograph{\monograph@true \let\headmark\rightheadtext
  \let\varindent@\indent \def\headfont@{\bf}\def\proclaimfont@{\smc}%
  \def\demofont@{\smc}}
\let\varindent@\noindent
\newbox\tocheadbox@    \newbox\tocsubheadbox@
\newbox\tocbox@
\def\toc{\toc@{Contents}}
\def\newtocdefs{%
   \def \title##1\endtitle
       {\penaltyandskip@\z@\smallskipamount
        \hangindent\wd\tocheadbox@\noindent{\bf##1}}%
   \def \chapter##1{%
        Chapter \uppercase\expandafter{\romannumeral##1.\unskip}\enspace}%
   \def \specialhead##1\endspecialhead
       {\par\hangindent\wd\tocheadbox@ \noindent##1\par}%
   \def \head##1 ##2\endhead
       {\par\hangindent\wd\tocheadbox@ \noindent
        \if\notempty{##1}\hbox to\wd\tocheadbox@{\hfil##1\enspace}\fi
        ##2\par}%
   \def \subhead##1 ##2\endsubhead
       {\par\vskip-\parskip {\normalbaselines
        \advance\leftskip\wd\tocheadbox@
        \hangindent\wd\tocsubheadbox@ \noindent
        \if\notempty{##1}\hbox to\wd\tocsubheadbox@{##1\unskip\hfil}\fi
         ##2\par}}%
   \def \subsubhead##1 ##2\endsubsubhead
       {\par\vskip-\parskip {\normalbaselines
        \advance\leftskip\wd\tocheadbox@
        \hangindent\wd\tocsubheadbox@ \noindent
        \if\notempty{##1}\hbox to\wd\tocsubheadbox@{##1\unskip\hfil}\fi
        ##2\par}}}
\def\toc@#1{\relaxnext@
   \def\page##1%
       {\unskip\penalty0\null\hfil
        \rlap{\hbox to\wd\pagesbox@{\quad\hfil##1}}\hfilneg\penalty\@M}%
 \DN@{\ifx\next\nofrills\DN@\nofrills{\nextii@}%
      \else\DN@{\nextii@{{#1}}}\fi
      \next@}%
 \DNii@##1{%
\ifmonograph@\bgroup\else\setbox\tocbox@\vbox\bgroup
   \centerline{\headfont@\ignorespaces##1\unskip}\nobreak
   \vskip\belowheadskip \fi
   \setbox\tocheadbox@\hbox{0.\enspace}%
   \setbox\tocsubheadbox@\hbox{0.0.\enspace}%
   \leftskip\indenti \rightskip\leftskip
   \setbox\pagesbox@\hbox{\bf\quad000}\advance\rightskip\wd\pagesbox@
   \newtocdefs
 }%
 \FN@\next@}
\def\endtoc{\par\egroup}
\let\pretitle\relax
\let\preauthor\relax
\let\preaffil\relax
\let\predate\relax
\let\preabstract\relax
\let\prepaper\relax
\def\dedicatory #1\enddedicatory{\def\preabstract{{\medskip
  \eightpoint\it \raggedcenter@#1\endgraf}}}
\def\thetranslator@{}
\def\translator#1\endtranslator{\def\thetranslator@{\nobreak\medskip
 \line{\eightpoint\hfil Translated by \uppercase{#1}\qquad\qquad}\nobreak}}
\outer\def\endtopmatter{\runaway@{abstract}%
 \edef\next{\the\leftheadtoks}\ifx\next\empty
  \expandafter\leftheadtext\expandafter{\the\rightheadtoks}\fi
 \ifmonograph@\else
   \ifx\thesubjclass@\empty\else \makefootnote@{}{\thesubjclass@}\fi
   \ifx\thekeywords@\empty\else \makefootnote@{}{\thekeywords@}\fi
   \ifx\thethanks@\empty\else \makefootnote@{}{\thethanks@}\fi
 \fi
  \pretitle
  \ifmonograph@ \topskip7pc \else \topskip4pc \fi
  \box\titlebox@
  \topskip10pt
  \preauthor
  \ifvoid\authorbox@\else \vskip2.5pc plus1pc \unvbox\authorbox@\fi
  \preaffil
  \ifvoid\affilbox@\else \vskip1pc plus.5pc \unvbox\affilbox@\fi
  \predate
  \ifx\thedate@\empty\else \vskip1pc plus.5pc \line{\hfil\thedate@\hfil}\fi
  \preabstract
  \ifvoid\abstractbox@\else \vskip1.5pc plus.5pc \unvbox\abstractbox@ \fi
  \ifvoid\tocbox@\else\vskip1.5pc plus.5pc \unvbox\tocbox@\fi
  \prepaper
  \vskip2pc plus1pc
}
\def\document{\let\fontlist@\relax\let\alloclist@\relax
  \tenpoint}
\newskip\aboveheadskip       \aboveheadskip\bigskipamount
\newdimen\belowheadskip      \belowheadskip6\p@
\def\headfont@{\smc}
\def\penaltyandskip@#1#2{\relax\ifdim\lastskip<#2\relax\removelastskip
      \ifnum#1=\z@\else\penalty@#1\relax\fi\vskip#2%
  \else\ifnum#1=\z@\else\penalty@#1\relax\fi\fi}
\def\nobreak{\penalty\@M
  \ifvmode\def\penalty@{\let\penalty@\penalty\count@@@}%
  \everypar{\let\penalty@\penalty\everypar{}}\fi}
\let\penalty@\penalty
\def\heading#1\endheading{\head#1\endhead}
\def\subheading#1{\subhead#1\endsubhead}
\def\specialheadfont@{\bf}
\outer\def\specialhead{\par\penaltyandskip@{-200}\aboveheadskip
  \begingroup\interlinepenalty\@M\rightskip\z@ plus\hsize \let\\\linebreak
  \specialheadfont@\noindent\ignorespaces}
\def\endspecialhead{\par\endgroup\nobreak\vskip\belowheadskip}
\outer\def\head#1\endhead{\par\penaltyandskip@{-200}\aboveheadskip
 {\headfont@\raggedcenter@\interlinepenalty\@M
 \ignorespaces#1\endgraf}\nobreak
 \vskip\belowheadskip
 \headmark{#1}}
\let\headmark\eat@
\newskip\subheadskip       \subheadskip\medskipamount
\def\subheadfont@{\bf}
\outer\def\subhead{\nofrills@{.\enspace}\subhead@
 \DNii@##1\endsubhead{\par\penaltyandskip@{-100}\subheadskip
  \varindent@{\usualspace@{{\subheadfont@\enspace}}%
 \subheadfont@\ignorespaces##1\unskip\subhead@}\ignorespaces}%
 \FN@\next@}
\outer\def\subsubhead{\nofrills@{.\enspace}\subsubhead@
 \DNii@##1\endsubsubhead{\par\penaltyandskip@{-50}\medskipamount
      {\usualspace@{{\it\enspace}}%
  \it\ignorespaces##1\unskip\subsubhead@}\ignorespaces}%
 \FN@\next@}
\def\proclaimheadfont@{\bf}
\outer\def\proclaim{\runaway@{proclaim}\def\envir@{proclaim}%
  \nofrills@{.\enspace}\proclaim@
 \DNii@##1{\penaltyandskip@{-100}\medskipamount\varindent@
   \usualspace@{{\proclaimheadfont@\enspace}}\proclaimheadfont@
   \ignorespaces##1\unskip\proclaim@
  \sl\ignorespaces}%
 \FN@\next@}
\outer\def\endproclaim{\let\envir@\relax\par\rm
  \penaltyandskip@{55}\medskipamount}
\def\demoheadfont@{\it}
\def\demo{\runaway@{proclaim}\nofrills@{.\enspace}\demo@
     \DNii@##1{\par\penaltyandskip@\z@\medskipamount
  {\usualspace@{{\demoheadfont@\enspace}}%
  \varindent@\demoheadfont@\ignorespaces##1\unskip\demo@}\rm
  \ignorespaces}\FN@\next@}

\def\qed{\ifhmode\unskip\nobreak\fi\quad\ifmmode\square\else$\m@th\square$\fi}

\def\definition{\runaway@{proclaim}%
  \nofrills@{.\proclaimheadfont@\enspace}\definition@
        \DNii@##1{\penaltyandskip@{-100}\medskipamount
        {\usualspace@{{\proclaimheadfont@\enspace}}%
        \varindent@\proclaimheadfont@\ignorespaces##1\unskip\definition@}%
        \rm \ignorespaces}\FN@\next@}

\newdimen\rosteritemwd
\newcount\rostercount@
\newif\iffirstitem@
\let\plainitem@\item
\newtoks\everypartoks@
\def\par@{\everypartoks@\expandafter{\the\everypar}\everypar{}}
\def\roster{\edef\leftskip@{\leftskip\the\leftskip}%
 \relaxnext@
 \rostercount@\z@  
 \def\item{\FN@\rosteritem@}%
 \DN@{\ifx\next\runinitem\let\next@\nextii@\else
  \let\next@\nextiii@\fi\next@}%
 \DNii@\runinitem  
  {\unskip  
   \DN@{\ifx\next[\let\next@\nextii@\else
    \ifx\next"\let\next@\nextiii@\else\let\next@\nextiv@\fi\fi\next@}%
   \DNii@[####1]{\rostercount@####1\relax
    \enspace{\rm(\number\rostercount@)}~\ignorespaces}%
   \def\nextiii@"####1"{\enspace{\rm####1}~\ignorespaces}%
   \def\nextiv@{\enspace{\rm(1)}\rostercount@\@ne~}%
   \par@\firstitem@false  
   \FN@\next@}%
 \def\nextiii@{\par\par@  
  \penalty\@m\smallskip\vskip-\parskip
  \firstitem@true}%
 \FN@\next@}
\def\rosteritem@{\iffirstitem@\firstitem@false\else\par\vskip-\parskip\fi
 \leftskip3\parindent\noindent  
 \DNii@[##1]{\rostercount@##1\relax
  \llap{\hbox to2.5\parindent{\hss\rm(\number\rostercount@)}%
   \hskip.5\parindent}\ignorespaces}%
 \def\nextiii@"##1"{%
  \llap{\hbox to2.5\parindent{\hss\rm##1}\hskip.5\parindent}\ignorespaces}%
 \def\nextiv@{\advance\rostercount@\@ne
  \llap{\hbox to2.5\parindent{\hss\rm(\number\rostercount@)}%
   \hskip.5\parindent}}%
 \ifx\next[\let\next@\nextii@\else\ifx\next"\let\next@\nextiii@\else
  \let\next@\nextiv@\fi\fi\next@}

\newif\ifnextRunin@
\def\endroster{\relaxnext@
 \par\leftskip@  
 \penalty-50 \vskip-\parskip\smallskip  
 \DN@{\ifx\next\Runinitem\let\next@\relax
  \else\nextRunin@false\let\item\plainitem@  
   \ifx\next\par 
    \DN@\par{\everypar\expandafter{\the\everypartoks@}}%
   \else  
    \DN@{\noindent\everypar\expandafter{\the\everypartoks@}}%
  \fi\fi\next@}%
 \FN@\next@}
\newcount\rosterhangafter@
\def\Runinitem#1\roster\runinitem{\relaxnext@
 \rostercount@\z@ 
 \def\item{\FN@\rosteritem@}%
 \def\runinitem@{#1}%
 \DN@{\ifx\next[\let\next\nextii@\else\ifx\next"\let\next\nextiii@
  \else\let\next\nextiv@\fi\fi\next}%
 \DNii@[##1]{\rostercount@##1\relax
  \def\item@{{\rm(\number\rostercount@)}}\nextv@}%
 \def\nextiii@"##1"{\def\item@{{\rm##1}}\nextv@}%
 \def\nextiv@{\advance\rostercount@\@ne
  \def\item@{{\rm(\number\rostercount@)}}\nextv@}%
 \def\nextv@{\setbox\z@\vbox  
  {\ifnextRunin@\noindent\fi  
  \runinitem@\unskip\enspace\item@~\par  
  \global\rosterhangafter@\prevgraf}%
  \firstitem@false  
  \ifnextRunin@\else\par\fi  
  \hangafter\rosterhangafter@\hangindent3\parindent
  \ifnextRunin@\noindent\fi  
  \runinitem@\unskip\enspace 
  \item@~\ifnextRunin@\else\par@\fi  
  \nextRunin@true\ignorespaces}%
 \FN@\next@}
\def\footmarkform@#1{$\m@th^{#1}$}
\let\thefootnotemark\footmarkform@
\def\makefootnote@#1#2{\insert\footins
 {\interlinepenalty\interfootnotelinepenalty
 \eightpoint\splittopskip\ht\strutbox\splitmaxdepth\dp\strutbox
 \floatingpenalty\@MM\leftskip\z@\rightskip\z@\spaceskip\z@\xspaceskip\z@
 \leavevmode{#1}\footstrut\ignorespaces#2\unskip\lower\dp\strutbox
 \vbox to\dp\strutbox{}}}
\newcount\footmarkcount@
\footmarkcount@\z@
\def\footnotemark{\let\@sf\empty\relaxnext@
 \ifhmode\edef\@sf{\spacefactor\the\spacefactor}\/\fi
 \DN@{\ifx[\next\let\next@\nextii@\else
  \ifx"\next\let\next@\nextiii@\else
  \let\next@\nextiv@\fi\fi\next@}%
 \DNii@[##1]{\footmarkform@{##1}\@sf}%
 \def\nextiii@"##1"{{##1}\@sf}%
 \def\nextiv@{\iffirstchoice@\global\advance\footmarkcount@\@ne\fi
  \footmarkform@{\number\footmarkcount@}\@sf}%
 \FN@\next@}
\def\footnotetext{\relaxnext@
 \DN@{\ifx[\next\let\next@\nextii@\else
  \ifx"\next\let\next@\nextiii@\else
  \let\next@\nextiv@\fi\fi\next@}%
 \DNii@[##1]##2{\makefootnote@{\footmarkform@{##1}}{##2}}%
 \def\nextiii@"##1"##2{\makefootnote@{##1}{##2}}%
 \def\nextiv@##1{\makefootnote@{\footmarkform@{\number\footmarkcount@}}{##1}}%
 \FN@\next@}
\def\footnote{\let\@sf\empty\relaxnext@
 \ifhmode\edef\@sf{\spacefactor\the\spacefactor}\/\fi
 \DN@{\ifx[\next\let\next@\nextii@\else
  \ifx"\next\let\next@\nextiii@\else
  \let\next@\nextiv@\fi\fi\next@}%
 \DNii@[##1]##2{\footnotemark[##1]\footnotetext[##1]{##2}}%
 \def\nextiii@"##1"##2{\footnotemark"##1"\footnotetext"##1"{##2}}%
 \def\nextiv@##1{\footnotemark\footnotetext{##1}}%
 \FN@\next@}
\def\adjustfootnotemark#1{\advance\footmarkcount@#1\relax}
\def\footnoterule{\kern-3\p@
  \hrule width 5pc\kern 2.6\p@} 
\def\captionfont@{\smc}
\def\topcaption#1#2\endcaption{%
  {\dimen@\hsize \advance\dimen@-\captionwidth@
   \rm\raggedcenter@ \advance\leftskip.5\dimen@ \rightskip\leftskip
  {\captionfont@#1}%
  \if\notempty{#2}.\enspace\ignorespaces#2\fi
  \endgraf}\nobreak\bigskip}
\def\botcaption#1#2\endcaption{%
  \nobreak\bigskip
  \setboxz@h{\captionfont@#1\if\notempty{#2}.\enspace\rm#2\fi}%
  {\dimen@\hsize \advance\dimen@-\captionwidth@
   \leftskip.5\dimen@ \rightskip\leftskip
   \noindent \ifdim\wdz@>\captionwidth@ 
   \else\hfil\fi 
  {\captionfont@#1}\if\notempty{#2}.\enspace\rm#2\fi\endgraf}}
\def\@ins{\par\begingroup\def\vspace##1{\vskip##1\relax}%
  \def\captionwidth##1{\captionwidth@##1\relax}%
  \setbox\z@\vbox\bgroup} 
\def\block{\RIfMIfI@\nondmatherr@\block\fi
       \else\ifvmode\vskip\abovedisplayskip\noindent\fi
        $$\def\endblock{\par\egroup$$}\fi
  \vbox\bgroup\advance\hsize-2\indenti\noindent}
\def\endblock{\par\egroup}
\def\cite#1{{\rm[{\citefont@\m@th#1}]}}
\def\citefont@{\rm}
\def\refsfont@{\eightpoint}
\outer\def\Refs{\runaway@{proclaim}%
 \relaxnext@ \DN@{\ifx\next\nofrills\DN@\nofrills{\nextii@}\else
  \DN@{\nextii@{References}}\fi\next@}%
 \DNii@##1{\penaltyandskip@{-200}\aboveheadskip
  \line{\hfil\headfont@\ignorespaces##1\unskip\hfil}\nobreak
  \vskip\belowheadskip
  \begingroup\refsfont@\sfcode`.=\@m}%
 \FN@\next@}

\newbox\nobox@            \newbox\keybox@           \newbox\bybox@
\newbox\paperbox@         \newbox\paperinfobox@     \newbox\jourbox@
\newbox\volbox@           \newbox\issuebox@         \newbox\yrbox@
\newbox\pagesbox@         \newbox\bookbox@          \newbox\bookinfobox@
\newbox\publbox@          \newbox\publaddrbox@      \newbox\finalinfobox@
\newbox\edsbox@           \newbox\langbox@
\newif\iffirstref@        \newif\iflastref@
\newif\ifprevjour@        \newif\ifbook@            \newif\ifprevinbook@
\newif\ifquotes@          \newif\ifbookquotes@      \newif\ifpaperquotes@
\newdimen\bysamerulewd@
\setboxz@h{\refsfont@\kern3em}
\bysamerulewd@\wdz@
\newdimen\refindentwd
\setboxz@h{\refsfont@ 00. }
\refindentwd\wdz@
\outer\def\ref{\begingroup \noindent\hangindent\refindentwd
 \firstref@true \def\nofrills{\def\refkern@{\kern3sp}}%
 \ref@}
\def\ref@{\book@false \bgroup\let\endrefitem@\egroup \ignorespaces}
\def\moreref{\endrefitem@\endref@\firstref@false\ref@}%
\def\transl{\endrefitem@\endref@\firstref@false
  \book@false
  \prepunct@
  \setboxz@h\bgroup \aftergroup\unhbox\aftergroup\z@
    \def\endrefitem@{\unskip\refkern@\egroup}\ignorespaces}%
\def\emptyifempty@{\dimen@\wd\currbox@
  \advance\dimen@-\wd\z@ \advance\dimen@-.1\p@
  \ifdim\dimen@<\z@ \setbox\currbox@\copy\voidb@x \fi}
\let\refkern@\relax
\def\endrefitem@{\unskip\refkern@\egroup
  \setboxz@h{\refkern@}\emptyifempty@}\ignorespaces
\def\refdef@#1#2#3{\edef\next@{\noexpand\endrefitem@
  \let\noexpand\currbox@\csname\expandafter\eat@\string#1box@\endcsname
    \noexpand\setbox\noexpand\currbox@\hbox\bgroup}%
  \toks@\expandafter{\next@}%
  \if\notempty{#2#3}\toks@\expandafter{\the\toks@
  \def\endrefitem@{\unskip#3\refkern@\egroup
  \setboxz@h{#2#3\refkern@}\emptyifempty@}#2}\fi
  \toks@\expandafter{\the\toks@\ignorespaces}%
  \edef#1{\the\toks@}}
\refdef@\no{}{. }
\refdef@\key{[\m@th}{] }
\refdef@\by{}{}
\def\bysame{\by\hbox to\bysamerulewd@{\hrulefill}\thinspace
   \kern0sp}
\def\manyby{\message{\string\manyby is no longer necessary; \string\by
  can be used instead, starting with version 2.0 of \styname.STY}\by}
\refdef@\paper{\ifpaperquotes@``\fi\it}{}
\refdef@\paperinfo{}{}
\def\jour{\endrefitem@\let\currbox@\jourbox@
  \setbox\currbox@\hbox\bgroup
  \def\endrefitem@{\unskip\refkern@\egroup
    \setboxz@h{\refkern@}\emptyifempty@
    \ifvoid\jourbox@\else\prevjour@true\fi}%
\ignorespaces}
\refdef@\vol{\ifbook@\else\bf\fi}{}
\refdef@\issue{no. }{}
\refdef@\yr{}{}
\refdef@\pages{}{}
\def\page{\endrefitem@\def\pp@{\def\pp@{pp.~}p.~}\let\currbox@\pagesbox@
  \setbox\currbox@\hbox\bgroup\ignorespaces}
\def\pp@{pp.~}
\def\book{\endrefitem@ \let\currbox@\bookbox@
 \setbox\currbox@\hbox\bgroup\def\endrefitem@{\unskip\refkern@\egroup
  \setboxz@h{\ifbookquotes@``\fi}\emptyifempty@
  \ifvoid\bookbox@\else\book@true\fi}%
  \ifbookquotes@``\fi\it\ignorespaces}
\def\inbook{\endrefitem@
  \let\currbox@\bookbox@\setbox\currbox@\hbox\bgroup
  \def\endrefitem@{\unskip\refkern@\egroup
  \setboxz@h{\ifbookquotes@``\fi}\emptyifempty@
  \ifvoid\bookbox@\else\book@true\previnbook@true\fi}%
  \ifbookquotes@``\fi\ignorespaces}
\refdef@\eds{(}{, eds.)}
\def\ed{\endrefitem@\let\currbox@\edsbox@
 \setbox\currbox@\hbox\bgroup
 \def\endrefitem@{\unskip, ed.)\refkern@\egroup
  \setboxz@h{(, ed.)}\emptyifempty@}(\ignorespaces}
\refdef@\bookinfo{}{}
\refdef@\publ{}{}
\refdef@\publaddr{}{}
\refdef@\finalinfo{}{}
\refdef@\lang{(}{)}

\let\refdef@\relax 
\def\ppunbox@#1{\ifvoid#1\else\prepunct@\unhbox#1\fi}
\def\nocomma@#1{\ifvoid#1\else\changepunct@3\prepunct@\unhbox#1\fi}
\def\changepunct@#1{\ifnum\lastkern<3 \unkern\kern#1sp\fi}
\def\prepunct@{\count@\lastkern\unkern
  \ifnum\lastpenalty=0
    \let\penalty@\relax
  \else
    \edef\penalty@{\penalty\the\lastpenalty\relax}%
  \fi
  \unpenalty
  \let\refspace@\ \ifcase\count@,
\or;\or.\or 
  \or\let\refspace@\relax
  \else,\fi
  \ifquotes@''\quotes@false\fi \penalty@ \refspace@
}
\def\transferpenalty@#1{\dimen@\lastkern\unkern
  \ifnum\lastpenalty=0\unpenalty\let\penalty@\relax
  \else\edef\penalty@{\penalty\the\lastpenalty\relax}\unpenalty\fi
  #1\penalty@\kern\dimen@}
\def\endref{\endrefitem@\lastref@true\endref@
  \par\endgroup \prevjour@false \previnbook@false }
\def\endref@{%
\iffirstref@
  \ifvoid\nobox@\ifvoid\keybox@\indent\fi
  \else\hbox to\refindentwd{\hss\unhbox\nobox@}\fi
  \ifvoid\keybox@
  \else\ifdim\wd\keybox@>\refindentwd
         \box\keybox@
       \else\hbox to\refindentwd{\unhbox\keybox@\hfil}\fi\fi
  \kern4sp\ppunbox@\bybox@
\fi 
  \ifvoid\paperbox@
  \else\prepunct@\unhbox\paperbox@
    \ifpaperquotes@\quotes@true\fi\fi
  \ppunbox@\paperinfobox@
  \ifvoid\jourbox@
    \ifprevjour@ \nocomma@\volbox@
      \nocomma@\issuebox@
      \ifvoid\yrbox@\else\changepunct@3\prepunct@(\unhbox\yrbox@
        \transferpenalty@)\fi
      \ppunbox@\pagesbox@
    \fi 
  \else \prepunct@\unhbox\jourbox@
    \nocomma@\volbox@
    \nocomma@\issuebox@
    \ifvoid\yrbox@\else\changepunct@3\prepunct@(\unhbox\yrbox@
      \transferpenalty@)\fi
    \ppunbox@\pagesbox@
  \fi 
  \ifbook@\prepunct@\unhbox\bookbox@ \ifbookquotes@\quotes@true\fi \fi
  \nocomma@\edsbox@
  \ppunbox@\bookinfobox@
  \ifbook@\ifvoid\volbox@\else\prepunct@ vol.~\unhbox\volbox@
  \fi\fi
  \ppunbox@\publbox@ \ppunbox@\publaddrbox@
  \ifbook@ \ppunbox@\yrbox@
    \ifvoid\pagesbox@
    \else\prepunct@\pp@\unhbox\pagesbox@\fi
  \else
    \ifprevinbook@ \ppunbox@\yrbox@
      \ifvoid\pagesbox@\else\prepunct@\pp@\unhbox\pagesbox@\fi
    \fi \fi
  \ppunbox@\finalinfobox@
  \iflastref@
    \ifvoid\langbox@.\ifquotes@''\fi
    \else\changepunct@2\prepunct@\unhbox\langbox@\fi
  \else
    \ifvoid\langbox@\changepunct@1%
    \else\changepunct@3\prepunct@\unhbox\langbox@
      \changepunct@1\fi
  \fi
}
\outer\def\enddocument{%
 \runaway@{proclaim}%
\ifmonograph@ 
\else
 \nobreak
 \thetranslator@
 \count@\z@ \loop\ifnum\count@<\addresscount@\advance\count@\@ne
 \csname address\number\count@\endcsname
 \csname email\number\count@\endcsname
 \repeat
\fi
 \vfill\supereject\end}
\def\folio{{\foliofont@\ifnum\pageno<\z@ \romannumeral-\pageno
 \else\number\pageno \fi}}
\def\foliofont@{\eightrm}
\def\headlinefont@{\eightpoint}
\def\leftheadline{\rlap{\folio}\hfill \iftrue\topmark\fi \hfill}
\def\rightheadline{\hfill \expandafter
  \hfill \llap{\folio}}
\newtoks\leftheadtoks
\newtoks\rightheadtoks
\def\leftheadtext{\nofrills@{\uppercasetext@}\lht@
  \DNii@##1{\leftheadtoks\expandafter{\lht@{##1}}%
    \mark{\the\leftheadtoks\noexpand\else\the\rightheadtoks}
    \ifsyntax@\setboxz@h{\def\\{\unskip\space\ignorespaces}%
        \headlinefont@##1}\fi}%
  \FN@\next@}
\def\rightheadtext{\nofrills@{\uppercasetext@}\rht@
  \DNii@##1{\rightheadtoks\expandafter{\rht@{##1}}%
    \mark{\the\leftheadtoks\noexpand\else\the\rightheadtoks}%
    \ifsyntax@\setboxz@h{\def\\{\unskip\space\ignorespaces}%
        \headlinefont@##1}\fi}%
  \FN@\next@}
\headline={\def\chapter#1{\chapterno@. }%
  \def\\{\unskip\space\ignorespaces}\headlinefont@
  \ifodd\pageno \rightheadline \else \leftheadline\fi}
\def\NoRunningHeads{\global\runheads@false\global\let\headmark\eat@}

\def\logo@{\baselineskip2pc \hbox to\hsize{\hfil\eightpoint Typeset by
 \AmSTeX}}
\newif\iffirstpage@     \firstpage@true
\newif\ifrunheads@      \runheads@true
\output={\output@}
\def\output@{\shipout\vbox{%
 \iffirstpage@ \global\firstpage@false
  \pagebody \logo@ \makefootline%
 \else \ifrunheads@ \makeheadline \pagebody
       \else \pagebody \makefootline \fi
 \fi}%
 \advancepageno \ifnum\outputpenalty>-\@MM\else\dosupereject\fi}
\tenpoint
\catcode`\@=\active

 \magnification=1200
 \NoBlackBoxes





 \define\zz{ {{\bold{Z}_2}}}


 \define\calc{\Cal C}

 \define\calq{\Cal Q}
 
 \define\calz{\Cal Z}



 \define\cycd#1#2{{\calz}_{#1}(#2)}
 \define\cyc#1#2{{\calz}^{#1}(#2)}
 \define\cych#1{{{\calz}^{#1}}}
 \define\cycp#1#2{{\calz}^{#1}(\bbp(#2))}
 \define\cyf#1#2{{\cyc{#1}{#2}}^{fix}}
 \define\cyfd#1#2{{\cycd{#1}{#2}}^{fix}}

 \define\crl#1#2{{\calz}_{\bbr}^{#1}{(#2)}}

 \define\crd#1#2{\widetilde{\calz}_{\bbr}^{#1}{(#2)}}

 \define\crld#1#2{{\calz}_{#1,\bbr}{(#2)}}
 \define\crdd#1#2{\widetilde{\calz}_{#1,{\bbr}}{(#2)}}

 \define\crlh#1{{\calz}_{\bbr}^{#1}}
 \define\crdh#1{{\widetilde{\calz}_{\bbr}^{#1}}}
 \define\cyav#1#2{{{\cyc{#1}{#2}}^{av}}}
 \define\cyavd#1#2{{\cycd{#1}{#2}}^{av}}

 \define\cyaa#1#2{{\cyc{#1}{#2}}^{-}}
 \define\cyaad#1#2{{\cycd{#1}{#2}}^{-}}

 \define\cyq#1#2{{\calq}^{#1}(#2)}
 \define\cyqd#1#2{{\calq}_{#1}(#2)}

 \define\cqt#1#2{{\calz}_{\bbh}^{#1}{(#2)}}
 \define\cqtav#1#2{{\calz}^{#1}{(#2)}^{av}}
 \define\cqtrd#1#2{\widetilde{\calz}_{\bbh}^{#1}{(#2)}}

 \define\cyct#1#2{{\calz}^{#1}(#2)_\zz}
 \define\cyft#1#2{{\cyc{#1}{#2}}^{fix}_\zz}
\define\cxg#1#2{G^{#1}_{\bbc}(#2)}
 \define\reg#1#2{G^{#1}_{\bbr}(#2)}

 \define\cyaat#1#2{{\cyc{#1}{#2}}^{-}_\zz}


 \define\fflag#1#2{{#1}={#1}_{#2} \supset {#1}_{{#2}-1} \supset
 \ldots \supset {#1}_{0} }
 
 \define\vect#1{ {\Cal{V}ect}_{#1}}


 \define\chv#1#2#3{{\calc}^{#1}_{#2}(#3)}
 \define\chvd#1#2#3{{\calc}_{#1,#2}(#3)}
 \define\chm#1#2{{\calc}_{#1}(#2)}



 \define\Claim#1{\subheading{Claim #1}}

 \define\xrightarrow#1{\overset{#1}\to{\rightarrow}}


\hfuzz1pc 


\define\bbz{\Bbb Z}

\define\bbr{\Bbb R}
\define\bbc{\Bbb C}
\define\bbh{\Bbb H}
\define\bbp{\Bbb P}

\define\cc{\Cal C}
\define\cd{\Cal D}
\define\ce{\Cal E}

\define\cf{\Cal F}

\define\cs{\Cal S}
\define\cz{\Cal Z}
\define\co#1{\Cal O_{#1}}
\define\ct{\Cal T}

\define\cR{\Cal R}


\define\a{\alpha}
\redefine\b{\beta}

\redefine\d{\delta}
\define\r{\rho}

\define\z{\zeta}

\define\e{\epsilon}
\redefine\D{\Delta}

\define\p#1{{\bbp}^{#1}}


\define\blbx{\hfill  $\square$}
\redefine\qed{\blbx}

\define\pf{\subheading{Proof}}
\define\Lemma#1{\subheading{Lemma #1}}
\define\Theorem#1{\subheading{Theorem #1}}
\define\Prop#1{\subheading{Proposition #1}}
\define\Cor#1{\subheading{Corollary #1}}
\define\Note#1{\subheading{Note #1}}
\define\Def#1{\subheading{Definition #1}}
\define\Remark#1{\subheading{Remark #1}}
\define\Ex#1{\subheading{Example #1}}
\define\arr{\longrightarrow}

\define\Hom{\text{Hom}}

\redefine\Xi{X_{\infty}}

\define\jac#1#2{\left(\!\!\!\left(
\frac{\partial #1}{\partial #2}
\right)\!\!\!\right)}
\define\restrict#1{\left. #1 \right|_{t_{p+1} = \dots = t_n = 0}}

\define\SP#1#2{{\roman SP}^{#1}(#2)}

\define\coc#1#2#3{\cc^{#1}(#2;\, #3)}
\define\zoc#1#2#3{\cz^{#1}(#2;\, #3)}
\define\zyc#1#2#3{\cz^{#1}(#2 \times #3)}

\define\ar#1{\overset{#1}\to{\longrightarrow}}

\define\th#1#2{{\Bbb H}^{#1}(#2)}
\define\hth#1#2{\widehat{\Bbb H}^{#1}(#2)}


\define\bad#1#2{\cf_{#2}(#1)}

\define\pch#1{\bbp_{\bbc}(\bbh^{#1})}

\def\l{\ell}

\def\I#1#2{I_{#1, #2}}

\def\D{D}

\def\L{L}

\def\z2t{\text{$\bbz_2\ct$}}

\def\hH{\widehat{I\!\!H}}

\def\oper{\operatorname}

\def\Lloc{L^1_{\oper{loc}}}

\def\cpt{\oper{cpt}}

\def\<{\left<}
\def\>{\right>}
\def\[{\left[}
\def\]{\right]}

\def\wt{\widetilde}

\def\supp{\operatorname{supp \ }}

\redefine\and{\qquad\text{and}\qquad}

\document



\ \bigskip
\centerline{\bf \titfont  Lefschetz-Pontrjagin Duality for
Differential Characters }

\vskip .2in
\centerline{by}
\vskip .3in
\def\bx{\partial X}
\def\xbx{(X,\bx)}
\def\du{^{\star\!\!\! \star}}
\def\sdu{^{{\star\!\!\! \star}_{\infty}}}

\centerline{\bf  \aufont   Reese Harvey \  and \
 Blaine Lawson\footnote{\footfont Research
of all authors was partially supported by the NSF. Research of the
second author was also partially supported by IHES and CMI.}}

\vskip .7in

\centerline{\bf Abstract} \medskip
  \font\abstractfont=cmr10 at 9 pt

{{\parindent=.7in\narrower\abstractfont \noindent 
A theory of differential characters is developed for manifolds with
boundary.  This is done from both the Cheeger-Simons and the 
deRham-Federer viewpoints. The central  result of the paper is the
formulation and proof of a Lefschetz-Pontrjagin Duality Theorem, which
asserts that the pairing
$$
\hH^k\xbx  \times \hH^{n-k-1}(X) \ \arr\ S^1
$$
given by
$$
(\a,\b)\ \mapsto\ (\a *\b)\,[X]
$$
induces isomorphisms
$$\aligned
\cd:\hH^k\xbx&\to \Hom_{\infty}(\hH^{n-k-1}(X),\, S^1)
\\
\cd':\hH^{n-k-1}(X)&\to \Hom_{\infty}(\hH^k\xbx,\, S^1)
\endaligned
$$
onto the smooth Pontrjagin duals.  In particular, $\cd$ and $\cd'$ are
injective with dense range in the group of all continuous
homomorphisms into the circle. A coboundary map is introduced which
yields a long sequence for the character groups associated to the pair
$\xbx$.  The relation of the sequence to the duality mappings is
analyzed.

 }}

\vfill\eject   
               
\centerline{Table of Contents}
\medskip  
       
\hskip 1in  0.  Introduction
 
\hskip 1in  1.  Differential Characters on Manifolds with Boundary    
  
\hskip 1in  2.  The Exact Sequences  
  
\hskip 1in  3.  The Star Product 

\hskip 1in  4.  Smooth Pontrjagin duals

\hskip 1in  5.  Lefschetz-Pontrjagin Duality

\hskip 1in  6.  Coboundary Maps
 
\hskip 1in  7.  Sequences and  Duality

\vskip .3in

\subheading{\S 0. \ \ Introduction}

The  theory of differential characters, introduced
by  Jim Simons and Jeff Cheeger in 1973, is of basic importance in
geometry.  It provides a wealth of invariants for bundles with
connection   starting with the classical one of  Chern-Simons 
 in dimension 3 and including large families of invariants for 
flat bundles and foliations.   Its cardinal
property is that it forms the natural receiving space for a refined
Chern-Weil theory.  This theory subsumes  integral characteristic
classes and the classical Chern-Weil characteristic forms.  It also
tracks certain ``transgression terms'' which give cohomologies between
smooth and singular cocycles and lead to interesting secondary
invariants.

Each  standard characteristic classes has a refinement in the group of
differential characters.  Thus for a complex bundle with unitary
connection, refined Chern classes ${\widehat c}_k$ are
defined and the total class 
$$\widehat c = 1 +  {\widehat c}_1 +
{\widehat c}_2+\dots: {I\!\!  K}(X) \arr
\hH^*(X)    
$$
gives a natural
transformation  from the $K$-theory of bundles {\sl with connection} to 
differential characters which satisfies the Whitney sum formula:
$\widehat c (E\oplus F) = \widehat c(E) * \widehat c(F)$.  This last
property leads to non-conformal immersion theorems in riemannian
geometry.

Differential characters form a highly structured theory with
certain aspects of cohomology: contravariant functoriality, ring
structure, and a pairing to cycles.  There are deRham-Federer
formulations of the theory [GS$_1$], [H], [HLZ], analogous to those
given for cohomology, which are useful for example in the theory of
singular connections [HL$_1$],[HL$_2$].  Furthermore, the groups
$\hH^k(X)$ of differential characters carry a natural topology.  The
connected component of 0 in this group consists of the {\bf smooth
characters},
 those which can be represented by smooth differential forms.

In [HLZ], where the deRham-Federer appoach is developed in detail,
the authors showed that differential characters satisfy {\bf 
Poincar\'e-Pontjagin duality}: On an oriented $n$ dimensional manifold
$X$ the pairing
$$
\hH^k(X)\times \hH^{n-k-1}_{\cpt}(X) \ \arr\ S^1
$$
given by 
$$
(\a,\b)\ \mapsto\ (\a *\b)[X]
$$
(where $\hH_{\cpt}^*$ denotes characters with compact support) induces
injective maps 
 $$
\hH^k(X)\ \to\ \Hom\left(\hH^{n-k-1}_{\cpt}(X),
S^1\right)
\and
\hH^{n-k-1}_{\cpt}(X)\ \to\ \Hom\left(\hH^k(X),
S^1\right)
 $$
with dense range in the groups of continuous homomorphisms into the
circle.  Moreover this range consists exactly of 
the {\bf smooth} homomorphisms.  These are defined precisely in \S 4
but can be thought of roughly as follows. The connected component 
of 0 in $\hH^k(X)$ consists essentially (i.e., up to a
finite-dimensional torus factor) of the exact $(k+1)$-forms
$d\ce^{k+1}(X)$ with the $C^{\infty}$-topology. Now
$\Hom(d\ce^{k+1}(X),S^1) = \Hom(d\ce^{k+1}(X),\bbr)$ is just the vector
space dual.  This is simply a quotient of the space of {\sl currents},
the $(n-k-1)$-forms with distribution coefficients.  The smooth
dual corresponds to those forms which have smooth coefficients.

In this paper we formulate the theory of differential characters for
compact manifolds with boundary $\xbx$ and prove a Lefschetz-Pontrjagin 
Duality Theorem analogous to the one above. To do this we introduce the
relative groups $\hH^*\xbx$ and develop the theory from [HLZ] for
this case.  The main theorem asserts the  existence of a pairing
$$
\hH^k(X)\times \hH^{n-k-1}\xbx \ \arr\ S^1
$$
given by 
$
(\a,\b)\ \mapsto\ (\a *\b)[X]
$
and inducing injective maps with dense range as above.

The two pairings above  have a formal similarity but are far
from  the same. The delicate part of these dualities comes from the
differential form component of characters.  In the first pairing     
(on possibly non-compact manifolds) we contrast  forms having no growth
restrictions at infinity with forms with compact support. The second
dualtiy (on compact manfiolds with boundary) opposes forms smooth up to
the bounday with forms which restrict to zero on the boundary.

In cohomology theory there are long exact sequences for the pair
$\xbx$ which interlace the Pontrjagin and Lefschetz Duality mappings.
In the last sections of this paper the parallel structure for
differential characters is studied. We introduce coboundary maps
$\partial:\hH^k(X)\to\hH^{k+1}\xbx$, yielding long sequences which
intertwine the duality mappings and reduce to the standard picture
under the natural transformation to integral cohomology.

\vskip .3in

\vfill\eject

\subheading{\bf \S 1.  Differential characters on manifolds with
boundary}
Let $X$ be  a compact oriented differentiable
$n$-manifold with boundary $\bx$.  Let $\ce^*(X)$ denote the de Rham
complex of differential forms which are smooth up to the boundary,  
and set 
$$
\ce^*(X,\bx)\ =\ \{\phi\in\ce^*(X)\ :\ \phi\bigl|_{\bx}=0\}.
$$
The cohomology of this complex is naturally isomorphic to
$H^*(X,\bx;\,\bbr)$.  Let $C_*(X)$ denote the complex of
$C^{\infty}$-singular chains on $X$ and  
$C_*(X,\bx)\ \equiv\  C_*(X)/C_*(\bx)$  the relative complex.
Denote by 
$$
Z_*(X,\bx)\ \equiv \ \{c\in C_*(X,\bx)\ :\ \partial c=0\}
$$
the cycles in this complex. We begin with definitions  of
differential characters in the spirit of Cheeger-Simons.

\Def{1.1}  The group of    {\bf differential characters} of degree
$k$ on $X$ is the set of homomorphisms
$$
\hat{H}^k(X;\bbr/\bbz) \ \equiv\ \{\a\in \Hom(Z_k(X), S^1)\ :\ \d(\a)\in
\ce^{k+1}(X)\} 
$$
where $\d$ denotes the coboundary.  Similarly the group of  {\bf
relative differential characters} of degree $k$ on $\xbx$ is defined
to be
$$
\hat{H}^k(X, \bx;\bbr/\bbz)\ \equiv \ \{\a\in \Hom(Z_k\xbx, S^1)\ :\
\d(\a)\in \ce^{k+1}\xbx\} 
$$
Inclusion and restriction give maps
$
\hH^k\xbx @>{j}>> \hH^k(X) @>{\r}>> \hH^k(\bx).
$
with $\r\circ j=0$.
\medskip

There is an alternative  de Rham-Federer approach to these groups.
Set
$$\aligned
\ce_{\Lloc}^k(X) \ &\equiv\ 
\text{ $k$-forms on $X$ with $\Lloc$-coefficients}\\
\cR^k(X) 
\ \equiv\ \text{ the }&\text{ rectifiable currents
of degree $k$ (dimension $n-k$) on $X$} \\
\ce_{\Lloc}^k\xbx \ \equiv\ \{a\in &\ce_{\Lloc}^k(X)\ :\ a
\text{ in smooth in a neighborhood of $\bx$ and  } a\bigl|_{\bx}=0\}\\
\cR^k_{\cpt}(X-\bx)
\ &\equiv\ \{R\in \cR^k(X)\ :\ \text{supp}(R)\subset X-\bx\}
\endaligned
$$

\Def{1.2} An element $a\in\ce_{\Lloc}^k(X)$ is called a {\bf
  spark} of degree $k$ on $X$ if
$$
da\ =\ \phi - R \qquad \text{where} \ \ \phi\in\ce^{k+1}(X)\text{ and
}  R\in \cR^{k+1}(X). 
\tag1.3
$$  
Denote by  $\cs^k(X)$  the group of all such sparks and by
$\ct^k(X)$  the subgroup of all $a\in \cs^k(X)$ such that $a=db +S$
where $b\in \ce^{k-1}_{\Lloc}(X)$ $S\in  \cR^{k}_{\cpt}(X)$. Then  the
group of { de Rham-Federer characters} of degree $k$ on $X$ is defined
 to be the quotient 
$$
\hH^k(X)\ \equiv\ \cs^k(X) / \ct^k(X).
$$
We define {\bf relative sparks} and {\bf relative de
Rham-Federer characters } on $\xbx$ by 
$$\aligned
\cs^k\xbx \equiv &\{a\in \ce_{\Lloc}^k\xbx\ :\ 
da=\phi-r, \ \ \phi\in\ce^{k+1}\xbx\text{ and
}  R\in \cR^{k+1}_{\cpt}(X-\bx) \}    \\
\ct^k\xbx &\equiv \{a\in\cs^k\xbx \ :\ a=db +S, \ \ b\in 
\ce^{k-1}_{\Lloc}\xbx \text{ and } S\in \cR^{k}_{\cpt}(X-\bx) \}    \\
&\ \qquad\hH^k\xbx\ \equiv\ \cs^k\xbx / \ct^k\xbx.
\endaligned
$$
\eject

The decomposition (1.3) is unique. In fact we have the following.
Recall that a current
$T$ is said to be {\bf integrally flat} if it can be written as
$T=R+dS$ where $R$ and $S$ are rectifiable. Then from [HLZ; 1.5] one
concludes:

 \Prop{1.4}{\sl Let $a$ be any current of degree $k$ on $X$
such that $da=\phi-R$ where $\phi\in \ce^{k+1}(X)$ and $R$ is
integrally flat.  If $da=\phi'-R'$ is a similar decomposition, then
$\phi=\phi'$ and $R=R'$. Furthermore, 
$$
d\phi\ =\ 0\and dR\bigl|_{X-\bx}\ =\ 0
$$ 
and $\phi$ has integral periods on cycles in $X$.  In the  case
that $\phi\in \ce^{k+1}\xbx$ and supp$(R)\subset X-\bx$,
one has that $dR=0$ and $\phi$ has integral periods on all relative
cycles in $\xbx$.}
\medskip

Set 
$$\aligned
\cz^{\ell}_0(X) &=\{\phi\in \ce^{\ell}(X): d\phi =0 \text{ and $\phi$
has integral periods}\}  \\
\cz^{\ell}_0\xbx &=  \{\phi\in \ce^{\ell}\xbx : d\phi =0 \text{ and
$\phi$ has} \\
&\ \qquad \text{integral periods on relative cycles in $\xbx$ }\} \\
Z_{\text{rect}}^{\ell}(X) &= \{R\in R^{\ell}(X): dR\bigl|_{X-\bx}= 0\}
\\ Z_{\text{rect}}^{\ell}\xbx &= \{R\in R^{\ell}_{\cpt}(X-\bx): dR= 0\}
\endaligned
\tag1.5$$

\Cor{1.6}{\sl Taking $d_1a=\phi$ and $d_2a=R$ from the decomposition
(1.3) gives well-defined mappings 
$$
d_1:S^k(X)\ \arr\ \cz^{k+1}_0(X), \qquad
d_2:S^k(X)\ \arr\ Z_{\text{rect}}^{k+1}(X),\qquad\text{and}
$$ 
$$
d_1:S^k\xbx\ \arr\ \cz^{k+1}_0\xbx, \qquad
d_2:S^k\xbx\ \arr\ Z_{\text{rect}}^{k+1}(X-\bx)
$$
}

\Prop{1.7} {\sl There are natural isomorphisms
$$
\Psi:\hH^k(X)@>{\cong}>> \hat{H}^k(X;\,\bbr/\bbz)
\and
\Psi:\hH^k\xbx@>{\cong}>> \hat{H}^k(X,\bx;\,\bbr/\bbz)
$$
induced by integration.
}
\pf  The first is proved in [HLZ].  The argument for the second
is exactly the same. \qed 

\Remark{1.8}In [HLZ] we showed that there are many different (but
equivalent) deRham-Federer definitions of differential characters on a
manifold $Y$.  Each of these different presentations has   obvious
analogues for $\hH^*(X)$ and $\hH^*\xbx$. The proof  of the equivlance
of these definitions closely follows the arguments in [HLZ, \S 2]
and will not be given here. However, this flexibility in definitions
is important in our treatment of the $*$-product.

To illustrate the point we give one
example. Recall that a current $R$ on $X$ is called {\bf integrally
flat} if  $R=S+dT$ where $S$ and $T$ are rectifiable. Denote by
${\cd'}^k(X)\equiv \{\ce^{n-k}(X)\}'$  the space of currents of degree
$k$ on $X$. Let $\cs_{\max}^k\xbx$ denote the set of $a\in {\cd'}^k(X)$
such that $a$ is smooth near $\bx$, $a\bigl|_{\bx}=0$, and $da=\phi-R$
where $\phi\in\ce^{k+1}\xbx$ and $R$ is integrally flat.
Let $\ct_{\max}^k\xbx$ denote the subgroup of elements of the form
$db+S$ where  $b$ is smooth near $\bx$, $b\bigl|_{\bx}=0$, and $S$
is integrally flat.  Then the inclusion $\cs^k\xbx\subset
\cs_{\max}^k\xbx$ induces an isomorphism
$$
\hH^k\xbx\ \cong\ \cs_{\max}^k\xbx/\ct_{\max}^k\xbx
$$
\vskip .3in

\subheading{\S 2. The exact sequences}
The fundamental exact sequences established by Cheeger and Simons in
[CS] carry over to the relative case.

\Def{2.1}    A character $\a\in \hH^k\xbx$ is said to be {\bf smooth}
if $\a=\langle a\rangle$ for a smooth  form $a\in\ce^k\xbx$. The set
of smooth characters is denoted $\hH^k_{\infty}\xbx$. There is a
natural isomorphism
$$
\hH^k_{\infty}\xbx\ \cong\ \ce^k\xbx/\cz^k_0\xbx
$$

\Prop{2.2}{\sl  The mappings $d_1$ and $d_2$  induce
functorial short exact sequences:
$$
0 @>>> H^k(X,\bx;\, S^1) @>{j_1}>>\hH^k\xbx @>{\d_1}>>
\cz^{k+1}_0\xbx@>>>0, \tag{A}
$$
$$
0@>>>\hH^k_{\infty}\xbx @>{j_2}>>\hH^k\xbx @>{\d_2}>>
H^{k+1}(X,\bx;\,\bbz)@>>>0.
\tag{B} $$
}

\pf Note that $\bx$ has a cofinal system of tubular neighborhoods
each of which is diffeomorphic to $\bx\times [0,1)$. We shall use the
following elementary result.
 \Lemma{2.3} {\sl For any $a\in \ce^k(\bx\times [0,1))$ such that 
$da=0$ and $a\bigl|_{\bx}=0$, there exists $b\in   \ce^{k-1}(\bx\times
[0,1))$ such that  $db=a$ and $b\bigl|_{\bx}=0$.}
\pf
Write $a = a_1 + dt\wedge a_2$ where $a_1$ and $a_2$ are  forms
on $X$ whose coefficients depend smoothly on $t\in [0,1)$, or in
other words,  $a_1(t)$, $a_2(t)$ are smooth curves in $\ce^k(X)$ and
$\ce^{k-1}(X)$ repsectively with $a_1(0)=0$. Now $da =
d_xa_1+dt\wedge \frac{\partial a_1}{\partial t} - dt\wedge d_x a_2 =0$.
We conclude that $d_xa_1=0$ and $d_xa_2= \frac{\partial a_1}{\partial
t}$. Since $a_1(0)=0$ we have
$$
a_1(t)\ =\ \int_0^t \frac{\partial a_1}{\partial t}(s)\,ds\ =\
\int_0^td_xa_2(s)\,ds\ =\ d_x\int_0^ta_2(s)\,ds 
$$
Set $b\equiv\int_0^ta_2(s)\,ds$, and note that: $b\bigl|_{\bx}=0$,
$d_xb=a_1$ and $\frac{\partial b}{\partial t}= a_2$.  Hence, $a =
db$.\qed
\medskip

We shall also need the following result.
On any manifold $Y$ let 
$$
\cf^k(Y)\equiv
\ce^k_{\Lloc}(Y)+d\ce^{k-1}_{\Lloc}(Y)
$$
denote {\bf flat currents} and
$\cf^k_{\cpt}(Y)$ those with compact support. 
Note that $d\cf^k(Y)=d\ce^k_{\Lloc}(Y)$. This definition of 
$\cf^k_{\cpt}(Y)$ arises naturally in sheaf theory.  However,
the following equivalent definition will also be useful here.

\Lemma {2.4}{ 
$\cf^k_{\cpt}(Y)\ =\
\ce^k_{\Lloc,\cpt}(Y)+d\ce^{k-1}_{\Lloc,\cpt}(Y)$ and so
$
 d\cf^k_{\cpt}(Y)\ =\ d\ce^k_{\Lloc,\cpt}(Y)
$
}
\pf
Fix $f\in\cf^k_{\cpt}(Y)$ and write  $f=a+db$
where $a\in \ce^k_{\Lloc}(Y)$ and $b\in \ce^{k-1}_{\Lloc}(Y)$.
Let $K=\supp(f)$, and note that in $N\equiv Y-K$ we have that $a=-db$.
By standard de Rham theory there exists an $\Lloc$-form $b_0$ on $N$
such that $a_{\infty}\equiv a+db_0$ is smooth on $N$. Furthermore since
$a_{\infty}$ is weakly exact on $N$ there exists a smooth form 
$b_{\infty}$ with $a_{\infty}=-db_{\infty}$ on $N$.  
Choose $\eta\in C_0^{\infty}(Y)$ with $\eta\equiv1$
in a neighbothood of $K$, let $\chi=1-\eta$ and set 
$\wt a = a +d(\chi b_0 + \chi b_{\infty})$ and 
$\wt b =b- \chi b_0 - \chi b_{\infty}$ with $\chi$ as above. 
Then $f=\wt a+d\wt b$ and $\wt a$ has compact support in $Y$.

Observe now that $f-\wt a$ is $d$-closed and has compact support 
in  $Y$. Since $H^*(\ce^*_{\cpt}(Y))\cong H^*(\cf^*_{\cpt}(Y))$
we conclude that there exist a smooth form $\omega$ and a flat form
$g$, both having compact support on $Y$ such that $f-\wt a=\omega +dg$.
Now by the paragraph above we can write $g= b+de$ where $b$ is $\Lloc$
 with compact support.  Hence $f=\wt a +\omega +db$.\qed
\medskip

We first prove the surjectivity of $\d_1$.
 Fix $\phi\in\cz^{k+1}_0\xbx$.  Then by Lemma 2.3 there is a
neighborhood $N  \cong \bx\times [0,1)$ of $\bx$ and a form
$A\in\ce^k(N)$ with $dA=\phi$ and $A \bigl|_{\bx}=0$. Choose $\chi\in
C_0^{\infty}(N)$ with $\chi\equiv1$ in a neighborhood of $\bx$, and set
$\phi_0=\phi-d(\chi A)$. Now supp$(\phi_0)\subset\subset X-\bx$ and
$\phi_0$ has integral periods, so there exists a cycle $R\in 
Z_{\text{rect}}^{\ell}\xbx$ with $[\phi_0-R]=0$ in 
$H^{*}_{\cpt}(X-\bx;\bbr)$.  By Lemma 2.4 there are $\Lloc$-forms
$a,b$ with compact support in $X-\bx$ such that $d(a+db)=da=\phi_0-R$.
Then  $d_1(\chi A + a) =\phi$ and surjectivity is proved.

We now construct the map $j_1$.
Recall that (cf. HLZ; \S1])
$$
H^k(X,\bx;\,S^1)\cong H^k_{\cpt}(X-\bx;\,S^1)\cong
\frac{\{f\in \cf^k_{\cpt}(X-\bx): df\in \cR^{k+1}_{\cpt}(X-\bx)\}}
{d \cf^{k-1}_{\cpt}(X-\bx)+\cR^{k}_{\cpt}(X-\bx)}
$$
Choose $f\in \cf^k_{\cpt}(X-\bx)$ with $df=R\in
\cR^{k+1}_{\cpt}(X-\bx)$, and write $f=a+db$ where $a$ and $b$ are
$\Lloc$-forms with compact support in $X-\bx$ (cf. Lemma 2.4).  Then
$ a \in \cs^k\xbx$ and we set
$j_1(f)\equiv\langle a\rangle \in\hH^k\xbx$. Note that if $f=a'+db'$
is another such decomposition, then  $a-a' = d(c'-c)$ and   
$\langle a\rangle =\langle a'\rangle$. Clearly $j_1=0$ on
$d \cf^{k-1}_{\cpt}(X-\bx)+\cR^{k}_{\cpt}(X-\bx)=d
\ce^{k-1}_{\cpt}(X-\bx)+\cR^{k}_{\cpt}(X-\bx)$, and so it descends to
the quotient $H^k(X,\bx;\,S^1)$.

To see that $j_1$ is injective, let $f=a+db$  as above and suppose 
$ a = dc +S\in \ct^k\xbx$ where $c$ is smooth and zero on $\bx$. By
Lemma 2.3 there exists an $\Lloc$-form $e$, smooth near $\bx$, such
that $c_0 = c-de \equiv 0$ near $\bx$.  Then $a = dc_0+S \equiv 0$
in $H^k(X,\bx;\,S^1)$.
 
We now prove the exactness of (A) in the middle.   
Suppose $a\in
\cs^k\xbx$ and $\d_1(\langle a \rangle)=0$.  Then $da = -R
\in  \cR^{k+1}_{\cpt}(X-\bx)$.  Thus, in a neighborhood $N$ of 
$\bx$ we have that $a$ is smooth, $da=0$ and $a\bigl|_{\bx}=0$.
By Lemma 2.3 there exists $b\in\ce^{k-1}(N)$ with $db=a$ and 
$b\bigl|_{\bx}=0$. Then $\wt a = a -d(\chi b)$, with $\chi$ as
above,  is equivalent to $a$ in
$\hH^k\xbx$. Since $\wt a$ has compact support in $X-\bx$ and $d\wt
a=-R$, we see that $ \langle\wt a\rangle$ lies in the
image of $j_1$. 

We now prove the surjectivity of $\d_2$.  Fix $u\in
H^{k+1}(X,\bx;\,\bbz)$ and choose a cycle $R\in u$. Then there is a
smooth form $\phi \in \cz^{k+1}_0\xbx$ such that $\phi-R=df$ for
$f\in\cf^k_{\cpt}(X-\bx)$. By Lemma 2.4 $f=a+db$ where $a$ is $\Lloc$
with compact support in $X-\bx$. Then $a\in\cs^k\xbx$ and $\d_2(\langle
a\rangle )= u$.

Now consider an element $a\in \cs^k\xbx$ with $\d_2(\langle
a\rangle )= 0$. Then $da=\phi-R$ where $\phi$ is smooth and $R=dS$ for
some $S\in\cR^k_{\cpt}(X-\bx)$. Then $\wt a=a-S \equiv a$ in
$\hH^k\xbx$ and $d\wt a=0$ on $X$. Since $\wt a$ is smooth near
$\bx$, standard de Rham theory shows that there is an $\Lloc$-form $b$ 
with compact support in $X-\bx$ such that $\wt a-db$ is smooth.
Hence, $\langle a\rangle=\langle \wt
a\rangle\in\hH^k_{\infty}\xbx$.\qed

\medskip

Note that 
$$
\ker(\d_1)\cap\ker(\d_2)\ \cong\ \frac{H^k(X,\bx;\,\bbr)}
{H^k(X,\bx;\,\bbz)_{\text{free}}}\ \cong\
\frac{H^k_{\cpt}(X-\bx;\,\bbr)} {H^k_{\cpt}(X-\bx;\,\bbz)}
\tag2.5
$$

\vskip.3in

\subheading{\S  3. The  star product}  In this section we  prove the
following.

\Theorem{3.1} {\sl There are functorial bilinear mappings
$$\aligned
\hH^k\xbx\times \hH^{\ell}\xbx\
&@>{*}>>\ \hH^{k+\ell+1}\xbx\qquad\text{and} 
\\
\hH^k\xbx\times \hH^{\ell}(X)\ &@>{*}>>\ \hH^{k+\ell+1}\xbx
\endaligned
$$
 which make $\hH^k\xbx$  a graded commutative ring
and $\hH^*(X)$  a graded $\hH^k\xbx$-module. With this structure
the maps $\d_1$, $\d_2$ are ring and module homomorphisms. 
}

\pf  Fix $\a\in \hH^k\xbx$ and $\b\in\hH^{\ell}(X)$.  Then from
[HLZ] we know that there exist sparks $a\in\a$ and $b\in\b$ with
$$
da\ =\ \phi-R \and db\ = \ \psi-S
$$
with $\phi\in\cz^{k+1}_0\xbx, \psi\in\cz^{\ell+1}_0(X), R\in
Z_{\text{rect}}^{k+1}\xbx$ and $S\in Z_{\text{rect}}^{\ell+1}(X)$, so
that the wedge-intersection products $R\wedge b$ and $R\wedge S$ are
well defined. Furthermore, if $\supp S\subset\subset X-\bx$ we can
also assume that $a\wedge S$ is well defined.  We then define
$$
a*b\ =\ a\wedge \psi +(-1)^{k+1}R\wedge b,
\tag3.2
$$
and if  $S\in Z_{\text{rect}}^{\ell+1}\xbx$ or if
$a\in\ce^k_{\cpt}(X-\bx)$, we can also define 
$$
a\wt*b\ =\ a\wedge S +(-1)^{k+1}\phi\wedge b.
\tag3.3
$$
Since $a$ is smooth near $\bx$ and $a\bigl|_{\bx}=0$, $a*b$
also has these properties (as well as $a\wt *b$ when
it is defined).  Note that 
$$
d(a*b) \ =\  d(a\wt*b)\ =\ \phi\wedge\psi - R\wedge S
\tag3.4
$$
The arguments from [HLZ] easily adapt to show that
$\langle a * b\rangle$ depends only on $\langle a \rangle$ and
$\langle b\rangle$, and that $\langle a * b\rangle=\langle a\wt *
b\rangle$ (when it is defined). Associativity, commutativity,
etc. are straightforward. Equation (3.4) establishes the homomorphism
propertes of $\d_1$ and $\d_2$. \qed

\vskip.3in 
\subheading{\S 4.  Smooth Pontrjagin Duals}
The exact sequences of Proposition 2.2 show that $\hH^*\xbx$ has a
natural topology making it a topological group  (in fact a topological
ring) for which $\d_1$ and $d_2$ are continuous homomorphisms.
Essentially it is a product of the standard $C^\infty$-topology on
forms with the standard topology on the torus 
$H^k(X,\bx;\,\bbr)/H^k(X,\bx;\,\bbz)$. It can also be defined as the
quotient of the topology induced on sparks by the embedding
$\cs^k\xbx\subset \cf^k(X)\times \ce^{k+1}\xbx\times
\cR^{k+1}_{\cpt}(X-\bx)$ sending $a\mapsto (a,d_1a,d_2a)$.
(Similar remarks apply to $\hH^*(X)$.)

It is natual to consider the dual to $\hH^*\xbx$ in the sense of
Pontrjagin. For an abelian topological group $A$ we denote by
$A\du\equiv\Hom_{\text{cont}}(A, S^1)$   the group of continuous
homomorphisms     $h:A\to S^1$.  Then 2.2(B) gives a dual sequence 
$$
0@>>>   H^{k+1}(X,\bx;\,\bbz)\du
 @>{}>>\hH^k\xbx\du @>{\r}>>  \hH^k_{\infty}\xbx\du
@>>>0.
\tag{4.1} $$
where $\r$ is the restriction mapping. 

\Def{4.2}  An element $f\in 
\hH^k_{\infty}\xbx\du$ is called {\bf smooth} if there exists a form 
$\omega\in\cz^{n-k}_0(X)$ such that 
$$
f(\a) \ \equiv\  \int_X\, a\wedge \omega\ \ (\text{mod} \ \  \bbz)
$$
for $a\in\a\in\hH^k_{\infty}\xbx$.
An element $f\in\hH^k\xbx\du$ is called {\bf smooth}  if $\r(f)$ is
smooth. The set of these is called the
{\bf smooth Pontrjagin dual} of $\hH^k\xbx$ and is denoted by 
$\hH^k\xbx\sdu = \Hom_{\infty}(\hH^k\xbx,S^1)$

\Prop{4.3}{\sl The  smooth Pontrjagin dual   $\hH^k\xbx\sdu$ is dense
in  $\hH^k\xbx\du$.}
\pf
 Applying $\d_1$ to 
$\hH^k_{\infty}\xbx$ gives an exact sequence
$$
0@>>> T @>>> \hH^k_{\infty}\xbx @>>> d\ce^k\xbx @>>> 0
$$
where $T=H^k(X,\bx;\,\bbr)/H^k_{\text{free}}(X,\bx;\,\bbz)$, 
with dual sequence $$
0@>>> d\ce^k\xbx\du @>>> \hH^k_{\infty}\xbx\du  @>>> T\du  @>>> 0
\tag4.4$$
Observe   that $T\du=H^k_{\text{free}}(X,\bx;\,\bbz)\cong
H^{n-k}_{\text{free}}(X;\,\bbz)$,  and that
$d\ce^k\xbx\du=\{d\ce^k\xbx\}'$ (the topological vector space dual)
which is exactly the space of currents of degree  $n-k-1$ on $X$
restricted to the closed subspace $d\ce^k\xbx$. This gives a
commutative diagram 
$$
 \CD 0  @>>> d\ce^{n-k-1}(X)@>{}>>
\cz^{n-k}_0(X)@>>>H^{n-k}_{\text{free}}(X;\,\bbz)@>>>0\\ 
@. @VVV@VVV@VV{\cong}V\\
0  @>>> d{\cd'}^{n-k-1}(X)@>{}>>  \hH^k_{\infty}\xbx\du @>>> T\du@>>>0\\
\endCD
$$
with exact rows. Since $\ce^{n-k-1}(X)$ is dense in ${\cd'}^{n-k-1}(X)$,
the result follows. \qed

\medskip

There is a parallel story for $\hH^*(X)$. The analogue of 2.2(B)
gives an exact sequence
$$
0@>>>   H^{k+1}(X;\,\bbz)\du
 @>{}>>\hH^k(X)\du @>{\r}>>  \hH^k_{\infty}(X)\du.
@>>>0.
\tag{4.5} 
$$

\Def{4.6}  An element $f\in 
\hH^k_{\infty}(X)\du$ is called {\bf smooth} if there exists a form 
$\omega\in\cz^{n-k}_0\xbx$ such that 
$$
f(\a) \ \equiv\  \int_X\, a\wedge \omega\ \ (\text{mod} \ \  \bbz)
$$
for $a\in\a\in\hH^k_{\infty}(X)$ 
An element $f\in\hH^k(X)\du$ is called {\bf smooth}  if $\r(f)$ is
smooth. The set of these is called the
{\bf smooth Pontrjagin dual} of $\hH^k(X)$ and is denoted 
$\hH^k(X)\sdu = \Hom_{\infty}(\hH^k(X),S^1)$

\Prop{4.7}{\sl The  smooth Pontrjagin dual  $\hH^k(X)\sdu$ is dense in 
$\hH^k(X)\du$.}
\pf
 Applying $\d_1$ to 
$\hH^k_{\infty}(X)$ gives an exact sequence
$$
0@>>> T @>>> \hH^k_{\infty}(X) @>>> d\ce^k(X) @>>> 0,
$$
where $T=H^k(X;\,\bbr)/H^k(X;\,\bbz)$, 
with dual sequence $$
0@>>> d\ce^k(X)\du @>>> \hH^k_{\infty}(X)\du  @>>> T\du  @>>> 0.
\tag4.8$$
Observe now that $T\du=H^k(X;\,\bbz)\cong H^{n-k}(X,\bx;\,\bbz)$, and
$d\ce^k(X)\du=\{d\ce^k(X)\}'$ is  the space of currents of degree
$n-k-1$ on $X$ restricted to the closed subspace $d\ce^k(X)$. This
gives a commutative diagram:
$$
\CD
 \ce^{n-k-1}\xbx@>{d}>>
\cz^{n-k}_0\xbx@>>>H^{n-k}(X,\bx;\,\bbz)@>>>0\\ @VVV@VVV@VV{\cong}V\\
{\cd'}^{n-k-1}(X)@>{d}>>  \hH^k_{\infty}(X)\du @>>> T\du@>>>0\\
\endCD
$$
with exact rows.  Since $\ce^{n-k-1}\xbx$ is dense in ${\cd'}^{n-k-1}(X)$,
the result follows. \qed

\vfill\eject

\vskip.3in 
\subheading{\S 5.  Lefschetz-Pontrjagin Duality} This brings us to the
main result of the paper.

\Theorem{5.1}  {\sl  Let $X$ be a compact, oriented $n$-manifold with
boundary  $\bx$.  Then the biadditive mapping
$$
\hH^k\xbx\times\hH^{n-k-1}(X)\ \arr\ S^1
$$
given by
$$
(\a,\b) \ \ \ \ \mapsto\ \ \ \ (\a*\b)\,[X]
$$
induces isomorphisms
$$
\cd:\hH^k\xbx@>{\cong}>>\hH^{n-k-1}(X)\sdu
\qquad\text{and} 
$$ 
$$
\cd':\hH^k(X)@>{\cong}>> \hH^{n-k-1}\xbx\sdu\ \ \ \
$$

}
\smallskip
\pf
Fix $\a\in\hH^k\xbx$ and suppose $(\a*\b)[X]=0$ for all
$\b\in\hH^{n-k-1}(X)$.  We shall show that $\a=0$. Choose a spark
$a\in \a$ and write $da=\phi-R$ as in 1.4.  Then for all smooth forms 
$b\in \ce^{n-k-1}(X)$ we have by (3.3) that
$$
\a*\<b\>\,[X]\ =\ (-1)^{k+1}\int_X\phi\wedge b\ \equiv\ 0\mod\bbz
$$
since $d_2b=0$. It follows that $\phi=0$.  

Hence, $da=-R\in
\cR^{k+1}_{\cpt}(X-\bx)$ is a cycle with
 $[R]\in H^{k+1}_{\cpt}(X-\bx;\,\bbz)_{\text{tor}}\cong
H_{n-k-1}(X-\bx;\,\bbz)_{\text{tor}}$.   Choose any
$u\in H^{n-k}(X;\,\bbz)_{\text{tor}}\cong
H_{k}(X,\bx;\,\bbz)_{\text{tor}}$, and choose a relative 
cycle $S\in u$. Let $m$ be the order of
$u$.  Then there is a $(k+1)$-chain $T$ on $X$ with $dT=mS$ rel $\bx$. 
Set $b=-\frac 1 m T$ and consider $b$ as a spark of degree $n-k-1$ on
$X$ with $db=-S$.  Now we may assume $S$ and $T$ to have been chosen so
that  $\supp(S)\cap \supp(R)=\emptyset$ and  $T$ meets $R$
properly. Then
$$
\aligned
0\ =\ \a*\<b\>\,[X]\ &\equiv\ (-1)^{k+1}R\wedge b\, [X]\mod \bbz \\
&\equiv\ (-1)^{k+1}\frac 1 m R\wedge T \,[X]\mod \bbz \\
&\equiv\ (-1)^{k+1}\text{Lk}([R], [S]) \mod \bbz\\
&\equiv\ (-1)^{k+1}\text{Lk}(\d_2\a, u) \mod \bbz
\endaligned
$$
where Lk denotes the de Rham-Seifert linking between  the groups
$H_{n-k-1}(X-\bx;\,\bbz)_{\text{tor}}$ and 
$H_{k}(X,\bx;\,\bbz)_{\text{tor}}$.   By the non-degeneracy of this
pairing we conclude that $\d_2\a=0$.

Therefore $\a \in \ker(\d_1)\cap\ker(\d_2)$ can be represented by a
smooth $d$-closed form $a\in \ce^k\xbx$. In fact by Lemma 2.3 we may
choose $a$ to have compact support in $X-\bx$.
Now for any cycle $S\in Z_{\text{rect}}^{n-k}(X)$,
i.e., any $k$-dimensional rectifiable current $S\in\cR_k(X)$
with $dS\in \cR_{k-1}(\bx)$, we can find $\psi\in\ce^{n-k}(X)$
and $b\in \ce^{n-k-1}_{\Lloc}(X)$ with $db=\psi-S$. Then by (3.3) we
have that 
$$
\aligned
0\ =\ \a*\<b\>\,[X]\ &\equiv\ a\wedge S\ [X]\mod \bbz \\
&\equiv\ \int_S\,a \mod \bbz.
\endaligned
$$
Hence, $a$ represents the zero class in
$\Hom(H_k(X,\bx;\,\bbz),\,\bbr)/ \Hom(H_k(X,\bx;\,\bbz),\,\bbz)\cong
{H^k(X,\bx;\,\bbr)}/ {H^k(X,\bx;\,\bbz)_{\text{free}}}$, and by (2.2)
and (2.5)  we conclude that $\a=0$.
Thus 
\eject

\noindent
the map $\cd$ is injective.

To see that  $\cd$ is surjective consider the commutative
diagram with exact rows:
$$
\CD
0@>>> H^k(X,\bx;\,S^1) @>{j_1}>> \hH^k\xbx @>{\d_1}>>\cz_0^{k+1}\xbx
@>>> 0\\ @.  @V{\cong}VV@VV{\cd}V@VV{\cd_0}V\\
0@>>> \Hom(H^{n-k}(X;\,\bbz),S^1) @>>> \hH^{n-k-1}(X)\du
@>{\r}>>\hH^{n-k-1}_{\infty}(X)\du @>>> 0 \endCD
$$
where the top row is 2.2(A) and the bottom row is the dual of 2.2(B).
By definition $\cd_0$ is onto the smooth elements in 
$\hH^{n-k-1}(X)\du$ and therefore the map $\cd$ is surjective.
\medskip

The proof that $\cd'$ is an isomorphism is parallel.
Fix $\b\in\hH^{n-k-1}(X)$ and suppose $(\a*\b)[X]=0$ for all
$\a\in\hH^k\xbx$.  We shall show that $\b=0$. Choose a spark
$b\in \b$ and write $db=\psi-S$ as in 1.4.  Then for all smooth forms 
$a\in \ce^k\xbx$ we have by (3.3) that
$$
\<a\>*\b\,[X]\ =\  \int_X\,a\wedge \psi\ \equiv\ 0\mod\bbz
$$
since $d_2a=0$. It follows that $\psi=0$.  

Hence, $db=-S\in
\cR^{n-k}(X)$ is a relative cycle with torsion homology class
 $[S]\in H^{n-k}(X;\,\bbz)_{\text{tor}}\cong
H_{k}(X,\bx;\,\bbz)_{\text{tor}}$.   Choose 
$u\in H^{k+1}(X, \bx;\,\bbz)_{\text{tor}}\cong
H_{n-k-1}(X;\,\bbz)_{\text{tor}}$, and choose a  
cycle $R\in u$ with support in $X-\bx$. Let $m$ be the order of
$u$.  Then there is a $(n-k-1)$-chain $T$ in $X-\bx$ with $dT=mR$.  Set
$a=-\frac 1 m T$ and consider $a$ as a spark of degree $k$ on $X$
with $da=-R$.  Now we may assume $R$ and $T$ to have been chosen so
that  $\supp(R)\cap \supp(S)=\emptyset$ and  $T$ meets $S$ properly.
Then 
$$
\aligned
0\ =\ \<a\>*\b\,[X]\ &\equiv\ (-1)^{k+1}a\wedge S\, [X]\mod \bbz \\
&\equiv\ (-1)^{k+1}\frac 1 m T\wedge S \,[X]\mod \bbz \\
&\equiv\ (-1)^{k+1}\text{Lk}([R], [S]) \mod \bbz\\
&\equiv\ (-1)^{k+1}\text{Lk}(u, \d_2\b) \mod \bbz
\endaligned
$$
where Lk denotes the de Rham-Seifert linking as before.   We conclude
that $\d_2\a=0$.

Therefore $\b \in \ker(\d_1)\cap\ker(\d_2)$ can be represented by a
smooth $d$-closed form $b\in \ce^{n-k-1}(X)$.
Now for any cycle $R\in Z_{\text{rect}}^{k+1}(X,\bx)$,
i.e., any $k$-dimensional rectifiable current $R\in\cR_{n-k-1}(X-\bx)$
with $dR=0$, we can find $\phi\in\ce^{k+1}\xbx$
and $a\in \ce^{k}_{\Lloc}\xbx$ with $da=\phi-R$. Then by (3.2) we
have that 
$$
\aligned
0\ =\ \<a\>*\b\,[X]\ &\equiv\ (-1)^{k+1}R\wedge b\ [X]\mod \bbz \\
&\equiv\ (-1)^{n(k+1)}\int_R\,b \mod \bbz.
\endaligned
$$
Hence, $b$ represents the zero class in
$\Hom(H_{n-k-1}(X;\,\bbz),\,\bbr)/
\Hom(H_{n-k-1}(X;\,\bbz),\,\bbz)\cong {H^{n-k-1}(X;\,\bbr)}/
{H^{n-k-1}(X;\,\bbz)_{\text{free}}}$, and by (2.2) and (2.5)   we
conclude that $\b=0$. Thus 
the map $\cd'$ is injective.

The surjectivity of  $\cd'$ follows as before from the commutative
diagram with exact rows:
$$
\CD
0@>>> H^{n-k-1}(X;\,S^1) @>{j_1}>> \hH^{n-k-1}(X)
@>{\d_1}>>\cz_0^{n-k}(X) @>>> 0\\ @.  @V{\cong}VV@VV{\cd}V@VV{\cd_0}V\\
0@>>> \Hom(H^{k+1}(X,\bx;\,\bbz),S^1) @>>> \hH^{k}\xbx\du
@>{\r}>>\hH^{k}_{\infty}\xbx\du @>>> 0. \endCD
$$
This completes the proof.\qed

\vskip .3in

\def\co{\partial}

\subheading{\S 6.  Coboundary maps} It is natural to ask if there is a
coboundary mapping $\co$ with the property that the sequence
$$
\dots@>{}>> \hH^{k-1}(\bx)
 @>{\co}>> \hH^k\xbx @>{j}>> \hH^k(X) @>{\r}>> \hH^k(\bx) 
@>{\co}>> \hH^{k+1}\xbx @>>>\dots
\tag6.1
$$
is exact. The differential-form-component of characters makes this
impossible.  However, there do exist natural coboundary maps 
$\co$ with the following properties:
\roster
\item Under $\d_2$ the sequence (6.1) becomes the standard 
long exact sequence in integral cohomology.
\item Under  $\d_1$ the sequence (6.1) becomes a sequence of smooth
$d$-closed forms which induces the standard long exact sequence in real
cohomology. \endroster

\def\onen{{\Bbb I}_N}
Recall that the
definitions of  Thom maps and Gysin maps for differential characters
depend essentially on a choice of ``normal geometry''. This will also
be true for our coboundary maps.  Fix a tubular neighborhood $N_0$ of 
$\bx$ in $X$ and an identification  $N_0\cong \bx\times
[0,2)$, and let  $\pi:N_0\to \bx$ be the projection.
Set $N=\bx\times [0,1)\subset N_0$ and let $\onen$ be the
characteristic function of this subset. Let $\chi$ be
a smooth approximation to $\onen$; specifically choose
$\chi(x,t)=\chi(t)$ where $\chi\equiv 1$ near 0  and $\chi(t)=0$
for $t\geq 1$.  Then set 
$$
\lambda \ \equiv\ \chi -\onen \in \hH^0(X)
$$
Note that $d\lambda=d\chi-[\partial N]$ has compact support in
$X-\bx$.

\Def{6.2}  We define the coboundary map $\co=\co_{\lambda}:\hH^k(\bx) 
@>{}>> \hH^{k+1}\xbx$ by
$$
\co (a) \ =\  \pi^*a * \lambda.
$$

\medskip

Verification of  (1) and (2) above is straightforward, and
the details are omitted.

\vskip .3in

\subheading{\S 7. Sequences and duality}  At the level of cohomology
the long exact sequences for the pair $\xbx$ are related by the
duality mappings.  There is an analogous diagram for differential
characters:
$$\CD
\hH^k\xbx @>{j}>> \hH^k(X) @>{\r}>> \hH^k(\bx) 
@>{\co}>> \hH^{k+1}\xbx    \\
@V{\cd}VV@V{\cd}VV@V{\cd}VV@V{\cd}VV   \\
  \hH^{n-k-1}(X)\du @>{j^*}>> \hH^{n-k-1}\xbx\du @>{\co^*}>>
\hH^{n-k-2}(\bx)\du  @>{\r^*}>> \hH^{n-k-2}(X)\du 
\endCD
$$
and it is natural to ask whether this diagram commutes (up to sign).
The square on the left is evidently commutative. 
The other two squares
commute up to an error  term which we now analyse.

We begin with the square on the right.  Fix $\a\in\hH^k(\bx)$
and $\b\in \hH^{n-k-2}(X)$ and choose $\Lloc$-sparks $a_0\in\a$ and
$b\in\b$ with $da_0=\phi_0-R_0$ and $db=\psi-S$ as usual.
Let $a=\pi^*a_0$, $\phi=\pi^*\phi_0$ and $R=\pi^*R_0$ denote the 
pull-backs to the collar neighborhood of $\bx$ via the projection
$\pi:N_0\to \bx$ defined in \S 6.
Then
$$
\{(\cd\circ\co)(\a)\}(\b) \ = \ (\pi^*a * b *\lambda)[X]
\ =\ \{(a * b)\wedge d\chi +(-1)^nd_2(a * b)\lambda\}[X].
\tag7.1$$
Now we may assume that $S\bigl|_{N_0} = 
\pi^*S_0$ for some $S_0\in \cR^{k+1}(\bx)$, and we may further assume
that $\supp(R_0)\cap\supp(S_0)=\emptyset$  because
$\dim(R_0)+\dim(S_0)=n-2$. Hence, $d_2(a*b) = \pi^*R_0\wedge\pi^*S_0
= \pi^*(R_0\wedge S_0)=0$, and from (7.1) we see that
$$
\aligned
(-1)^{n-1}\{(\cd\circ\co)(\a)\}(\b) \ &= \ (-1)^{n-1}(a * b)\wedge
d\chi[X]\\ 
&=\ (a * b)[\bx] - \chi d(a*b) [X] \\
&=\ \{(\r^*\circ \cd)(\a)\}(\b) - \chi d(a*b) [X].
\endaligned
$$
Now $d(a*b)=\phi\wedge\psi-R\wedge S = \phi\wedge \psi$ and we can
write $\psi = \psi_1+dt\wedge \psi_2$ as in the proof of Lemma 2.3.
Since $\phi=\pi^*\phi_0$ we see that $\phi\wedge\psi_1=0$ and we
conclude that
$$\aligned
\{(\r^*\circ \cd)(\a)\}(\b)+(-1)^{n}\{(\cd\circ\co)(\a)\}(\b)
\ &=\ \int_N\phi\wedge\chi dt \wedge  \psi_2  \\
  &=\ \int_{\bx}\phi\wedge  \int_0^1 \chi(t)\,dt \wedge \psi_2 \\
&= \ \int_{\bx}\phi\wedge\pi_*\left\{\chi(t)\,dt \wedge \psi_2 \right\}
\equiv E(\lambda).
\endaligned
\tag7.2
$$
Thus for example we see that\ \ 
$(\r^*\circ \cd)(\a)=(-1)^{n-1}(\cd\circ\co)(\a)$ on all $\b$ which
are $\pi^*$-pull backs in $N$.  Furthermore, we can consider the 
family of sparks $\lambda_{\e}\equiv r_{\e}^*\lambda$ where
$r_{\e}:\bx\times[0,\e)\to\bx[0,1)$ is given by $r_{\e}(x,t)
=(x,t/\e)$.  From (7.2) we see that
$$\lim_{\e\to0}E(\lambda_{\e})=0.$$

A similar analysis applies to middle square in the diagram and we
have the following. 
\Prop{7.3}{\sl The duality 
diagram above commutes in the limit as $\e\to0$. }
\medskip

This is the best one can expect.  The ``commutators'' in this diagram
do not lie in the smooth dual.  Of course by Propositions 4.3 and
4.7 they do lie in its closure.

Here is an explicit example of this non-commutativity. Let
$X=S^2\times D^3$ be the product of the 2-sphere and the 3-disk.
Choose sparks $\a\in \cs^1(S^2)$ and $b\in\cs^2(D^3)$ with
$da=\omega - [x_0]$ and $db=\Omega-[0]$ for some $x_0\in S^2$, where
$\omega$ and $\Omega$ are unit volume forms on $S^2$ and $D^3$
respectively.  Direct calculation shows that 
$$
(a*b)[\bx]\ =\ 1\qquad\text{but}\qquad
(a*\lambda*b)[X]\ =\ \int_{D^3}(1-\chi)\Omega \ <\ 1.
$$

 

\centerline{\bf References}

\vskip .2in

\nobreak

\ref\key C \by J. Cheeger\paper Multiplication of 
Differential
Characters \jour Instituto Nazionale di Alta Mathematica, 
Symposia Mathematica \vol XI \yr 1973 \pages 441--445
\endref

\ref\key  CS \by J. Cheeger and J. Simons\paper
Differential Characters and Geometric Invariants
\jour Lect. Notes in Math.\vol 1167\publ Springer--Verlag
\publaddr New York\yr 1985 \pages 50--80\endref

\ref\key  deR \by  G. de Rham \book Vari\'et\'es Diff\'erentiables, formes,
courants, formes harmoniques\publ 
 Hermann\publaddr Paris \yr 1955\endref

\ref\key  F \by   H. Federer\book Geometric Measure 
Theory\publ 
 Springer--Verlag\publaddr New York \yr 1969\endref

\ref \key  GS$_1$  \by H. Gillet and C. Soul\'e \paper  
Arithmetic chow groups and differential characters
 \pages 30-68 \inbook Algebraic K-theory;  Connections with Geometry and
Topology   \publ Jardine and Snaith (eds.), Kluwer Academic Publishers \yr
1989 \endref

\ref \key  GS$_2$  \bysame \paper  
Arithmetic intersection theory
 \pages 94-174 \jour Publ. I.H.E.S.  \vol 72 \yr 1990 \endref

\ref \key  H  \by B. Harris \paper  
Differential characters and the Abel-Jacobi map
 \pages 69-86 \inbook Algebraic K-theory;  Connections with Geometry and
Topology   \publ Jardine and Snaith (eds.), Kluwer Academic Publishers \yr
1989 \endref

\ref \key  HL$_1$  \by F.R. Harvey and H.B. Lawson, Jr. \paper A theory of
 characteristic currents associated with a singular connection
 \pages 1--268 \jour Ast\'erisque \vol 213\yr 1993 \endref

\ref\key  HL$_2$  \bysame \paper Geometric residue
 theorems\jour Amer. J. Math. \vol 117 \yr 1995  \pages 829--
873\endref

\ref\key  HLZ  \by F.R. Harvey, H.B. Lawson, Jr., and John Zweck \paper The de Rham-Federer Theory
of Differential Characters and Character Duality   \jour Amer. J. Math. \vol 125 \yr 2003 \pages  
791-847.  ArXiv:math.DG/0512251
\endref

\ref \key  S  \by J. Simons \paper Characteristic forms and
transgression: characters associated to a connection \jour Stony
Brook preprint,   \yr 1974 \endref

\end